\documentclass[12pt,a4paper]{article}
\usepackage[no-math]{fontspec}
\usepackage[margin=25mm,headheight=27pt]{geometry}
\usepackage{needspace}
\usepackage{enumitem}
\usepackage{fancyhdr}
\setlist{itemsep=3pt,topsep=5pt}
\AddToHook{cmd/subsection/before}{\Needspace{3\baselineskip}}
\usepackage{amsmath,amsfonts,amssymb,amsthm}
\usepackage{mathtools}
\DeclareRobustCommand{\qed}{%
  \ifmmode \mathqed
  \else
    \leavevmode\unskip\nobreak\hbox{}\nobreak\hfill
    \quad\hbox{\qedsymbol}%
  \fi
}
\usepackage{tikz-cd}
\usepackage{booktabs}
\usepackage{xcolor}
\usepackage[colorlinks=true,linkcolor=blue,citecolor=blue,urlcolor=blue]{hyperref}
\allowdisplaybreaks[2]
\numberwithin{equation}{section}
\newtheorem{theorem}{Theorem}[subsection]
\newtheorem{proposition}[theorem]{Proposition}
\newtheorem{lemma}[theorem]{Lemma}
\newtheorem{corollary}[theorem]{Corollary}

\theoremstyle{definition}
\newtheorem{definition}[theorem]{Definition}

\newtheorem{remark}[theorem]{Remark}
\newcommand{\C}{\mathbb C}
\newcommand{\N}{\mathbb N}
\newcommand{\Z}{\mathbb Z}
\newcommand{\one}{\mathbf 1}
\newcommand{\wt}{\operatorname{wt}}
\newcommand{\Span}{\operatorname{span}}
\newcommand{\Sym}{\operatorname{Sym}}
\newcommand{\cA}{\mathcal A}
\newcommand{\cL}{\mathcal L}
\newcommand{\cF}{\mathcal F}
\newcommand{\cT}{\mathcal T}
\newcommand{\m}{\mathfrak m}
\newcommand{\Crel}[2]{C_1^{#1}(#2)}
\title{$C_2$-Cofiniteness of Permutation Orbifolds I:\\Relative $C_1$ Descent}
\author{}
\date{}
\makeatletter
\renewcommand{\maketitle}{%
  \begingroup
  \centering
  {\Large\bfseries\@title\par}
  \medskip
  {\normalsize\renewcommand{\and}{\quad}\@author\par}
  \endgroup
  \medskip
  \thispagestyle{plain}%
}
\makeatother

\hypersetup{pdftitle={C2-Cofiniteness of Permutation Orbifolds I: Relative C1 Descent},pdfauthor={Yuhan Gu and Jiancai Sun}}
\begin{document}
\maketitle
\begin{center}
{Yuhan Gu\\
Department of Mathematics\\
Shanghai University, Shanghai 200444, China\\
Email address: gyh.guyuhan@outlook.com}
\end{center}

\begin{center}
{Jiancai Sun\\
Department of Mathematics and Newtouch Center for Mathematics\\
Shanghai University, Shanghai 200444, China\\
Email address: jcsun@shu.edu.cn}
\end{center}

\begin{abstract}
We establish a relative $C_1$-cofiniteness criterion for
$C_2$-cofiniteness of finite group orbifolds. For permutation
actions on a simple, self contragredient, $C_2$-cofinite vertex operator algebra of
CFT type, it suffices to check a single relation for a derivative
of the conformal vector in the Virasoro vertex operator algebra.
This reduces the problem to a calculation depending only on the
central charge and tensor degree. An exact calculation
then proves cofiniteness for a family of nonsolvable permutation
orbifolds without restricting the central charge.
\end{abstract}


\begingroup
\small
\linespread{0.92}\selectfont
\tableofcontents
\endgroup

\Needspace{7\baselineskip}
\section{Introduction}\label{sec:introduction}

\subsection{Statements of main results}
\label{subsec:introductionresults}

This is the first of two papers on $C_2$-cofiniteness of
permutation orbifolds. We establish a descent criterion for
arbitrary finite automorphism groups, reduce permutation
cofiniteness to a translation relation in the Virasoro vertex operator algebra,
and prove that relation in the fivefold case.
For a vertex operator algebra $U=\bigoplus_{d\ge0}U_d$ of
CFT type and a $U$-module $M$, put $U_{>0}=\bigoplus_{d>0}U_d$.
Following the module convention in \cite{Mfus}, we set
\[
C_1^U(M)=\Span_{\C}\{u_{-1}m:u\in U_{>0},\ m\in M\}.
\]
We do not add $L(-1)M$ to this subspace.

\Needspace{10\baselineskip}
\begin{theorem}\label{thm:descent}
Let $V=(V,Y,\one,\omega)$ be a $C_2$-cofinite vertex operator algebra of CFT type,
let $G$ be a finite group
of vertex operator algebra automorphisms, and set $U=V^G$. If
\[
\dim V/\Crel{U}{V}<\infty,
\]
then $U$ is $C_2$-cofinite.
\end{theorem}

Write $T^n(V)=V^{\otimes n}$ and
$U^{(n)}=U^{(n)}(V)=T^n(V)^{S_n}$.
We omit the argument of $U^{(n)}$ when the underlying vertex operator algebra is clear,
and display it when different vertex operator algebras must be distinguished.
The translation operator is $D=L(-1)$. Put $V_+=\bigoplus_{d>0}V_d$. The graded subspace
$T_n(V)\subset V_+$ is the span of the symmetric operations
$B_n(a_1,\ldots,a_n;b)$ defined in
Definition~\ref{def:symmetricoperations}. These spaces are preserved by
translation and the operators $v_0$.

For a simple vertex operator algebra of CFT type,
Theorem~\ref{prop:permutationsubgroups}, proved in
Section~\ref{sec:permutationsubgroups}, transfers
$C_2$-cofiniteness from $U^{(n)}$ to $T^n(V)^G$ for
every subgroup $G\le S_n$.
We state the following results in this form.

Let $\mathcal V_c$ denote the Virasoro vertex operator algebra of central
charge $c$, with conformal vector $\omega_c=L(-2)\one$ and
translation operator $D=L(-1)$. Its construction is recalled in
Section~\ref{subsec:virasoro-heisenberg}.

We call the following statement the Virasoro criterion.

\begin{theorem}\label{thm:criterion}
Fix $c\in\C$ and an integer $n\ge2$. Suppose that in
the Virasoro vertex operator algebra $\mathcal V_c$ one has
\[
D^r\omega_c\in T_n(\mathcal V_c)
\qquad\text{for some integer }r\ge0.
\]
Then, for every simple, self contragredient, $C_2$-cofinite vertex operator algebra
$V$ of CFT type and central charge $c$, the orbifold $T^j(V)^G$
is $C_2$-cofinite for every $1\le j\le n$ and every subgroup
$G\le S_j$.
\end{theorem}

The hypothesis is a single linear membership test in the
finite dimensional weight-$(r+2)$ subspace of $\mathcal V_c$.
The same calculation applies to every $V$ with central charge $c$;
it requires only the Virasoro relations and the operations $B_n$.

The fivefold case gives the following consequence.

\begin{corollary}
\label{ex:fivefold}
Let $V$ be a simple, self contragredient, $C_2$-cofinite vertex operator algebra
of CFT type. Then
\[
\dim_{\C}\bigl(T^5(V)^{S_5}/C_2(T^5(V)^{S_5})\bigr)<\infty.
\]
The conclusion also holds for $T^j(V)^G$ whenever
$1\le j\le5$ and $G\le S_j$.
\end{corollary}

The Virasoro calculation in Section~\ref{sec:t5}
gives $D^{16}\omega_c\in T_5(\mathcal V_c)$ for every $c\in\C$.
The Virasoro criterion
then proves the assertion, as recorded in
Theorem~\ref{thm:t5-main}.

\subsection{Background and motivation}
\label{subsec:introductionbackground}

The orbifold construction associates to a vertex operator algebra
$V$ and a group of automorphisms $G$ the fixed point algebra
$V^G$. It is a basic way of constructing vertex operator algebras
and studying the effect of symmetry on their representations.
A central finiteness question is whether $C_2$-cofiniteness passes
from $V$ to $V^G$ when $G$ is finite. This paper addresses that
question for permutation actions on tensor powers and develops a
descent criterion that also applies to arbitrary finite actions.

The importance of the $C_2$ condition originates in Zhu's work
on modular invariance \cite{Zhu}. The quotient
$R_V=V/C_2(V)$ is a graded commutative Poisson algebra, so
$C_2$-cofiniteness expresses a finite dimensional constraint on
the vertex algebra itself. It also provides spanning set and
finiteness properties in representation theory
\cite{GN,ABD}. Its significance extends to logarithmic conformal
field theory: Miyamoto proved modular invariance for generalized
trace functions under $C_2$-cofiniteness \cite{Mmod}, while the
triplet algebras furnish $C_2$-cofinite examples with
nonsemisimple representation theory \cite{AMtriplet}.
Accordingly, the results here assume cofiniteness without imposing
rationality, and apply to both rational and logarithmic examples
under the stated hypotheses.

The difficulty with fixed points is already visible in the
$C_2$ quotient. For $U=V^G$, one has
\[
C_2(U)\subseteq U\cap C_2(V),
\]
and averaging over $G$ gives a surjection
\[
R_U\longrightarrow (R_V)^G
\quad\text{with kernel}\quad
(U\cap C_2(V))/C_2(U).
\]
Even when the target is finite dimensional, this does not control
the kernel. Thus ordinary invariant theory on $R_V$ alone does
not settle the orbifold problem. One must account for the
relations that can be formed using invariant vectors.

Permutation actions provide a natural setting for this problem:
the tensor power retains the structure of the original vertex operator algebra,
while its symmetry groups include nonsolvable groups in
arbitrarily large degrees. Barron, Dong, and Mason constructed
permutation twisted sectors from modules for the original vertex operator algebra
\cite{BDM}. This makes permutation orbifolds accessible from the perspective of representation theory and motivates a corresponding
finiteness theorem for their fixed point algebras. Abe proved
$C_2$-cofiniteness for twofold permutation orbifolds and
developed symmetric tensor identities \cite{Abe}. Miyamoto's
theorem for solvable groups establishes $C_2$-cofiniteness of fixed
points of a simple $C_2$-cofinite vertex operator algebra under every finite solvable
group \cite{Msol}. In particular, it covers all permutation
subgroups in tensor degrees at most four. Degree five is the
first degree in which nonsolvable groups occur, namely $A_5$
and $S_5$. Particular nonsolvable actions have also been studied:
Xu's preprint proves $C_2$-cofiniteness and strong rationality
of the icosahedral orbifold $V_{L_2}^{A_5}$ for the rank one
lattice $L_2=\Z\alpha$, $(\alpha,\alpha)=2$ \cite{XuA5}.

Recently, associated varieties of finite group orbifolds were studied by
Riesen \cite{Riesen}. For a simple, strongly finitely generated vertex
algebra with a connected nonnegative integer grading, he identifies
the reduced associated variety of its fixed point algebra under a finite group preserving the grading with the corresponding quotient variety.
As a consequence, $C_2$-cofiniteness descends to finite group orbifolds
of simple $C_2$-cofinite vertex operator algebras of CFT type.

Theorem~\ref{thm:descent} identifies a sufficient condition for
controlling the kernel above: finiteness of the relative module
quotient $V/C_1^{V^G}(V)$. It therefore converts the fixed point
problem into a question about $V$ as a module over its invariant
subalgebra. The proof combines spanning and filtration methods
\cite{GN,LiAb} with the use of $C_1$ bounds in differential
equations and intertwining operators \cite{Huang,Mfus}.
The theorem applies to arbitrary finite $G$, without a
solvability assumption. Its conclusion remains conditional on
relative $C_1$-cofiniteness; the subsequent permutation arguments
provide a method for verifying that condition.

Theorem~\ref{thm:criterion} performs the next reduction in the
Virasoro vertex operator algebra. A single relation
$D^r\omega_c\in T_n(\mathcal V_c)$ establishes cofiniteness
in every tensor degree at most $n$ for every simple,
self contragredient, $C_2$-cofinite vertex operator algebra of CFT type and central
charge $c$. This is a finite dimensional linear membership test
once $r$ is fixed, and its validity depends only on $c$ and $n$.
The canonical homomorphism $\mathcal V_c\to V$ transfers it to
$D^r\omega\in T_n(V)$ by functoriality of the operations $B_n$.
The two quotient descriptions
in Section~\ref{sec:symmetricquotients} explain why such a
small amount of data is enough. They identify a relative
$C_1$ quotient for a position stabilizer with $V/T_n(V)$ and
identify its quotient by translation, after removing the vacuum,
with the strong indecomposables in the preceding tensor degree.
The relation for the conformal vector gives a uniform bound
$D^NV\subseteq T_n(V)$; induction on tensor degree supplies
finite dimensionality modulo $D$. Together these yield the
required relative finiteness. Miyamoto's theorem on
$C_1$-cofinite intertwining images \cite{Mfus} then propagates
the condition from the position stabilizers to the full tensor
power, where Theorem~\ref{thm:descent} applies.

Theorem~\ref{prop:permutationsubgroups} makes these
reductions useful for every $G\le S_n$. Quantum Galois
decomposition \cite{DLM} and module $C_2$-cofiniteness
\cite{ABD} show that cofiniteness of $T^n(V)^{S_n}$ implies
cofiniteness of $T^n(V)^G$. Thus the symmetric group calculation
also settles all permutation subgroups in the same degree,
including groups that are themselves nonsolvable.

Corollary~\ref{ex:fivefold} gives a concrete application at the
first such degree. Section~\ref{sec:t5} proves the Virasoro relation
\[
D^{16}\omega_c\in T_5(\mathcal V_c)\qquad(c\in\C).
\]
The calculation is carried out in weight eighteen with $c$
formal, and functoriality transfers the relation to every $V$
under the hypotheses of Theorem~\ref{thm:criterion}. It follows
that all $T^j(V)^G$ with $j\le5$ and $G\le S_j$ are
$C_2$-cofinite, with no restriction on the central charge.
This resolves the cofiniteness question for permutation orbifolds through
degree five for this class of vertex operator algebras, including the $S_5$ and
$A_5$ cases.
Appendix~\ref{app:t5-verification} gives the exact integer
rank certificate and the complete verification program.

These finiteness results also have a rationality consequence.
Theorem~\ref{thm:external-orbifold-rationality} states that a finite group fixed point algebra of a
strongly rational vertex operator algebra is strongly rational once it is
$C_2$-cofinite. Thus, whenever the original
vertex operator algebra is strongly rational, the cofiniteness conclusions proved
here supply this remaining hypothesis; the fivefold case is
recorded in Corollary~\ref{cor:t5-rationality}. 
A subsequent article will present cofiniteness results for permutation orbifolds in every tensor degree when $V_1\ne0$, and outside
a finite set of possible rational central charges in each
fixed degree when $V_1=0$. It will also treat framed vertex operator algebras
and the moonshine vertex operator algebra in every tensor degree.
The present paper proves the common descent mechanism and
settles the fivefold case for every central charge.

\subsection{Structure of the paper}
\label{subsec:introductionstructure}

Section~\ref{sec:prelim} collects notation, vertex operator identities
and invariant forms, commutative algebra, strong generation and
cofiniteness, rational and algebraic functions, tensor products and
group actions, and vertex operator algebras associated with
the Virasoro and Heisenberg algebras.
Section~\ref{sec:relativedescent} proves the criteria for graded algebras
and the lemmas on rational functions, establishes finite generation and
algebraicity of averaged coefficients, and proves
Theorem~\ref{thm:descent}.
Section~\ref{sec:permutationsubgroups} proves the passage from
symmetric groups to permutation subgroups.
Section~\ref{sec:conformalcriterion} develops the symmetric
operations, identifies the two quotient spaces, proves the module
and translation lemmas, and establishes Theorem~\ref{thm:criterion}.
Section~\ref{sec:t5} derives the Virasoro product identities,
specializes the polynomial relation in weight eighteen, and proves
Corollary~\ref{ex:fivefold}.
Appendix~\ref{app:t5-verification} establishes the finite
spanning certificate and gives the implementation and verification
program.

\Needspace{14\baselineskip}
\section{Preliminaries}\label{sec:prelim}

We collect the conventions and background results used in the descent
argument and in the permutation calculations. The module convention for
relative cofiniteness is kept separate from strong generation and the
Poisson quotient.

Sections~\ref{subsec:conventions}--\ref{subsec:algebraicfunctions}
cover vertex operator identities, commutative algebra, and rational
and algebraic functions. Sections~\ref{subsec:groupactions}--\ref{subsec:virasoro-heisenberg}
recall group actions, intertwining images, tensor products, and the
Virasoro and Heisenberg constructions.

\subsection{Notation, vertex operator identities, and invariant forms}
\label{subsec:conventions}\label{subsec:coefficientidentities}\label{subsec:invariantforms}\label{subsec:strong}

All vector spaces and tensor products are over $\C$, and
$\N=\Z_{\geq0}$. We follow \cite{FHL,LL} for vertex operators, modules,
and the identities below, and retain the quadruple notation used in
\cite{DLM}.

\begin{definition}\label{def:vertexoperatoralgebra}
A \emph{vertex operator algebra} is a quadruple $(V,Y,\one,\omega)$
with a $\Z$-graded complex vector space
\[
V=\bigoplus_{n\in\Z}V_n,
\]
where each $V_n$ is finite dimensional and $V_n=0$ for all sufficiently
negative $n$, a nonzero vacuum vector $\one\in V_0$, a conformal vector
$\omega\in V_2$, and a linear map
\[
Y:V\longrightarrow (\operatorname{End}_{\C}V)[[z,z^{-1}]],\qquad
a\longmapsto Y(a,z)=\sum_{n\in\Z}a_nz^{-n-1},
\]
satisfying the following axioms.
\begin{enumerate}[label=\textup{(\roman*)},leftmargin=*]
\item \emph{Truncation.} For every $a,b\in V$, one has $a_nb=0$ for all
sufficiently large integers $n$; equivalently, $Y(a,z)b\in V((z))$.

\item \emph{Vacuum and creation.} For every $a\in V$,
\[
Y(\one,z)=\operatorname{id}_V,\qquad
Y(a,z)\one\in V[[z]],\qquad
\left.Y(a,z)\one\right|_{z=0}=a.
\]

\item \emph{Jacobi identity.} Let $z_0,z_1,z_2$ be commuting formal
variables and put $\delta(z)=\sum_{n\in\Z}z^n$. For every $a,b\in V$,
\[
\begin{aligned}
z_0^{-1}\delta\!\left(\frac{z_1-z_2}{z_0}\right)
 Y(a,z_1)Y(b,z_2)&\\
-z_0^{-1}\delta\!\left(\frac{z_2-z_1}{-z_0}\right)
 Y(b,z_2)Y(a,z_1)
&=z_2^{-1}\delta\!\left(\frac{z_1-z_0}{z_2}\right)
 Y(Y(a,z_0)b,z_2).
\end{aligned}
\]
Each integral power of a difference is expanded in nonnegative powers
of its second variable, with the differences ordered as displayed.

\item \emph{Conformal structure.} Write
\[
Y(\omega,z)=\sum_{n\in\Z}L(n)z^{-n-2}.
\]
There is a scalar $c\in\C$, called the central charge, such that,
for all $m,n\in\Z$,
\[
[L(m),L(n)]=(m-n)L(m+n)
 +\frac{c}{12}(m^3-m)\delta_{m+n,0}\operatorname{id}_V.
\]
The grading and translation satisfy
\[
L(0)a=na\quad(a\in V_n),\qquad
Y(L(-1)a,z)=\frac{d}{dz}Y(a,z)\quad(a\in V).
\]
\end{enumerate}
\end{definition}

Here $a_n$ is the coefficient with index $n$ of $Y(a,z)$. We put
$D=L(-1)$. The translation and vacuum identities give
\[
\begin{aligned}
Y(Da,z)&=\frac{d}{dz}Y(a,z)=[D,Y(a,z)],\\
(Da)_r&=[D,a_r]=-r a_{r-1},\\
Y(a,z)\one&=e^{zD}a.
\end{aligned}
\]

A vertex operator algebra is called of \emph{CFT type} if
$V=\bigoplus_{d\in\N}V_d$ and $V_0=\C\one$.
The vertex operator algebras considered here are of CFT type. Put
$V_+=V_{>0}=\bigoplus_{d>0}V_d$ and
$DV=\{Da:a\in V\}$.
For homogeneous $a,b\in V$,
$\wt(a_n b)=\wt(a)+\wt(b)-n-1$ whenever the product is nonzero.
Here $\wt(a)=d$ for $a\in V_d$.
The contragredient module is denoted by
$V'=\bigoplus_{d\geq0}V_d^*$, where $V_d^*=\operatorname{Hom}_{\C}(V_d,\C)$,
with its contragredient vertex operators.
Self contragredience means $V\cong V'$ as $V$-modules.
For $n\in\Z$ and $k\geq0$, we use the generalized binomial
coefficient $\binom nk=n(n-1)\cdots(n-k+1)/k!$, with
$\binom n0=1$.
The symbol $\delta_{i,j}$ denotes the Kronecker delta.
Formal Laurent series in $z$ with a finite lower bound form $\C((z))$.
We write $\C(z)$ for the field of rational functions in $z$ and
$[z^j]f$ for the coefficient of $z^j$ in the specified expansion of $f$.
For a finite automorphism group $G$, the averaging operator is
$E_G=|G|^{-1}\sum_{g\in G}g$, abbreviated to $E$ when $G$ is fixed.
For $q\geq2$, put
\[
C_q(V)=\Span\{a_{-q}b:a,b\in V\}.
\]
The notation $\Span$ denotes complex linear span. For operators,
$[A,B]=AB-BA$ denotes the commutator.
Square brackets denote classes in the stated quotient.
A vertex operator algebra is $C_2$-cofinite when $\dim V/C_2(V)<\infty$.

For a module $M$ over a vertex operator algebra $U$ of CFT type, we use the module convention
\begin{equation}\label{eq:c1}
\Crel{U}{M}=\Span\{u_{-1}m:u\in U_{>0},\ m\in M\}.
\end{equation}
We do not add $L(-1)M$ to this subspace. A module is called $C_1$-cofinite if
its quotient by \eqref{eq:c1} is finite dimensional. This is the module
convention used in \cite{Mfus}; compare also
\cite{Huang}. In particular, $C_1^U(U)=U_+$ under
\eqref{eq:c1}.
An $\N$-gradable weak $U$-module has a decomposition
$M=\bigoplus_{r\in\N}M(r)$ satisfying
$u_jM(r)\subset M(r+\wt(u)-j-1)$ for homogeneous $u$,
where $M(t)=0$ for $t<0$. For submodules of a conformal vertex operator algebra extension,
we use the $L(0)$ grading.

A homogeneous subspace $S\subset U_+$ strongly generates $U$ if
$U$ is spanned by $\one$ and the vectors
$s^1_{-m_1}\cdots s^k_{-m_k}\one$, with $s^i\in S$ and $m_i\geq1$.

We recall the commutator, iterate, and skew symmetry identities from \cite{FHL,LL}.

\begin{proposition}
\label{prop:external-identities}
For vectors $u,a$ of a vertex algebra and integers $p,q$, one has
\begin{align}
[u_p,a_q]&=\sum_{j\geq0}\binom pj(u_ja)_{p+q-j},
\label{eq:commutator}\\
(u_{-1}a)_q&=\sum_{i\geq0}
\bigl(u_{-1-i}a_{q+i}+a_{q-1-i}u_i\bigr),
\label{eq:iterate}\\
Y(a,z)b&=e^{zD}Y(b,-z)a.
\label{eq:skew}
\end{align}
When applied to a fixed vector, all sums that require truncation are
finite. \end{proposition}

The weak commutativity and weak associativity identities are recalled from \cite{LL}.

\begin{proposition}\label{prop:external-weakidentities}
For $a,b\in V$, there is an integer $k\geq0$,
depending only on $a,b$, such that, for every $c\in V$,
\begin{equation*}\label{eq:elementcommutativity}
(z_1-z_2)^kY(a,z_1)Y(b,z_2)c
=(z_1-z_2)^kY(b,z_2)Y(a,z_1)c.
\end{equation*}
For $a,c\in V$, choose $\ell\geq0$ with
$z^\ell Y(a,z)c\in V[[z]]$. Then, for every $b\in V$,
\begin{equation}\label{eq:elementassociativity}
(z_0+z_2)^\ell Y(a,z_0+z_2)Y(b,z_2)c
=(z_0+z_2)^\ell Y(Y(a,z_0)b,z_2)c.
\end{equation}
On the left of \eqref{eq:elementassociativity}, integer powers of
$z_0+z_2$ are expanded in nonnegative powers of $z_2$.
These are the weak commutativity and weak associativity identities.
\end{proposition}

We use the following description of invariant bilinear forms \cite{FHL,LiForm}.

\begin{proposition}
\label{prop:external-invariant-form}
For a vertex operator algebra $V$, invariant bilinear forms correspond to
$V$-module homomorphisms $\varphi:V\to V'$ by
$(u,v)=\langle\varphi(u),v\rangle$. Every such form is symmetric;
invariance means
\[
(Y(a,z)u,v)
=(u,Y(e^{zL(1)}(-z^{-2})^{L(0)}a,z^{-1})v).
\]
\end{proposition}

\subsection{Products and the Poisson quotient}
\label{subsec:products}\label{subsec:c2generation}

Let $V$ be a vertex operator algebra of CFT type.
Define the $-1$ product by $a*b=a_{-1}b$.
For $a\in V$, write $[a]=a+C_2(V)$.

The product satisfies the following identities \cite{DA,LL}.

\begin{proposition}
\label{prop:external-product-identities}
The vacuum $\one$ is a two sided unit for $*$. For $a,b,c\in V$,
\begin{equation}\label{eq:product-associator}
(a*b)*c-a*(b*c)
=\sum_{j\geq0}\bigl(a_{-2-j}b_jc+b_{-2-j}a_jc\bigr).
\end{equation}
The right hand side is symmetric in $a,b$. Moreover,
\begin{equation}\label{eq:trace}
a*b-b*a\in DV.
\end{equation}
For $v\in V$, the derivation identities are
\begin{equation}
\label{eq:productderivations}
D(a*b)=(Da)*b+a*(Db),\qquad v_0(a*b)=(v_0a)*b+a*(v_0b).
\end{equation}
\end{proposition}

Zhu's quotient has the following algebra structure \cite{Zhu}.

\begin{proposition}
\label{prop:external-c2-algebra}
For a vertex operator algebra $V$, the formulas
\[
[a][b]=[a*b]=[a_{-1}b],\qquad \{[a],[b]\}=[a_0b]
\]
make $R_V=V/C_2(V)$ a commutative associative Poisson algebra
with unit $[\one]$.
\end{proposition}

\subsection{Commutative algebra}
\label{subsec:descenttools}\label{subsec:descentauxiliary}

For a commutative unital $\C$-algebra $A$, an \emph{ideal}
is a linear subspace $J\subseteq A$ such that $aj\in J$ for every
$a\in A$ and $j\in J$. For a subset $S\subseteq A$, the ideal
generated by $S$ is
\[
(S)_A=\left\{\sum_{i=1}^m a_i s_i:
m\in\N,\ a_i\in A,\ s_i\in S\right\},
\]
where the empty sum is $0$. This is the smallest ideal containing
$S$. A subset $S\subseteq J$ is a set of \emph{ideal generators}
of $J$ if $J=(S)_A$. For a finite set, we write
\[
(s_1,\ldots,s_m)_A=A s_1+\cdots+A s_m.
\]

A \emph{unital character} of a commutative unital $\C$-algebra
$R$ is a unital $\C$-algebra homomorphism $\chi:R\to\C$.
Explicitly, for $a,b\in R$ and $\lambda,\mu\in\C$, it satisfies
\[
\chi(\lambda a+\mu b)=\lambda\chi(a)+\mu\chi(b),\qquad
\chi(ab)=\chi(a)\chi(b),\qquad \chi(1_R)=1.
\]
For a graded algebra $R$, no compatibility with the grading is
required.

For any nonnegatively graded commutative algebra
$R=\C\one\oplus R_+$, the ideal $R_+^2$ is the span of products
of two elements of $R_+$.
Thus $R_+/R_+^2$ is the quotient of $R_+$ by these products.

For an ideal $J$ in a commutative unital algebra $A$, write
\[
\sqrt{J}=\{a\in A:a^m\in J\text{ for some integer }m\ge1\}.
\]

We use the following form of the Nullstellensatz \cite{AM}.

\begin{proposition}
\label{prop:external-nullstellensatz}
Let $k$ be algebraically closed and $J$ an ideal of
$k[x_1,\ldots,x_t]$. For
\[
Z(J)=\{a\in k^t:f(a)=0\text{ for all }f\in J\},
\]
one has
\[
\{f\in k[x_1,\ldots,x_t]: f(a)=0\text{ for all }a\in Z(J)\}
=\sqrt J.
\]
\end{proposition}

\subsection{Cofiniteness}
\label{subsec:finiteness}

The base algebra in each statement is specified separately;
we use the relative module convention of \eqref{eq:c1}.

The finiteness property below follows from \cite{GN,ABD}.

\begin{lemma}\label{lem:cq}
If a vertex operator algebra $V$ of CFT type is $C_2$-cofinite, then $V/C_q(V)$ is
finite dimensional for every integer $q\geq2$.
\end{lemma}

Let $U$ be a vertex operator algebra of CFT type.
For an $\N$-gradable weak $U$-module $M$ and $m\ge1$, put
\[
C_m^U(M)=\Span\{u_{-m}v:u\in U\text{ homogeneous},\
\wt(u)>1-m,\ v\in M\}.
\]
At $m=1$ this is exactly \eqref{eq:c1}.

For irreducible weak modules, we use the cofiniteness result in \cite{ABD}.

\begin{proposition}
\label{prop:external-module-c2}
If $V$ is $C_2$-cofinite and $M$ is an irreducible weak $V$ module,
then $M$ is $C_2$-cofinite.
Here
\[
C_2^V(M)=\Span\{a_{-2}m:a\in V,\ m\in M\}.
\]
\end{proposition}

\subsection{Expansions of rational functions}
\label{subsec:rationality}

Recall the restricted dual from Section~\ref{subsec:conventions}
and its evaluation pairing:
\[
V'=\bigoplus_{m\geq0}V_m^*,\qquad
\langle\theta,v\rangle=\theta(v)
\quad(\theta\in V',\ v\in V);
\]
Each element of $V_m^*$ is extended by zero
on all other weight spaces. We extend the pairing coefficientwise
to vector valued formal series by applying the functional to
each coefficient.

Write $\iota_{s,t}$ for expansion of mixed
factors in nonnegative powers of $t/s$.

The following rationality result for matrix coefficients is recalled from \cite{FHL}.

\begin{proposition}\label{prop:external-rationality}
For $a,b,c\in V$ and
$\theta\in V'$, there is a rational function
\[
f(z_1,z_2)\in\C[z_1^{\pm1},z_2^{\pm1},(z_1-z_2)^{-1}]
\]
such that
\[
\iota_{z_1,z_2}f(z_1,z_2)
=\langle\theta,Y(a,z_1)Y(b,z_2)c\rangle.
\]
Commutativity and associativity identify its other expansions as
\begin{align}
\iota_{z_2,z_1}f(z_1,z_2)
&=\langle\theta,Y(b,z_2)Y(a,z_1)c\rangle,
\notag\\
\iota_{z_2,z_0}f(z_2+z_0,z_2)
&=\langle\theta,Y(Y(a,z_0)b,z_2)c\rangle.
\label{eq:rationalassociativity}
\end{align}
\end{proposition}

We record the product and associativity results of
\cite{FHL,HuangCoh} in the case of the adjoint module.

\begin{proposition}
\label{prop:external-multiproduct}
Let $n\ge1$, $v_1,\ldots,v_n,w\in V$ and $\theta\in V'$.
There is a rational function
\[
F\in\C[z_1^{\pm1},\ldots,z_n^{\pm1},
             (z_i-z_j)^{-1}:1\le i<j\le n]
\]
such that, for every $\sigma\in S_n$, the series
\[
\bigl\langle\theta,
Y(v_{\sigma(1)},z_{\sigma(1)})\cdots
Y(v_{\sigma(n)},z_{\sigma(n)})w\bigr\rangle
\]
converges absolutely to $F(z_1,\ldots,z_n)$ whenever
$|z_{\sigma(1)}|>\cdots>|z_{\sigma(n)}|>0$.
\end{proposition}

The corresponding identity for an adjacent iterate is given in \cite{HuangCoh}.

\begin{proposition}
\label{prop:external-adjacentassociativity}
With the notation of Proposition~\ref{prop:external-multiproduct},
suppose that $n\ge2$ and $1\le i<n$. The series
\[
\bigl\langle\theta,\,Y(v_1,z_1)\cdots Y(v_{i-1},z_{i-1}) \cdot Y\bigl(Y(v_i,z_i-z_{i+1})v_{i+1},z_{i+1}\bigr) \cdot Y(v_{i+2},z_{i+2})\cdots Y(v_n,z_n)w\bigr\rangle
\]
converges absolutely to $F(z_1,\ldots,z_n)$ on the region
\[
\begin{gathered}
|z_1|>\cdots>|z_{i-1}|>|z_{i+1}|>\cdots>|z_n|>0,\\
|z_{i+1}|>|z_i-z_{i+1}|>0,\\
|z_k-z_{i+1}|>|z_i-z_{i+1}|\qquad(k\ne i,i+1).
\end{gathered}
\]
Empty products are interpreted as the identity.
\end{proposition}

For three independent formal variables, put
\[
\mathcal R_{1;23}=\C[z_1^{\pm1},z_2^{\pm1},z_3^{\pm1}, (z_1+z_2)^{-1},(z_1-z_3)^{-1},(z_1+z_2-z_3)^{-1}].
\]
Define the $\C$-algebra homomorphism
\[
\iota_{z_1;z_2,z_3}:\mathcal R_{1;23}
\longrightarrow
\C[z_1^{\pm1}][[z_2,z_3]][z_2^{-1},z_3^{-1}]
\]
by fixing Laurent monomials and, for $q\ge1$ and
$(\alpha,\beta)\in\{(1,0),(0,-1),(1,-1)\}$, setting
\[
\iota_{z_1;z_2,z_3}(z_1+\alpha z_2+\beta z_3)^{-q}
=\sum_{r\ge0}(-1)^r\binom{q+r-1}{r}
z_1^{-q-r}(\alpha z_2+\beta z_3)^r.
\]
Thus the two variables after the semicolon are expanded jointly
in nonnegative powers relative to $z_1$, while their Laurent
monomial factors are retained.

\begin{lemma}\label{lem:matrixassociativity}
For $a,b,c,d\in V$ and
$\theta\in V'$, let $f_\theta(z_1,z_2,z_3)$ be the rational function
whose expansion for $|z_1|>|z_2|>|z_3|>0$ is
\[
\left\langle\theta,Y(a,z_1)Y(b,z_2)Y(c,z_3)d\right\rangle.
\]
For independent formal variables $x,z',z''$, this gives
\begin{equation}\label{eq:matrixassociativity}
\iota_{x;z',z''}f_\theta(x+z',x,z'')
=\left\langle\theta,
Y\bigl(Y(a,z')b,x\bigr)Y(c,z'')d\right\rangle.
\end{equation}
\end{lemma}

To obtain \eqref{eq:matrixassociativity} from
Proposition~\ref{prop:external-adjacentassociativity}, take $n=3$, $i=1$,
$(v_1,v_2,v_3,w;\theta)=(a,b,c,d;\theta)$, and
$(z_1,z_2,z_3)=(x+z',x,z'')$.
The region $|x|>|z'|+|z''|$, $|z'|,|z''|>0$ lies in the
convergence region of that proposition. Uniqueness of Laurent
expansions identifies its series there with
$\iota_{x;z',z''}f_\theta(x+z',x,z'')$.
See also \cite{Qi}.

\subsection{Algebraic functions}
\label{subsec:algebraicfunctions}

We use the local expansion theorem for algebraic functions in \cite{Forster}.

\begin{proposition}
\label{prop:external-puiseux}
Let $f$ be algebraic over $\C(z)$ and let $a\in\mathbb P^1$.
Each local branch admits a convergent expansion
\[
f(z)=\sum_{j\ge j_0}b_jt^{j/e},\qquad
j_0\in\Z,\quad e\ge1,
\]
on a sufficiently small punctured neighborhood, with
$t=z-a$ if $a\ne\infty$ and $t=z^{-1}$ if $a=\infty$.
\end{proposition}

\begin{lemma}\label{lem:algebraicODE}
Let $f\in\C((z))^r$ have components algebraic over $\C(z)$.
If $f'=A(z)f$ with $A(z)\in\operatorname{Mat}_r(\C[z,z^{-1}])$,
then every component of $f$ is a Laurent polynomial.
\end{lemma}
\begin{proof}
By Proposition~\ref{prop:external-puiseux}, each algebraic
Laurent series $f_i$ converges in a punctured disk about $0$.
Its integral exponents and finite principal part make it
single valued and meromorphic there.

The matrix $A$ is holomorphic on $\C^*$. To describe the local
solutions at $a\ne0$, write $t=z-a$ and
$A(a+t)=\sum_{j\ge0}A_jt^j$. A fundamental matrix with $H(0)=I$
has coefficients determined by
\[
H(t)=\sum_{m\ge0}H_mt^m,\qquad H_0=I,\qquad
(m+1)H_{m+1}=\sum_{j=0}^m A_jH_{m-j}.
\]
The local existence and uniqueness theorem for analytic linear
systems gives convergence and continuation along every path
in $\C^*$ \cite{NY}. The recurrence also gives uniqueness for
each initial value. Since the fundamental group of $\C^*$ is
generated by a small loop about $0$, and the original Laurent
series is unchanged by that loop, continuation yields a
single valued holomorphic vector $f$ on $\C^*$.

At infinity, set $t=z^{-1}$. Proposition~\ref{prop:external-puiseux}
gives, for each component,
\[
f_i(z)=\sum_{j\ge j_i}b_{i,j}t^{j/e_i}.
\]
Single valuedness under $t\mapsto e^{2\pi\mathrm i}t$ implies
\[
\bigl(e^{2\pi\mathrm i j/e_i}-1\bigr)b_{i,j}=0.
\]
Hence only integral powers of $t$ occur, and $f_i$ has at most
a pole at infinity. Choose $N_i\ge0$ so that
$g_i(z)=z^{N_i}f_i(z)$ is holomorphic at $0$. Then $g_i$ is entire
and satisfies $|g_i(z)|\le C_i|z|^{M_i}$ for large $|z|$, for some
integer $M_i\ge0$. Cauchy's estimate on the circle $|z|=R$ gives
\[
\frac{|g_i^{(k)}(0)|}{k!}\le C_iR^{M_i-k}.
\]
For $k>M_i$, letting $R\to\infty$ gives $g_i^{(k)}(0)=0$.
Thus $g_i$ is a polynomial and $f_i=z^{-N_i}g_i$ is a Laurent
polynomial, as required.
\end{proof}

\subsection{Group actions}
\label{subsec:groupactions}\label{subsec:exampleinputs}\label{subsec:solvablebackground}

Miyamoto's theorem for solvable groups gives the following result \cite{Msol}.

\begin{theorem}
\label{thm:solvable}
For a simple $C_2$-cofinite vertex operator algebra $V$ of CFT type and a finite
solvable subgroup $G\le\operatorname{Aut}(V)$,
\[
\dim V^G/C_2(V^G)<\infty.
\]
\end{theorem}

The following is the Schur--Weyl duality theorem for compact
automorphism groups of simple vertex operator algebras \cite{DLM}.

\begin{theorem}
\label{thm:external-galois}
Let $X$ be a simple vertex operator algebra and let a compact Lie group $G$ act
faithfully and continuously by automorphisms. Index its
finite dimensional irreducible continuous representations
$M_\chi$ by $\chi$. Then, as a $G\times X^G$-module,
\begin{equation*}
X=\bigoplus_\chi M_\chi\otimes X_\chi,
\end{equation*}
where every $X_\chi$ is nonzero and irreducible, and
$X_\chi\cong X_\psi$ as $X^G$-modules if and only if $\chi=\psi$.
\end{theorem}

For finite $G$, there are only finitely many irreducible representations.

A vertex operator algebra is called \emph{strongly rational} if it is simple,
self contragredient, rational, and $C_2$-cofinite of CFT type.

The following rationality result is due to McRae \cite{McRaeOrb}.

\begin{theorem}
\label{thm:external-orbifold-rationality}
For a strongly rational vertex operator algebra $V$ and a finite subgroup
$G\le\operatorname{Aut}(V)$, $C_2$-cofiniteness of $V^G$
implies that $V^G$ is strongly rational.
\end{theorem}

\subsection{Intertwining images}
\label{subsec:intertwiningimages}

Let $U\subseteq V$ be a conformal inclusion of vertex operator
algebras of CFT type, with the same conformal vector.
For $U$-submodules $B,C$ of $V$, the restriction
\[
\mathcal Y_{B,C}(b,z)c=Y_V(b,z)c
=\sum_{n\in\Z}(b_nc)z^{-n-1}
\]
is an ordinary intertwining operator for $U$ \cite{FHL,LL}.
Miyamoto \cite{Mfus} constructs the target of a restricted
intertwining operator by taking the linear span of all its
coefficients. Applying that construction to $\mathcal Y_{B,C}$,
we write
\begin{equation}\label{eq:diamond}
B\diamond C=\Span\{b_nc:b\in B,\ c\in C,\ n\in\mathbb Z\}.
\end{equation}
Thus \eqref{eq:diamond} is the definition of this coefficient
image,\newline
obtained by extracting $[z^{-n-1}]\mathcal Y_{B,C}(b,z)c$.
The symbol $\diamond$ is our notation for that image.

The bound on coefficient images is supplied by the following theorem \cite{Mfus}.

\begin{theorem}
\label{thm:image}
Fix $m\ge1$. Let $B,C$ be $C_m$-cofinite $\N$-gradable
weak $U$-modules. There is a bound $f_m(B,C)$ such that, for
every $\N$-gradable weak $U$-module $P$ and every surjective
logarithmic intertwining operator $\mathcal Y$ of type
$\binom{P}{B\ C}$,
\[
\dim P/C_m^U(P)<f_m(B,C).
\]
\end{theorem}
Surjectivity means that the coefficients span the target.

\subsection{Tensor products}
\label{subsec:tensorproducts}\label{subsec:tensorproductproof}\label{subsec:permutationbackground}

For vertex operator algebras $V^{(1)},\ldots,V^{(n)}$, the tensor product construction
has vertex operators and conformal vector
\[
\begin{aligned}
Y_{\otimes_i V^{(i)}}(a_1\otimes\cdots\otimes a_n,z)
&=Y_{V^{(1)}}(a_1,z)\otimes\cdots\otimes Y_{V^{(n)}}(a_n,z),\\
\omega_{\otimes_i V^{(i)}}
&=\sum_{i=1}^n\one^{\otimes(i-1)}\otimes\omega_i
                  \otimes\one^{\otimes(n-i)}.
\end{aligned}
\]
Here $\omega_i$ is the conformal vector of $V^{(i)}$.
For factors of CFT type, its grading is
\[
\left(\bigotimes_{i=1}^n V^{(i)}\right)_d
=\bigoplus_{d_1+\cdots+d_n=d}\bigotimes_{i=1}^n V^{(i)}_{d_i},
\qquad
\left(\bigotimes_{i=1}^n V^{(i)}\right)_0=\C\one^{\otimes n},
\]
so the tensor product is again of CFT type \cite{FHL}.

Simplicity of tensor products is characterized as follows \cite{FHL}.

\begin{proposition}\label{prop:tensorsimple}
The vertex operator algebras $V^{(1)},\ldots,V^{(n)}$ are all simple if and only if
their tensor product $V^{(1)}\otimes\cdots\otimes V^{(n)}$ is simple.
\end{proposition}

Cofiniteness of tensor products is treated in \cite{ABD}.

\begin{lemma}\label{lem:tensorc2}
If the vertex operator algebras $V^{(1)},\ldots,V^{(n)}$ are $C_2$-cofinite, then
$V^{(1)}\otimes\cdots\otimes V^{(n)}$ is $C_2$-cofinite.
\end{lemma}

We now fix a vertex operator algebra $V$ of CFT type and consider
its tensor powers $T^n(V)=V^{\otimes n}$.
For $a\in V$, let $a^{(i)}$ be supported in position $i$, with
$\one$ in the other positions, and put $\eta^{(n)}(a)=\sum_{i=1}^n a^{(i)}$.
The conformal vector of $T^n(V)$ is $\eta^{(n)}(\omega)$ and its translation
operator is denoted by
\[
D_{(0)}=\eta^{(n)}(\omega)_0=\sum_{i=1}^nD^{(i)},
\]
where $D^{(i)}$ acts as $D$ in position $i$ and as the identity
on the other tensor factors.

We have $D_{(0)}\eta^{(n)}(a)=\eta^{(n)}(Da)$ and
$\eta^{(n)}(\one)=n\one^{\otimes n}$.

For $n\geq1$, retain $T^n(V)$ and let $S_n$ be the
symmetric group on $n$ positions, and set
\[
U^{(n)}=T^n(V)^{S_n},\qquad R^{(n)}=U^{(n)}/C_2(U^{(n)}),\qquad \m_n=(R^{(n)})_+=\bigoplus_{d>0}(R^{(n)})_d.
\]
Here $(R^{(n)})_d$ denotes the homogeneous component of degree $d$,
and $\m_n$ is the ideal formed by the components of positive degree.
The action is
$\sigma(v_1\otimes\cdots\otimes v_n)
=v_{\sigma^{-1}(1)}\otimes\cdots\otimes v_{\sigma^{-1}(n)}$.
For $n\geq2$, let $S_{n-1}$ fix the last position and set
\[
W^{(n)}=U^{(n-1)}\otimes V=T^n(V)^{S_{n-1}},
\]
using the canonical inclusion into the first $n-1$ tensor positions
and the last position. Thus $U^{(n)}\subset W^{(n)}\subset T^n(V)$.

We write $\Sym^k(E)$ for the $k$th symmetric power of a vector
space $E$, and $\Sym(E)=\bigoplus_{k\geq0}\Sym^k(E)$.
Following Abe's notation and normalization
\cite{Abe}, for $a_1,\ldots,a_k\in V$
and $1\leq k\leq n$ set
\begin{equation*}
\phi_k^{(n)}(a_1,\ldots,a_k) =\frac{1}{(n-k)!}\sum_{\sigma\in S_n} \sigma\bigl(a_1\otimes\cdots\otimes a_k \otimes\one^{\otimes(n-k)}\bigr),\qquad \eta^{(n)}(a)=\phi_1^{(n)}(a).
\end{equation*}
This is exactly the sum over injections of the $k$ labeled arguments
into the $n$ positions: each injection occurs $(n-k)!$ times in the
permutation sum. We also use the convention for zero arguments
$\phi_0^{(n)}=\one^{\otimes n}$ and Abe's convention
$\phi_{n+1}^{(n)}=0$. Restricting the arguments to $V_+$ gives a
vector space identification
\[
U^{(n)}\cong\bigoplus_{k=0}^n\Sym^k(V_+).
\]
Abe's multiplication identity takes the following form \cite{Abe}.

\begin{proposition}
\label{prop:external-abe}
For $a,a_1,\ldots,a_k\in V$ and $1\le k\le n$,
\begin{equation}
\label{eq:Phimultiply}
\eta^{(n)}(a)_{-1}\phi_k^{(n)}(a_1,\ldots,a_k) ={}\phi_{k+1}^{(n)}(a,a_1,\ldots,a_k) +\sum_{j=1}^k\phi_k^{(n)}(a_1,\ldots,a*a_j,\ldots,a_k).
\end{equation}
\end{proposition}

The corresponding strong generation statement is also due to Abe \cite{Abe}.

\begin{proposition}
\label{prop:external-abe-generation}
$U^{(n)}$ is strongly generated by $\eta^{(n)}(V)$.
\end{proposition}

\subsection{Virasoro and Heisenberg vertex operator algebras}
\label{subsec:virasoro-heisenberg}\label{subsec:heisenbergmodules}

The constructions in this subsection follow \cite{LL}.

For $c\in\C$, let $\mathcal V_c=V_{\mathrm{Vir}}(c,0)$ be
the vertex operator algebra associated with the Virasoro algebra
in \cite{LL}. It has conformal vector $\omega_c=L(-2)\one$ and
translation operator $D=L(-1)$.
We also write $\omega$ for $\omega_c$ when
only $\mathcal V_c$ is under consideration.
As a module over the Virasoro algebra, it is induced from the
one dimensional module $\C\one$ of the subalgebra spanned by the
$L(m)$ with $m\geq-1$ and the central element. On $\C\one$ these $L(m)$
act by zero and the central element acts by $c$. It has basis
\[
v_\lambda=L(-\lambda_1)\cdots L(-\lambda_r)\one,
\qquad \lambda_1\geq\cdots\geq\lambda_r\geq2.
\]
The empty partition gives $\one$. Write
$|\lambda|=\sum_i\lambda_i$ and $\ell(\lambda)=r$.
The defining relations are
\begin{equation*}
[L(m),L(n)]=(m-n)L(m+n)
 +\frac{m^3-m}{12}c\,\delta_{m+n,0},
\qquad L(m)\one=0\quad(m\geq-1).
\end{equation*}

Every vertex operator algebra of central charge $c_V$ receives a graded vertex algebra
homomorphism $\mathcal V_{c_V}\to V$, with
$L(-2)\one\mapsto\omega$.

Let $\mathcal H=M(1)=V_{\mathfrak h}(1,0)$ be the rank one
vertex operator algebra associated with the Heisenberg algebra
of $\mathfrak h=\C h$, with $\langle h,h\rangle=1$, as in \cite{LL}.
We identify $h$ with $h_{-1}\one$ and write
$Y(h,z)=\sum_{m\in\Z}h_mz^{-m-1}$. Its normalization is
\begin{equation*}
[h_m,h_n]=m\delta_{m+n,0}\operatorname{id},
\qquad h_m\one=0\quad(m\geq0).
\end{equation*}
The standard conformal vector $\omega_{\mathcal H}=\frac12h_{-1}^2\one$
makes $\mathcal H$ a vertex operator algebra of CFT type and central charge one, with
$\wt(h)=1$ and translation operator $D=L_{\mathcal H}(-1)$.

The following is the rank one consequence of the affine algebra
construction and the homomorphism criterion for vacuum like vectors
in \cite{LL,LiLocal}.

\begin{proposition}\label{prop:heisenberg-homomorphism}
Let $V$ be a vertex operator algebra of CFT type and let $a\in V_1$ satisfy
\[
a_ja=\delta_{j,1}\one\qquad(j\ge0).
\]
There is a unique graded vertex algebra homomorphism
$f_a:\mathcal H\to V$ with $f_a(h)=a$. It is given by
\[
f_a(h_{-n_1}\cdots h_{-n_s}\one)
=a_{-n_1}\cdots a_{-n_s}\one
\qquad(s\ge0,\ n_1,\ldots,n_s\ge1),
\]
where $s=0$ gives $f_a(\one)=\one$. In particular,
\[
f_aD_{\mathcal H}=D_Vf_a,\qquad
f_a(\omega_{\mathcal H})=\tfrac12a_{-1}a.
\]
\end{proposition}
\begin{proof}
The commutator formula gives
\[
[a_m,a_n]=m\delta_{m+n,0}\operatorname{id}_V,\qquad
[D_V,a_m]=-ma_{m-1}.
\]
The affine algebra construction of \cite{LL,LiLocal}, applied to
$\mathfrak h=\C h$ at level one, makes $V$ an $\mathcal H$-module
with the field $Y_V(a,z)$ acting for $h$.
The creation property gives $a_m\one=0$ for $m\ge0$ and $D_V\one=0$.
The homomorphism criterion for vacuum like vectors therefore gives
the displayed map. The iterated normally ordered products in this
module are those in $V$, so the map preserves all vertex products.
It preserves weights because $\wt(a)=1$, and uniqueness follows
because $h$ generates $\mathcal H$. Preservation of the vacuum and
vertex products gives the two final identities.
\end{proof}

Write $\mathcal H(n)=\mathcal H^{\otimes n}$ and
$A_n^{\mathcal H}=\mathcal H(n)^{O(n,\C)}$. The orthogonal
group acts on the weight one generating space in its standard
representation.

We use the decomposition described by Linshaw \cite{LinHeis}.

\begin{theorem}
\label{thm:heisenbergmodules}
There is an algebraic $O(n,\C)\times A_n^{\mathcal H}$-decomposition
\[
\mathcal H(n)=\bigoplus_\lambda L_\lambda\otimes M_\lambda,
\]
where $L_\lambda$ is finite dimensional and irreducible, and
$M_\lambda$ is an irreducible highest weight
$A_n^{\mathcal H}$-module.
\end{theorem}

The decomposition and irreducibility in
Theorem~\ref{thm:heisenbergmodules} are the Heisenberg specialization
of Theorem~\ref{thm:external-galois}. Apply the latter with
$X=\mathcal H(n)$ and the compact group $K=O(n,\mathbb R)$.
Each homogeneous subspace is a finite dimensional algebraic
$O(n,\C)$-module, and $K$ is Zariski dense in $O(n,\C)$. Hence
\[
\mathcal H(n)^K=\mathcal H(n)^{O(n,\C)}=A_n^{\mathcal H},
\qquad
\operatorname{Hom}_K(L_\lambda,\mathcal H(n))
=\operatorname{Hom}_{O(n,\C)}(L_\lambda,\mathcal H(n)).
\]
Thus the compact group decomposition gives exactly the displayed
algebraic $O(n,\C)\times A_n^{\mathcal H}$-decomposition.
The highest weight assertion additionally uses the commutativity
of the Zhu algebra of $A_n^{\mathcal H}$ established in
\cite{LinHeis}: each $M_\lambda$ is generated by its
one dimensional lowest weight subspace.

The following finiteness result uses \cite{LinHeis,LinAff}.

\begin{proposition}\label{prop:heisenberg-wick}
For every irreducible highest weight $A_n^{\mathcal H}$-submodule
$M$ of $\mathcal H(n)$, the subspace
\[
M_{\mathrm{Wick}}=
\Span\{a_{-1}m:a\in(A_n^{\mathcal H})_{>0},\ m\in M\}
\]
has finite codimension in $M$.
\end{proposition}
The statement is conditional in \cite{LinHeis};
\cite{LinAff} proves the required strong finite generation and
explicitly removes that condition.

\Needspace{16\baselineskip}
\section{Theorem~\ref*{thm:descent}: Relative \texorpdfstring{$C_1$}{C1} descent}\label{sec:relativedescent}

We prove Theorem~\ref{thm:descent}. In the fixed point applications
below, we retain its hypotheses and write $U=V^G$; the statements about graded algebras and rational functions are formulated with their own assumptions.

Section~\ref{subsec:gradedalgebra} gives the commutative algebra
criteria, and Section~\ref{sec:finitegeneration} proves finite
generation. Section~\ref{sec:productclosure} passes from a finite
differential system to multiplicative closure and algebraicity, using
uniform denominator bounds before summing homogeneous components.
Section~\ref{sec:descent} obtains the Laurent polynomial bound and
proves character vanishing.

\subsection{Graded commutative algebras}
\label{subsec:gradedalgebra}

\begin{theorem}\label{thm:gradedgeneration}
Let $R=\C\one\oplus R_+$ be a nonnegatively graded commutative
algebra with finite dimensional homogeneous pieces.
If $R_+/R_+^2$ is finite dimensional, then $R$ is finitely generated.
\end{theorem}
\begin{proof}
Choose homogeneous lifts $x_1,\ldots,x_s\in R_+$ of a basis of
$R_+/R_+^2$, and put $d_i=\deg x_i>0$ and
$R'=\C[x_1,\ldots,x_s]\subseteq R$. We show $R_d\subseteq R'$
by induction on $d$. The assertion holds for $d=0$. For $d>0$
and $r\in R_d$, the choice of the lifts gives
\[
r-\sum_{d_i=d}\lambda_i x_i
\in (R_+^2)_d=\sum_{e=1}^{d-1}R_eR_{d-e}.
\]
Every factor on the right has degree less than $d$, so induction
places the right side in $R'$. Thus $r\in R'$, and
$R=\C[x_1,\ldots,x_s]$.
If $s=0$, this gives $R=\C$.
\end{proof}

\begin{theorem}\label{lem:gradedalgebra}
Let $R=\C\one\oplus R_+$ be a finitely generated, nonnegatively
graded commutative algebra with finite dimensional homogeneous pieces.
If every unital character $R\to\C$ vanishes on $R_+$,
then $R$ is finite dimensional.
\end{theorem}
\begin{proof}
Choose homogeneous algebra generators
$x_1,\ldots,x_s\in R_+$. Such a set is obtained by
taking the positive degree homogeneous components of any
finite generating set. If $s=0$, then $R=\C$ and there
is nothing to prove.
Let $P=\C[X_1,\ldots,X_s]$ and let
$\pi:P\twoheadrightarrow R$ send $X_i$ to $x_i$. Write
$I=\ker\pi$. Each point $a\in Z(I)$ defines a unital character
$\chi_a(\pi(f))=f(a)$. Under the additional hypothesis, for $1\le i\le s$, we have
\[
a_i=\chi_a(x_i)=0.
\]
The augmentation $R\to R/R_+=\C$ supplies the point $0$, so
$Z(I)=\{0\}$. Proposition~\ref{prop:external-nullstellensatz} gives
\[
\sqrt I=(X_1,\ldots,X_s).
\]
For each $i$ choose $N_i\ge1$ with $X_i^{N_i}\in I$.
Then $x_i^{N_i}=0$, and every monomial reduces to one in the set
\[
\{x_1^{e_1}\cdots x_s^{e_s}:0\le e_i<N_i\}.
\]
Consequently $\dim R\le\prod_{i=1}^sN_i<\infty$.
\end{proof}

\begin{remark}\label{rem:ideal-generators}
The variables $X_i$ in the proof of Theorem~\ref{lem:gradedalgebra}
generate $\sqrt I$ as an ideal of $P$, and their images $x_i$
generate the ideal $R_+$ of $R$:
\[
\sqrt I=(X_1,\ldots,X_s)_P=\sum_{i=1}^s P X_i,\qquad R_+=(x_1,\ldots,x_s)_R=\sum_{i=1}^s R x_i.
\]
These statements concern ideal generation in commutative algebras,
not vertex operator algebra generation. The unit $\one$ belongs to neither ideal and
is not included among the generators. The notation
$R=\C[x_1,\ldots,x_s]$ already includes the scalars.
\end{remark}

\subsection{Finite generation of the fixed point \texorpdfstring{$C_2$}{C2} algebra}
\label{sec:finitegeneration}

For the descent argument, assume the hypotheses of Theorem~\ref{thm:descent}. We use
the vertex operator coefficients, translation operator $D$,
positive weight spaces, averaging operator, and quotient algebra
conventions of Section~\ref{subsec:conventions}, and write
\[
U=V^G,\qquad E=|G|^{-1}\sum_{g\in G}g,\qquad R=U/C_2(U).
\]
By Proposition~\ref{prop:external-c2-algebra}, $R$ is a graded
commutative associative algebra. To apply
Theorem~\ref{lem:gradedalgebra}, we first prove that
$R$ is finitely generated.

Since $C_2(U)\subseteq U\cap C_2(V)$, the inclusion $U\hookrightarrow V$
induces the graded algebra homomorphism
\[
\iota:R=U/C_2(U)\longrightarrow V/C_2(V),\qquad u+C_2(U)\longmapsto u+C_2(V).
\]
Set
\[
K=\ker\iota=(U\cap C_2(V))/C_2(U),
\qquad Q=K/(R_+K).
\]
Then $K$ is a graded ideal of $R$, and
\[
R/K\cong\operatorname{im}\iota\subseteq V/C_2(V),
\qquad
\dim_{\C}(R/K)\leq\dim_{\C}\bigl(V/C_2(V)\bigr)<\infty.
\]
Since $R/K$ is already finite dimensional, the remaining task is
to control generators coming from $K$. We consider $Q$ because it
discards products of positive degree elements of $R$ with elements
of $K$. In particular, the classes of $u_{-1}E(a_{-m}b)$ vanish
in $Q$ for $u\in U_+$ and $m\geq2$. This allows the averaged relations between vertex operator
coefficients to reduce to finitely many representatives of
$V/C_1^U(V)$.

The $G$-stability of $C_2(V)$ and $E|_U=\operatorname{id}_U$ give
\[
E(C_2(V))=U\cap C_2(V).
\]

\begin{lemma} For $a,b\in V$ and $m\geq2$, define
\[
q_m(a,b)=\bigl(E(a_{-m}b)+C_2(U)\bigr)+R_+K\in Q.
\]
Then
    \[
Q=\Span_{\C}\{q_m(a,b):a,b\in V,\ m\geq2\}.
\]
\end{lemma}
\begin{proof}
For $m\geq2$, translation covariance from
Section~\ref{subsec:conventions} gives
\[
a_{-m}b=\frac{1}{(m-1)!}(D^{m-2}a)_{-2}b\in C_2(V).
\]
Thus $E(a_{-m}b)\in U\cap C_2(V)$, so $q_m(a,b)$ is well defined.
Linearity of $E$ gives
\[
U\cap C_2(V)=E(C_2(V))
=\Span_{\C}\{E(a_{-2}b):a,b\in V\}.
\]
The composite quotient map
$U\cap C_2(V)\twoheadrightarrow K\twoheadrightarrow Q$
is surjective and sends $E(a_{-2}b)$ to $q_2(a,b)$. Hence
$Q=\Span_{\C}\{q_2(a,b):a,b\in V\}$.
The family with $m\geq2$ contains these generators, and every
$q_m(a,b)$ already lies in $Q$, proving the displayed equality.
\end{proof}

\begin{lemma}\label{lem:qreducefirst}
For homogeneous $u\in U_+$, $a,b\in V$, and $m\geq2$,
\begin{equation}\label{eq:qreduce1}
q_m(u_{-1}a,b)=\sum_{j\geq0}q_{m+1+j}(a,u_jb).
\end{equation}
\end{lemma}
\begin{proof}
Since $gu=u$ for
every $g\in G$, averaging commutes with its coefficients:
\[
E(u_rw)=u_rE(w)\qquad(w\in V,\ r\in\Z).
\]
Write $[v]_R=v+C_2(U)$ for $v\in U$.
Since $m\geq2$, we have $[E(a_{-m}b)]_R\in K$, while
$[u]_R\in R_+$. The product in $R$ therefore satisfies
\[
[u_{-1}E(a_{-m}b)]_R
=[u]_R[E(a_{-m}b)]_R\in R_+K.
\]
Consequently this term has zero image in $Q=K/(R_+K)$.
For $r\geq2$, translation covariance also gives
\[
u_{-r}E(w)=\frac{1}{(r-1)!}(D^{r-2}u)_{-2}E(w)
\in C_2(U),
\]
because $D^{r-2}u,E(w)\in U$. These terms already vanish in $R$.

Taking the coefficient with index $-m$ in the iterate formula
\eqref{eq:iterate} and applying $E$ gives
\[
E((u_{-1}a)_{-m}b) ={}u_{-1}E(a_{-m}b) +\sum_{j\geq1}u_{-1-j}E(a_{-m+j}b) +\sum_{j\geq0}E(a_{-m-1-j}u_jb).
\]
The first term has zero image in $Q$ by the product calculation
above. Every term in the first sum belongs to $C_2(U)$, since
$1+j\geq2$. The terms in the last sum represent
$q_{m+1+j}(a,u_jb)$. Passing to $Q$ gives
\eqref{eq:qreduce1}. All sums are finite on the fixed arguments
by lower truncation.
\end{proof}

\begin{lemma}\label{lem:qreducesecond}
For homogeneous $u\in U_+$, $a,b\in V$, and $m\geq2$,
\begin{equation}\label{eq:qreduce2}
q_m(a,u_{-1}b)=-\sum_{j\geq0}(-1)^j
q_{m+1+j}(u_ja,b).
\end{equation}
\end{lemma}
\begin{proof}
Substitute $p=-1$ and $q=-m$ into
\eqref{eq:commutator}. Since $\binom{-1}{j}=(-1)^j$, this gives
\[
[u_{-1},a_{-m}]
=\sum_{j\geq0}(-1)^j(u_ja)_{-m-1-j}.
\]
Rearranging $a_{-m}u_{-1}=u_{-1}a_{-m}-[u_{-1},a_{-m}]$,
applying both sides to $b$, and averaging gives
\[
E(a_{-m}u_{-1}b)
=u_{-1}E(a_{-m}b)
-\sum_{j\geq0}(-1)^jE((u_ja)_{-m-1-j}b).
\]
The product calculation in the proof of Lemma~\ref{lem:qreducefirst}
shows that the first term has zero image in $Q$. Each summand has
index $-m-1-j\leq-3$ and represents
$q_{m+1+j}(u_ja,b)$. Passing to $Q$ proves \eqref{eq:qreduce2}.
The sum is finite because $u_ja=0$ for sufficiently large $j$.
\end{proof}

Recall the hypotheses of Theorem~\ref{thm:descent}.
With $C_1^U(V)$ defined by \eqref{eq:c1}, put
$d=\dim_{\C}(V/C_1^U(V))$ and choose homogeneous vectors
$s_1,\ldots,s_d\in V$ such that
$s_1+C_1^U(V),\ldots,s_d+C_1^U(V)$ form a basis of $V/C_1^U(V)$.
Set
\[
w_i=\wt(s_i)\quad(1\leq i\leq d),\qquad
h=\max_{1\leq i\leq d}w_i.
\]
Thus $h$ is the largest weight among the chosen representatives.

For $1\leq i,j\leq d$ and an integer $N\geq w_i+w_j+1$, define
\[
q_{ij}(N)=q_{N-w_i-w_j+1}(s_i,s_j)\in Q_N,
\]
where $Q_N$ is the homogeneous component of $Q$ of weight $N$.
The restriction on $N$ ensures that $N-w_i-w_j+1\geq2$, as required
in the definition of $q_m$. In particular, every $q_{ij}(N)$ is
defined when $N\geq2h+2$, the range used below.
Since $K$ and $R_+K$ are graded and $V$ is of CFT type,
\[
Q=\bigoplus_{N\geq0}Q_N.
\]

\begin{lemma}\label{lem:Qgradedpieces}
For every $N\geq0$, the space $Q_N$ is finite dimensional.
For $N\geq2h+2$, it is spanned by the vectors
$q_{ij}(N)$, $1\leq i,j\leq d$.
\end{lemma}
\begin{proof}
Each $Q_N$ is a quotient of a subspace of $R_N$, and $R_N$ is
a quotient of $U\cap V_N$. Hence
$\dim Q_N\leq\dim V_N<\infty$, including the degrees $N<2h+2$.
To prove the spanning assertion, let $a,b\in V$ be homogeneous
of weights $A,B$, respectively.
Since $C_1^U(V)$ is graded and the classes of $s_1,\ldots,s_d$
form a basis of $V/C_1^U(V)$, there exist finite decompositions
\begin{equation*}
a=\sum_{w_i=A}\alpha_i s_i +\sum_{\lambda}u^\lambda_{-1}x^\lambda,\qquad b=\sum_{w_j=B}\beta_j s_j +\sum_{\mu}v^\mu_{-1}y^\mu,
\end{equation*}
where all arguments are homogeneous,
$u^\lambda,v^\mu\in U_+$, and
\[
\wt(u^\lambda)+\wt(x^\lambda)=A,\qquad
\wt(v^\mu)+\wt(y^\mu)=B.
\]
These decompositions follow from the definition of the quotient and
of $C_1^U(V)$ as a linear span. Fix such decompositions independently
of the index parameter $m$.

By bilinearity and Lemma~\ref{lem:qreducefirst},
\begin{equation}\label{eq:qexpandfirst}
q_m(a,b)=\sum_{w_i=A}\alpha_i q_m(s_i,b)
 +\sum_{\lambda}\sum_{t\geq0}
 q_{m+1+t}(x^\lambda,u^\lambda_t b).
\end{equation}
For each term in the first sum, Lemma~\ref{lem:qreducesecond} gives
\begin{equation}\label{eq:qexpandsecond}
q_m(s_i,b)=\sum_{w_j=B}\beta_jq_m(s_i,s_j)
 -\sum_{\mu}\sum_{t\geq0}(-1)^t
 q_{m+1+t}(v^\mu_t s_i,y^\mu).
\end{equation}
Each sum over $t$ is finite by lower truncation. Every nonzero
remainder pair in these two formulas has total weight
$A+B-t-1<A+B$. Thus induction on the nonnegative integer $A+B$
reduces every remainder to pairs of representatives. The base case
$A+B=0$ is immediate, since $C_1^U(V)$ has no weight zero part.
The parameter $m$ only increases, so it remains at least two
throughout the reduction.

Every step preserves the total weight $A+B+m-1$. Consequently,
for all $m\geq2$ the reduction has the form
\begin{equation}\label{eq:qrepresentativereduction}
q_m(a,b)=\sum_{\substack{1\leq i,j\leq d\\w_i+w_j\leq A+B}}
c_{ij}(a,b)\,q_{m+A+B-w_i-w_j}(s_i,s_j).
\end{equation}
The coefficients $c_{ij}(a,b)$ depend on the chosen decompositions
but not on $m$: both the decompositions and the coefficients in
\eqref{eq:qexpandfirst}--\eqref{eq:qexpandsecond} are independent
of $m$. Since the homogeneous $q_m(a,b)$ span $Q$, this proves
\[
Q_N=\Span_{\C}\{q_{ij}(N):1\leq i,j\leq d\}
\qquad(N\geq2h+2).\qedhere
\]
\end{proof}

\begin{lemma}\label{lem:Qfinite}
The graded vector space
$Q=K/(R_+K)=\bigoplus_{N\geq0}Q_N$ is finite dimensional.
\end{lemma}
\begin{proof}
Fix $N\geq2h+2$ and $1\leq i,j\leq d$, and put
$m=N-w_i-w_j\geq2$. If $Ds_i\ne0$, its weight is $w_i+1$.
Applying \eqref{eq:qrepresentativereduction} to $(Ds_i,s_j)$ gives
\[
q_m(Ds_i,s_j) =\sum_{k,l}c_{kl}(Ds_i,s_j)\, q_{m+w_i+1+w_j-w_k-w_l}(s_k,s_l) =\sum_{k,l}c_{kl}(Ds_i,s_j)\,q_{kl}(N),
\]
where coefficients outside the range in
\eqref{eq:qrepresentativereduction} are set to zero.
These coefficients do not depend on $N$, because that reduction
is independent of $m$ and the pair $(Ds_i,s_j)$ is fixed.
Let $B_{ij,kl}=c_{kl}(Ds_i,s_j)$; when $Ds_i=0$, set
$B_{ij,kl}=0$ for all $k,l$.

On the other hand, $(Ds_i)_{-m}=m(s_i)_{-m-1}$ implies
\[
q_m(Ds_i,s_j)=m q_{m+1}(s_i,s_j)
=(N-w_i-w_j)q_{ij}(N).
\]
Combining these two expressions yields
\begin{equation}\label{eq:spectral}
(N-w_i-w_j)q_{ij}(N)
=\sum_{k,l}B_{ij,kl}q_{kl}(N).
\end{equation}
Order the pairs $(i,j)$ lexicographically and write
\[
\mathbf q(N)=(q_{ij}(N))_{(i,j)}\in Q_N^{d^2},\qquad
M=\operatorname{diag}_{(i,j)}(w_i+w_j)+B
\in\operatorname{Mat}_{d^2}(\C).
\]
The matrix $M$ is fixed, and \eqref{eq:spectral} becomes
\[
(N I_{d^2}-M)\mathbf q(N)=0.
\]
If $N$ is not an eigenvalue of $M$, multiplying by
$(N I_{d^2}-M)^{-1}$ gives $q_{ij}(N)=0$ for every $i,j$.
Since these vectors span $Q_N$ by
Lemma~\ref{lem:Qgradedpieces}, it follows that $Q_N=0$.

A finite matrix has only finitely many eigenvalues. Choose an
integer $N_0\geq2h+2$ larger than the absolute value of each
eigenvalue of $M$. Then $Q_N=0$ for every $N\geq N_0$, and hence
$Q$ is finite dimensional.
\end{proof}

\begin{proposition}\label{prop:finitegeneration}
Under the hypotheses of Theorem~\ref{thm:descent}, the graded
commutative algebra $R=U/C_2(U)$ is finitely generated.
\end{proposition}
\begin{proof}
Since $U$ is of CFT type, $R_0=\C[\one]$. The map $\iota$ sends
$[\one]$ to the nonzero class $\one+C_2(V)$, so its restriction
to $R_0$ is injective. Thus $K_0=0$. As $K$ is graded, this gives
\[
K\subseteq R_+,\qquad R_+K\subseteq R_+^2.
\]

Lemma~\ref{lem:Qfinite} shows that $Q=K/(R_+K)$ is
finite dimensional. Choose homogeneous lifts
$k_1,\ldots,k_e\in K$ of a basis of $Q$.
Also, $R/K\cong\operatorname{im}\iota\subseteq V/C_2(V)$ is
finite dimensional by $C_2$-cofiniteness of $V$.
Choose homogeneous $t_1,\ldots,t_f\in R_+$ whose classes form a
basis of $(R/K)_+$. All the $k_i$ and $t_j$ have positive degree.

Let $x\in R_N$ with $N>0$. There exists $y\in K_N$ such that
\[
x=y+\sum_{j=1}^f\beta_jt_j.
\]
Next express the class of $y$ in $K/(R_+K)$ in the basis given
by the $k_i$. Then we have
\[
x=z+\sum_{i=1}^e\alpha_i k_i+\sum_{j=1}^f\beta_jt_j,
\]
where $z\in R_{+}K$.
Since the chosen lifts are homogeneous, the coefficients can be
chosen so that $\alpha_i=0$ when $\deg k_i\ne N$ and $\beta_j=0$
when $\deg t_j\ne N$. Hence
\[
z=x-\sum_{j=1}^f\beta_jt_j-\sum_{i=1}^e\alpha_i k_i
\in(R_+K)_N\subseteq(R_+^2)_N.
\]
Thus the classes of the finitely many $k_i,t_j$ span
$R_+/R_+^2$, and
\[
\dim_{\C}(R_+/R_+^2)
\leq \dim_{\C}Q+\dim_{\C}(R/K)_+<\infty.
\]

Theorem~\ref{thm:gradedgeneration} now proves finite generation.
\end{proof}

\subsection{Scalar reduction and algebraicity of averaged coefficients}
\label{sec:productclosure}

We reduce scalar matrix coefficients to a finite differential
system, then prove algebraicity by an averaged product identity.
The reduction uses the identities of
Proposition~\ref{prop:external-identities};
we give the argument for a functional in the full algebraic dual.
The subsequent uniform denominator argument justifies summing
its homogeneous components. We use the vertex operator and
$C_1$ conventions of Section~\ref{subsec:conventions} and the
restricted dual pairing and rational expansions of
Section~\ref{subsec:rationality}.

\begin{theorem}\label{thm:scalarreduction}
Let $U\subset V$ be vertex operator algebras of CFT type sharing their conformal
vector, and suppose
$\dim V/C_1^U(V)<\infty$.
Let $\lambda:U\to\C$ and $\varphi:V\to\C$ be linear maps satisfying
\begin{equation}\label{eq:scalarhyp}
\varphi(u_{-1}w)=\lambda(u)\varphi(w),\qquad
\varphi(u_{-n}w)=0\quad(n\geq2)
\end{equation}
for $u\in U$ and $w\in V$. Choose homogeneous lifts
$s_1,\ldots,s_d$ of a basis of $V/C_1^U(V)$ and put
\[
F(a,b;z)=\varphi(Y(a,z)b)\in\C((z)).
\]
Then
\begin{equation}\label{eq:L}
\cL:=\sum_{a,b\in V}\C[z^{-1}]F(a,b;z)
=\sum_{i,j=1}^d\C[z^{-1}]F(s_i,s_j;z).
\end{equation}
Moreover, the vector $\cF=(F(s_i,s_j;z))_{i,j}$ satisfies
\begin{equation}\label{eq:ODE}
\cF'(z)=A(z)\cF(z),\qquad
A(z)\in\operatorname{Mat}_{d^2}(\C[z^{-1}]).
\end{equation}
\end{theorem}
\begin{proof}
This proof follows an approach parallel to
Section~\ref{sec:finitegeneration}.
Equations \eqref{eq:iterate}, \eqref{eq:commutator}, and
\eqref{eq:scalarhyp} give, for positive weight homogeneous $u\in U$,
\begin{align*}
F(u_{-1}a,b;z)&=\lambda(u)F(a,b;z)
+\sum_{j\geq0}z^{-j-1}F(a,u_jb;z),
\\
F(a,u_{-1}b;z)&=\lambda(u)F(a,b;z)
-\sum_{j\geq0}(-1)^jz^{-j-1}F(u_ja,b;z).
\end{align*}
The sums are finite. For the first equation, all terms with a
coefficient $u_{-r}$, $r\geq2$, of the outer vertex operator vanish,
while its coefficient $u_{-1}$ gives the scalar term. The second follows by commuting $u_{-1}$
past $Y(a,z)$.
Every pair of arguments on the right has smaller total weight than the
pair on the left. Decompose each argument into the chosen lifts and
$C_1^U(V)$, and induct on total weight to obtain \eqref{eq:L}.
Finally, translation covariance from Section~\ref{subsec:conventions}
gives $F'(a,b;z)=F(Da,b;z)$, and the same reduction of the
finitely many pairs $(Ds_i,s_j)$ gives \eqref{eq:ODE}.
\end{proof}

Use the expansion conventions of
Section~\ref{subsec:rationality}, and put
\[
\mathcal R:=\C[z_1^{\pm1},z_2^{\pm1},z_3^{\pm1}, (z_1-z_2)^{-1},(z_1-z_3)^{-1},(z_2-z_3)^{-1}].
\]
Write $\iota_{s_1,s_2,s_3}$ for the $\C$-algebra homomorphism
from $\mathcal R$, with $z_i$ replaced by $s_i$, into
$\C((s_1))((s_2))((s_3))$ that fixes Laurent monomials and satisfies
\[
\iota_{s_1,s_2,s_3}(s_i-s_j)^{-q}
=\sum_{r\ge0}\binom{q+r-1}{r}s_i^{-q-r}s_j^r
\quad(1\le i<j\le3,\ q\ge1).
\]

Here $\iota_{x;u,v}$ is the joint expansion defined in
Section~\ref{subsec:rationality}.

The next lemma applies
Proposition~\ref{prop:external-multiproduct} and
Lemma~\ref{lem:matrixassociativity} to the specified formal
expansion maps.

\begin{lemma}\label{lem:permutedexpansions}
Let $V$ be a vertex operator algebra of CFT type, let $a_1,a_2,a_3,a_4\in V$,
and let $\theta\in V'$. There is a unique rational function
$f_\theta\in\mathcal R$ satisfying the following identities.
For $\sigma\in S_3$, set
\[
f_{\theta,\sigma}(s_1,s_2,s_3)
=f_\theta\bigl(s_{\sigma^{-1}(1)},s_{\sigma^{-1}(2)},
                 s_{\sigma^{-1}(3)}\bigr).
\]
For every $\sigma\in S_3$,
\[
\iota_{s_1,s_2,s_3}f_{\theta,\sigma}(s_1,s_2,s_3) =\bigl\langle\theta, Y(a_{\sigma(1)},s_1)Y(a_{\sigma(2)},s_2) Y(a_{\sigma(3)},s_3)a_4\bigr\rangle,
\]
and
\[
\iota_{x;u,v}f_{\theta,\sigma}(x+u,x,v) =\bigl\langle\theta, Y\bigl(Y(a_{\sigma(1)},u)a_{\sigma(2)},x\bigr) Y(a_{\sigma(3)},v)a_4\bigr\rangle.
\]
\end{lemma}
\begin{proof}
Apply Proposition~\ref{prop:external-multiproduct} with
$n=3$, $(v_1,v_2,v_3,w)=(a_1,a_2,a_3,a_4)$ and functional
$\theta$, and denote its rational function by $f_\theta$.
For every $\sigma\in S_3$, the proposition identifies
\[
\iota_{z_{\sigma(1)},z_{\sigma(2)},z_{\sigma(3)}}f_\theta =\langle\theta, Y(a_{\sigma(1)},z_{\sigma(1)}) Y(a_{\sigma(2)},z_{\sigma(2)}) Y(a_{\sigma(3)},z_{\sigma(3)})a_4\rangle.
\]
Substituting $z_{\sigma(i)}=s_i$ gives the first identity.
The injectivity of $\iota_{s_1,s_2,s_3}$ makes its case
$\sigma=\mathrm{id}$ determine $f_\theta$ uniquely.

Apply Lemma~\ref{lem:matrixassociativity} with
\[
(a,b,c,d;\theta)
=(a_{\sigma(1)},a_{\sigma(2)},a_{\sigma(3)},a_4;\theta),
\qquad (z',z'')=(u,v).
\]
Its rational function is $f_{\theta,\sigma}$ by the first
identity, so \eqref{eq:matrixassociativity} gives the second identity.
\end{proof}

Let $f(x,z)\in\C((z))(x)$ be a rational function of $x$, with
coefficients in $\C((z))$. To extract coefficients at $x=0$,
take its unique formal Laurent expansion
\[
f(x,z)=\sum_{j\geq j_0}f_j(z)x^j
\quad\text{in }\C((z))((x)),\qquad f_j(z)\in\C((z)).
\]
For each integer $k$, define
\[
[x^k]_0(f(x,z)):=f_k(z),
\]
where $f_k(z)=0$ if $k<j_0$. Thus this notation denotes the
coefficient of $x^k$, with $z$ retained in the coefficient field.
The subscript $0$ specifies the expansion point, and negative powers
of $x$ are allowed when $f(x,z)$ has a pole there.
Similarly, if the Laurent expansion at infinity is
\[
f(x,z)=\sum_{j\leq j_1}f_j^{\infty}(z)x^j
\quad\text{in }\C((z))((x^{-1})),
\]
we define $[x^k]_\infty(f(x,z)):=f_k^{\infty}(z)$, taking this
coefficient to be zero for $k>j_1$. For a polynomial in $x$, both
expansions are the same finite sum, and we write $[x^k](f(x,z))$
for its coefficient of $x^k$.

The same definitions apply to a rational function in $K(t)$
for any coefficient field $K$, with $t$ as the expansion variable
and all other variables retained in $K$.

\begin{lemma}\label{lem:uniformdenominator}
Use the notation of Lemma~\ref{lem:permutedexpansions}, and assume
that each $a_i$ is homogeneous of weight $w_i$ for $1\le i\le4$. Set
\begin{equation*}
\mathcal D(z_1,z_2,z_3)=
\prod_{1\leq i<j\leq3}(z_i-z_j)^{w_i+w_j}
\prod_{i=1}^3z_i^{w_i+w_4}.
\end{equation*}
Then
\[
\mathcal D f_\theta\in\C[z_1,z_2,z_3].
\]
\end{lemma}
\begin{proof}
Put $q_{ij}=w_i+w_j$. The weight formula gives
\[
\wt((a_i)_m a_j)=q_{ij}-m-1<0
\quad(m\ge q_{ij}),\qquad (a_i)_m a_j=0.
\]
Write $f_\theta=A/B$ in reduced form, where
\[
A\in\C[z_1,z_2,z_3],\qquad
B=\prod_{i=1}^3z_i^{p_i}
  \prod_{1\le i<j\le3}(z_i-z_j)^{p_{ij}},\qquad
p_i,p_{ij}\in\N.
\]
For the coefficient extractions below, use $K=\C(x,y)$
for expansion in $t$ and $K=\C(x,v)$ for expansion in $u$.

Fix $i$ and let $(\sigma(1),\sigma(2),\sigma(3))=(j,k,i)$,
where $\{i,j,k\}=\{1,2,3\}$. The product identity of
Lemma~\ref{lem:permutedexpansions} gives, for $\ell<-q_{i4}$,
\[
\iota_{x,y}\bigl([t^\ell]_0(f_{\theta,\sigma}(x,y,t))\bigr)=[t^\ell]\iota_{x,y,t}f_{\theta,\sigma}(x,y,t) =\langle\theta,Y(a_j,x)Y(a_k,y)(a_i)_{-\ell-1}a_4\rangle =0.
\]
The map $\iota_{x,y}$ is injective. Thus the indicated Laurent
coefficients vanish, and the reduced denominator satisfies
\[
p_i\le q_{i4}=w_i+w_4\qquad(1\le i\le3).
\]

For $i<j$, take $(\sigma(1),\sigma(2),\sigma(3))=(i,j,k)$ and put
\[
g_{ij}(x,u,v)=f_{\theta,\sigma}(x+u,x,v).
\]
The second identity of Lemma~\ref{lem:permutedexpansions} gives,
for $\ell<-q_{ij}$,
\[
\iota_{x,v}\bigl([u^\ell]_0(g_{ij}(x,u,v))\bigr)=[u^\ell]\iota_{x;u,v}g_{ij}(x,u,v) =\langle\theta, Y((a_i)_{-\ell-1}a_j,x)Y(a_k,v)a_4\rangle =0.
\]
The compatibility of the coefficient expansions follows from
\[
\iota_{x;u,v}(x+u-v)^{-p} =\sum_{r\ge0}(-1)^r\binom{p+r-1}{r} u^r\iota_{x,v}(x-v)^{-p-r}\qquad(p\in\Z_{\ge1}),
\]
and its specialization $v=0$, since the denominator factors
other than $u$ are
\[
x+u,\quad x,\quad v,\quad x-v,\quad x+u-v,
\]
which are units in $\C(x,v)[[u]]$. Injectivity of $\iota_{x,v}$
and reducedness of $A/B$ therefore give
\[
p_{ij}\le q_{ij}=w_i+w_j\qquad(1\le i<j\le3).
\]
Consequently
\[
\mathcal D f_\theta
=A\prod_{i=1}^3z_i^{q_{i4}-p_i}
  \prod_{1\le i<j\le3}(z_i-z_j)^{q_{ij}-p_{ij}}
\in\C[z_1,z_2,z_3].\qedhere
\]
\end{proof}

\begin{lemma}\label{lem:homogeneousnumerator}
Retain $V,a_i,w_i,\mathcal D$ and $f_\theta$ from
Lemma~\ref{lem:uniformdenominator}, and set
\[
\tau=w_1+w_2+w_3+w_4.
\]
Let $N\ge0$ and let $\theta_N\in V_N^*$ be extended by zero
on every other weight space. Then the polynomial
\[
\mathcal N_N(z_1,z_2,z_3)
=\mathcal D(z_1,z_2,z_3)f_{\theta_N}(z_1,z_2,z_3)
\]
is zero or homogeneous of total degree $N+2\tau$.
\end{lemma}
\begin{proof}
Lemma~\ref{lem:uniformdenominator} gives
$\mathcal N_N\in\C[z_1,z_2,z_3]$.\\
A nonzero coefficient
$\langle\theta_N,(a_1)_{r_1}(a_2)_{r_2}(a_3)_{r_3}a_4\rangle$
satisfies
\[
N=\tau-r_1-r_2-r_3-3.
\]
Its monomial $z_1^{-r_1-1}z_2^{-r_2-1}z_3^{-r_3-1}$ has
total degree $N-\tau$. Moreover,
\[
\deg\mathcal D
=2(w_1+w_2+w_3)+(w_1+w_2+w_3+3w_4)=3\tau.
\]
For $t\in\C^*$, these degree calculations give
\[
\begin{aligned}
f_{\theta_N}(tz_1,tz_2,tz_3)
&=t^{N-\tau}f_{\theta_N}(z_1,z_2,z_3),\\
\mathcal D(tz_1,tz_2,tz_3)
&=t^{3\tau}\mathcal D(z_1,z_2,z_3),\\
\mathcal N_N(tz_1,tz_2,tz_3)
&=t^{N+2\tau}\mathcal N_N(z_1,z_2,z_3).
\end{aligned}
\]
Writing $\mathcal N_N=\sum_{\alpha\in\N^3}c_\alpha z^\alpha$,
where $z^\alpha=z_1^{\alpha_1}z_2^{\alpha_2}z_3^{\alpha_3}$ and
$|\alpha|=\alpha_1+\alpha_2+\alpha_3$, we obtain
\[
\bigl(t^{|\alpha|}-t^{N+2\tau}\bigr)c_\alpha=0
\quad(t\in\C^*),\qquad
c_\alpha=0\quad\text{if }|\alpha|\ne N+2\tau.
\]
Hence
\[
\mathcal N_N
=\sum_{\substack{\alpha\in\N^3\\|\alpha|=N+2\tau}}
c_\alpha z^\alpha,
\]
with the zero polynomial allowed.
\end{proof}

For the rest of this subsection, retain the hypotheses of
Theorem~\ref{thm:descent}, write $U=V^G$, and use the averaging
map $E$ of Section~\ref{subsec:conventions}. Let $\varphi:V\to\C$
be a $G$-invariant linear functional with $\varphi(\one)=1$
satisfying \eqref{eq:scalarhyp} with $\lambda=\varphi|_U$.
In particular, $\varphi\circ E=\varphi$ and
$\varphi(C_2(U))=0$. Set $F(a,b;z)=\varphi(Y(a,z)b)$ and
define $\cL$ as in \eqref{eq:L}.
Theorem~\ref{thm:scalarreduction} shows that $\cL$ is finitely
generated and that its generating vector satisfies \eqref{eq:ODE}.

The functional $\varphi$ may have infinitely many nonzero
homogeneous components. We first apply restricted dual rationality
to each component, then justify their sum using a bound independent
of the output weight. The following lemmas establish
the product formula needed for multiplicative closure.

Recall Lemmas~\ref{lem:uniformdenominator} and~\ref{lem:homogeneousnumerator}.
Throughout the following lemmas, fix homogeneous $a,b,c,d\in V$ of weights
$A_0,B_0,C_0,D_0$, respectively, and set
$\tau=A_0+B_0+C_0+D_0$.
For $N\geq0$, let $\varphi_N$ agree with $\varphi$ on $V_N$ and
vanish on every other homogeneous component. It belongs to the
restricted dual of $V$. Define
\begin{equation}\label{eq:HN}
H_N(x,z)=\frac1{|G|^2}\sum_{g,h\in G}
\left\langle\varphi_N,
Y(ga,x+z)Y(gb,x)Y(hc,z)hd\right\rangle.
\end{equation}
Here we first take the common rational function in independent
coordinates $z_1,z_2,z_3$ and then substitute
$(z_1,z_2,z_3)=(x+z,x,z)$. Thus \eqref{eq:HN} defines a rational
function; its formal expansions will be specified below. Put
\begin{equation}\label{eq:P}
P(x,z)=x^\tau(x-z)^{B_0+C_0}(x+z)^{A_0+D_0},
\qquad M=\deg_x P(x,z)=2\tau.
\end{equation}

\begin{lemma}\label{lem:averageddenominator}
For every $N\geq0$, the expression
$z^\tau P(x,z)H_N(x,z)$ is a homogeneous polynomial in $x,z$
of total degree $N+2\tau$, possibly zero. The degree $M$ is independent of $N$.
\end{lemma}
\begin{proof}
vertex operator algebra automorphisms preserve weights, so the arguments $ga,gb,hc,hd$
have weights $A_0,B_0,C_0,D_0$ for every $g,h\in G$.
Lemma~\ref{lem:uniformdenominator} therefore gives the same
denominator for all summands in \eqref{eq:HN}:
\[
\mathcal D(z_1,z_2,z_3) ={}(z_1-z_2)^{A_0+B_0}(z_1-z_3)^{A_0+C_0} (z_2-z_3)^{B_0+C_0} \cdot z_1^{A_0+D_0}z_2^{B_0+D_0}z_3^{C_0+D_0}.
\]
Hence $\deg\mathcal D=3\tau.$

To compute the degree of a matrix coefficient, consider a nonzero
term paired with $\varphi_N$ in the product expansion. If the
indices of its vertex operator coefficients are $r,s,t$, then
\[
\wt\bigl((ga)_r(gb)_s(hc)_t hd\bigr)
=\tau-r-s-t-3=N.
\]
Its scalar monomial $z_1^{-r-1}z_2^{-s-1}z_3^{-t-1}$ has total
degree $-r-s-t-3=N-\tau$. Thus the common rational function
is homogeneous of degree $N-\tau$. By Lemma~\ref{lem:homogeneousnumerator}, multiplication by
$\mathcal D(z_1,z_2,z_3)$ gives a polynomial of degree
\[
(N-\tau)+3\tau=N+2\tau,
\]
unless the polynomial is zero. The finite average over $g,h$
has the same property.

Now substitute $(z_1,z_2,z_3)=(x+z,x,z)$. Then
\[
\mathcal D(x+z,x,z)=z^\tau P(x,z).
\]
This polynomial is nonzero, so the substitution is valid in the
rational function. Substituting homogeneous linear forms in the
numerator preserves its total degree unless it becomes zero.
Hence $z^\tau P(x,z)H_N(x,z)$ has the claimed polynomial form.
The denominator depends only on the weights of the four vectors, not on $N$.
\end{proof}

\begin{lemma}\label{lem:averagedinfinity}
For every $N\geq0$, the expansion of $H_N(x,z)$ at $x=\infty$
has no strictly positive powers of $x$.
\end{lemma}
\begin{proof}
In the rational function underlying \eqref{eq:HN}, use the
coordinates $(x+z',x,z'')$ with $z',z''$ independent.
Associativity \eqref{eq:rationalassociativity}, in the form
\eqref{eq:matrixassociativity} of Lemma~\ref{lem:matrixassociativity},
and the identities
$Y(ga,z')gb=gY(a,z')b$ and $Y(hc,z'')hd=hY(c,z'')d$ show that
averaging independently over $g,h\in G$ gives the expansion
\[
\left\langle\varphi_N,
Y\bigl(EY(a,z')b,x\bigr)EY(c,z'')d\right\rangle.
\]

Put $u_{(m)}=E(a_m b)$ and $v_{(n)}=E(c_n d)$ for $m,n\in\Z$.
Both vectors belong to $U$. A term arising from the coefficient
with index $r$ of the outer vertex operator $Y(u_{(m)},x)$ is
\[
\varphi_N((u_{(m)})_r v_{(n)})
x^{-r-1}(z')^{-m-1}(z'')^{-n-1}.
\]
Since $\varphi_N$ is supported in weight $N$, a nonzero term
must satisfy
\[
N=(A_0+B_0-m-1)+(C_0+D_0-n-1)-r-1
=\tau-m-n-r-3.
\]
Nonnegative weights also imply $u_{(m)}=0$ for $m>A_0+B_0-1$
and $v_{(n)}=0$ for $n>C_0+D_0-1$. Thus, for fixed $N,r$,
the contributing pairs $(m,n)$ satisfy
\[
A_0+B_0-N-r-2\leq m\leq A_0+B_0-1,
\qquad n=\tau-N-r-3-m.
\]
There are only finitely many such pairs, and none if $r<-N-1$.
Hence the substitution $z'=z''=z$ is legitimate coefficientwise
and produces a Laurent series in $x^{-1}$. With $H_N(x,z)$
expanded at $x=\infty$, we obtain
\begin{equation}\label{eq:internal}
\begin{aligned}
H_N(x,z)
&=\left\langle\varphi_N,
Y\bigl(EY(a,z)b,x\bigr)EY(c,z)d\right\rangle\\
&=\sum_{r\in\Z}
\left(\sum_{\substack{m,n\in\Z\\m+n=\tau-N-r-3}}
\varphi_N((u_{(m)})_r v_{(n)})\right)
x^{-r-1}z^{N+r+1-\tau}.
\end{aligned}
\end{equation}

Positive powers of $x$ have $r=-s$ with $s\geq2$.
Translation covariance gives
\[
(u_{(m)})_{-s}v_{(n)}
=\frac1{(s-1)!}(D^{s-2}u_{(m)})_{-2}v_{(n)}\in C_2(U).
\]
Since $C_2(U)$ is homogeneous and $\varphi(C_2(U))=0$,
each $\varphi_N$ also annihilates $C_2(U)$.
All positive powers of $x$ in \eqref{eq:internal} therefore vanish.
\end{proof}

\begin{lemma}\label{lem:averagedsummation}
For every $N\geq0$, the product $P(x,z)H_N(x,z)$ is a polynomial
in $x$ with coefficients in $\C((z))$, of the form
\begin{equation}\label{eq:homogeneouscoeff}
P(x,z)H_N(x,z)=\sum_{j=0}^M c_{N,j}x^jz^{N+\tau-j},
\qquad c_{N,j}\in\C.
\end{equation}
The sum over $N$ is well defined coefficientwise in $z$ for each
power of $x$, and
\[
\sum_{N\geq0}P(x,z)H_N(x,z)
=\sum_{j=0}^M x^j\left(\sum_{N\geq0}c_{N,j}z^{N+\tau-j}\right)
\in\C((z))[x].
\]
Consequently,
\[
H(x,z):=\frac1{P(x,z)}\sum_{N\geq0}P(x,z)H_N(x,z)
\in\C((z))(x)
\]
is well defined, and $P(x,z)H(x,z)$ is a polynomial in $x$ of
degree at most $M$ over $\C((z))$.
\end{lemma}
\begin{proof}
By Lemma~\ref{lem:averageddenominator},
$z^\tau P(x,z)H_N(x,z)$ is a homogeneous polynomial of total
degree $N+2\tau$. Lemma~\ref{lem:averagedinfinity} shows that
the expansion of $H_N(x,z)$ at infinity has no positive powers
of $x$.
Homogeneity then forces its coefficient of $x^j$ to have the form
$c_{N,j}z^{N+2\tau-j}$. Dividing by $z^\tau$ proves
\eqref{eq:homogeneouscoeff}.

Now we justify the sum over $N$. For every fixed
$0\leq j\leq M$,
\[
\sum_{N\geq0}c_{N,j}z^{N+\tau-j}
=z^{\tau-j}\sum_{N\geq0}c_{N,j}z^N
\in\C((z)).
\]
Indeed, the exponents are bounded below by $\tau-j$, and the
coefficient of $z^t$ receives a contribution only from
$N=t-\tau+j$, if this integer is nonnegative.
Since there are only $M+1$ possible powers of $x$, these Laurent
series are the coefficients of a polynomial in $x$ of degree
at most $M$. This polynomial is precisely
$\sum_{N\geq0}P(x,z)H_N(x,z)$, with the sum interpreted as above.
Dividing by the nonzero polynomial $P(x,z)$ defines the asserted
element $H(x,z)$ of $\C((z))(x)$.
\end{proof}

\begin{lemma}\label{lem:averagedleading}
For the rational function $H(x,z)$ of
Lemma~\ref{lem:averagedsummation}, one has
\begin{equation}\label{eq:leading}
[x^M](P(x,z)H(x,z))=F(a,b;z)F(c,d;z).
\end{equation}
\end{lemma}
\begin{proof}
Set $u_{(m)}=E(a_m b)$ and $v_{(n)}=E(c_n d)$ for $m,n\in\Z$, as in the
proof of Lemma~\ref{lem:averagedinfinity}.
By Lemma~\ref{lem:averageddenominator}, $P(x,z)$ is monic of degree
$M$. Lemma~\ref{lem:averagedinfinity} shows that $H_N(x,z)$ has
no positive powers in its expansion at infinity. In the product
$P(x,z)H_N(x,z)$, a term of degree $M$ can therefore arise only
from the leading term $x^M$ of $P(x,z)$ and the constant term of
$H_N(x,z)$. Thus
\[
[x^M](P(x,z)H_N(x,z))=[x^0]_\infty(H_N(x,z)).
\]
The constant term in \eqref{eq:internal} has outer index $r=-1$,
so
\[
[x^0]_\infty(H_N(x,z))
=\sum_{\substack{m,n\in\Z\\m+n=\tau-N-2}}
\varphi_N((u_{(m)})_{-1}v_{(n)})z^{-m-n-2}.
\]
The sum is finite by the same truncation argument. Summing over
$N$ is also legitimate coefficient by coefficient: for a fixed
power of $z$, the sum $m+n$, and therefore the output weight
$N=\tau-m-n-2$, is fixed. Consequently
\[
[x^M](P(x,z)H(x,z))=\sum_{m,n\in\Z}
\varphi((u_{(m)})_{-1}v_{(n)})z^{-m-n-2}.
\]
Since $u_{(m)},v_{(n)}\in U$, the first identity in
\eqref{eq:scalarhyp}, with $\lambda=\varphi|_U$, and
$\varphi\circ E=\varphi$ give
\[
\varphi((u_{(m)})_{-1}v_{(n)})=\varphi(u_{(m)})\varphi(v_{(n)})
=\varphi(a_m b)\varphi(c_n d).
\]
The last double sum is precisely the product of the two Laurent
series $F(a,b;z)$ and $F(c,d;z)$. Each coefficient of that product
is a finite sum, because the two series have lower bounds on
their exponents. This proves \eqref{eq:leading}.
\end{proof}

\begin{lemma}\label{lem:averagedcrossedcomponent}
For every $N\geq0$ and every integer $k$, one has
\begin{equation}\label{eq:crossedcomponent}
[x^k]_0(H_N(x,z))=\frac1{|G|^2}\sum_{g,h\in G}
\sum_{m+n=-k-2}
\left\langle\varphi_N,
Y((ga)_m(hc),z)(gb)_n(hd)\right\rangle.
\end{equation}
The inner sum is finite: its nonzero terms satisfy
$m\leq A_0+C_0-1$ and $n\leq B_0+D_0-1$.
\end{lemma}
\begin{proof}
Fix $N,g,h$, and let $f_{N,g,h}(z_1,z_2,z_3)$ be the rational
function associated with
\[
\left\langle\varphi_N,
Y(ga,z_1)Y(gb,z_2)Y(hc,z_3)hd\right\rangle
\]
by Lemma~\ref{lem:permutedexpansions}.
Keep $x',x''$ independent
and substitute $(z_1,z_2,z_3)=(z+x',x'',z)$.

Apply Lemma~\ref{lem:permutedexpansions} with
\[
(a_1,a_2,a_3,a_4;\theta)=(ga,gb,hc,hd;\varphi_N),
\qquad \sigma=(2\,3),
\]
and $(x,u,v)=(z,x',x'')$. Its second identity gives
\[
\iota_{z;x',x''}f_{N,g,h}(z+x',x'',z)
=\left\langle\varphi_N,
Y(Y(ga,x')hc,z)Y(gb,x'')hd\right\rangle.
\]

The term with inner indices $m,n$ is
\[
\left\langle\varphi_N,
Y((ga)_m(hc),z)(gb)_n(hd)\right\rangle
(x')^{-m-1}(x'')^{-n-1}.
\]
Nonnegative weights imply that it vanishes unless
$m\leq A_0+C_0-1$ and $n\leq B_0+D_0-1$.
After setting $x'=x''=x$, the coefficient of $x^k$ therefore
receives contributions only from
\[
-k-B_0-D_0-1\leq m\leq A_0+C_0-1,
\qquad n=-k-2-m.
\]
There are finitely many such pairs, and none for $k<-\tau$.
This justifies the diagonal substitution coefficientwise. The
mixed factors now expand in nonnegative powers of $x/z$, so the
result is precisely the expansion at $x=0$ of
$f_{N,g,h}(z+x,x,z)$. Averaging over the finite set of $g,h$
proves \eqref{eq:crossedcomponent}.
\end{proof}

\begin{lemma}\label{lem:averagedcrossed}
For every integer $k$, the function $H(x,z)$ of
Lemma~\ref{lem:averagedsummation} satisfies
\begin{equation}\label{eq:crossedcoefficient}
[x^k]_0(H(x,z))=\frac1{|G|^2}\sum_{g,h\in G}
\sum_{m+n=-k-2}
F((ga)_m(hc),(gb)_n(hd);z).
\end{equation}
For fixed $k$, the right hand side is a finite sum of elements
of $\cL$.
\end{lemma}
\begin{proof}
We sum \eqref{eq:crossedcomponent} from
Lemma~\ref{lem:averagedcrossedcomponent} over $N$, first checking
that coefficient extraction commutes with this sum.
From
\eqref{eq:P}, the expansion of the reciprocal denominator at zero
has the form
\[
\frac1{P(x,z)}
=x^{-\tau}z^{-\tau}
\sum_{\ell\geq0}\gamma_\ell x^\ell z^{-\ell},
\qquad \gamma_\ell\in\C,\quad
\gamma_0=(-1)^{B_0+C_0}.
\]
Combining it with \eqref{eq:homogeneouscoeff} yields
\[
[x^k]_0(H_N(x,z))
=\left(\sum_{\substack{0\leq j\leq M\\k+\tau-j\geq0}}
c_{N,j}\gamma_{k+\tau-j}\right)z^{N-\tau-k}.
\]
The sum over $j$ is finite. For fixed $k$, the $z$ exponents
are bounded below by $-\tau-k$, and a given exponent determines
at most one $N$. Thus the definition of $H(x,z)$ in
Lemma~\ref{lem:averagedsummation} implies
\[
[x^k]_0(H(x,z))=\sum_{N\geq0}[x^k]_0(H_N(x,z))
\quad\text{in }\C((z)).
\]
In \eqref{eq:crossedcomponent}, the two vectors in the vertex operator product have
weights whose sum is
$\tau-m-n-2=\tau+k$. Pairing their product with $\varphi_N$
therefore selects precisely the coefficient with power
$z^{N-\tau-k}$. Summing over $N$ recovers the full series
$F((ga)_m(hc),(gb)_n(hd);z)$. Since there are only finitely many
$g,h,m,n$ for the fixed $k$, this summation proves
\eqref{eq:crossedcoefficient}. Each summand belongs to $\cL$
by its definition in \eqref{eq:L}.
\end{proof}

\begin{proposition}\label{prop:productclosure}
The finite dimensional $\C(z)$-space $\C(z)\cL$ is closed under
multiplication.
\end{proposition}
\begin{proof}
By bilinearity, it suffices to prove the assertion for
$F(a,b;z)F(c,d;z)$ with homogeneous $a,b,c,d$. Retain the notation
of the preceding lemmas and write
$P(x,z)=\sum_{i=0}^M p_i(z)x^i$, where $p_i(z)\in\C[z]$.
Since $P(x,z)H(x,z)$ is a polynomial in $x$, its coefficient of $x^M$ can
be computed using the expansion of $H(x,z)$ at zero. Hence
Lemma~\ref{lem:averagedleading} gives
\[
F(a,b;z)F(c,d;z)=[x^M](P(x,z)H(x,z))
=\sum_{i=0}^M p_i(z)[x^{M-i}]_0(H(x,z)).
\]
Applying Lemma~\ref{lem:averagedcrossed} to $k=M-i$ yields
\begin{equation*}
F(a,b;z)F(c,d;z) ={}\frac1{|G|^2}\sum_{g,h\in G}\sum_{i=0}^M p_i(z) \mathrel{}\cdot\sum_{m+n=i-M-2} F((ga)_m(hc),(gb)_n(hd);z).
\end{equation*}
Every sum here is finite, and its terms belong to
$\C(z)\cL$. These products therefore lie in $\C(z)\cL$.
Every element of that space is a finite $\C(z)$-linear combination
of the series $F(a,b;z)$, so distributivity proves closure for
arbitrary products.
\end{proof}

\begin{proposition}\label{prop:descentfield}
Let $d=\dim V/C_1^U(V)$ and set $\cA=\C(z)\cL$.
Then $\cA\subseteq\C((z))$ is a finite field extension of
$\C(z)$, with $\dim_{\C(z)}\cA\leq d^2$. In particular, every
$F(a,b;z)$ is algebraic over $\C(z)$.
\end{proposition}
\begin{proof}
Choose homogeneous lifts $s_1,\ldots,s_d$ of a basis of
$V/C_1^U(V)$. The functional $\varphi$ satisfies the hypotheses of
Theorem~\ref{thm:scalarreduction} with $\lambda=\varphi|_U$, so
\[
\cA=\sum_{i,j=1}^d\C(z)F(s_i,s_j;z)
\subseteq\C((z)),\qquad
\dim_{\C(z)}\cA\leq d^2.
\]
Proposition~\ref{prop:productclosure} shows that $\cA$ is closed
under multiplication. Moreover,
\[
F(\one,\one;z)=\varphi(\one)=1,
\]
so $\cA$ contains the unit and therefore all of $\C(z)$.
For any nonzero $f\in\cA$, multiplication by $f$ is an injective
$\C(z)$-linear map $\cA\to\cA$, because $\C((z))$ is a field. Finite dimensionality makes this map
surjective, so there is $g\in\cA$ with $fg=1$. Thus $\cA$ is a
finite field extension of $\C(z)$, proving both assertions.
\end{proof}

\subsection{Proof of Theorem~\ref*{thm:descent}}\label{sec:descent}

Retain the hypotheses of Theorem~\ref{thm:descent} and the notation
$U=V^G$, $R=U/C_2(U)$, and $E$, with the averaging and quotient
conventions of Section~\ref{subsec:conventions}. Let $\chi:R\to\C$ be an
arbitrary unital character, also denoting its pullback to $U$ by
$\chi$, and set $\varphi=\chi\circ E$.
Then $\varphi$ is $G$-invariant, $\varphi(\one)=1$, and
$\varphi|_U=\chi$. Averaging is a $U$-module map, so
\begin{equation*}
\varphi(u_{-1}w)=\chi(u)\varphi(w),\qquad
\varphi(u_{-n}w)=0\quad(n\geq2)
\end{equation*}
for $u\in U$ and $w\in V$. Indeed, average first and then use
multiplication in $R$ and $u_{-n}Ew\in C_2(U)$ for $n\geq2$.
Thus $\varphi$ satisfies the conditions used in
Section~\ref{sec:productclosure}. The character need not preserve
the grading, so $\varphi$ may initially have infinitely many
nonzero homogeneous components.

Choose homogeneous lifts $s_1,\ldots,s_d$ of a basis of
$V/C_1^U(V)$; the hypothesis of the theorem ensures that $d$ is
finite. Set $F(a,b;z)=\varphi(Y(a,z)b)$, define $\cL$ as in
\eqref{eq:L}, and put $\cA=\C(z)\cL$.
Proposition~\ref{prop:finitegeneration} already shows that $R$ is
finitely generated. By Theorem~\ref{lem:gradedalgebra}, it remains
to show that $\chi(R_+)=0$. We establish this through the following
three lemmas.

\begin{lemma}
\label{lem:descentlaurent}
There is an integer $B\geq0$ such that for all $a,b\in V$,
\[
F(a,b;z)\in z^B\C[z^{-1}].
\]
Consequently, for $q=B+2\geq2$, one has $\varphi(C_q(V))=0$.
\end{lemma}
\begin{proof}
Theorem~\ref{thm:scalarreduction} supplies the differential system
\[
\cF'(z)=A(z)\cF(z),\qquad
\cF(z)=(F(s_i,s_j;z))_{i,j},\qquad
A(z)\in\operatorname{Mat}_{d^2}(\C[z^{-1}]).
\]
By Proposition~\ref{prop:descentfield}, every component is algebraic
over $\C(z)$. Each component belongs to $\C((z))$, and the entries
of $A(z)$ are Laurent polynomials. Thus all hypotheses of
Lemma~\ref{lem:algebraicODE} are satisfied, and every
$F(s_i,s_j;z)$ is a Laurent polynomial.

There are only finitely many such Laurent polynomials. Choose an
integer $B\geq0$ so that none contains a power of $z$ greater
than $B$. For arbitrary $a,b\in V$, \eqref{eq:L} gives
\[
F(a,b;z)=\sum_{i,j=1}^d
p_{ij}^{a,b}(z^{-1})F(s_i,s_j;z),\qquad
p_{ij}^{a,b}(z^{-1})\in\C[z^{-1}].
\]
Multiplying by a polynomial in $z^{-1}$ cannot increase the
largest exponent of $z$. Hence the same bound $B$ applies to
every $F(a,b;z)$. When $n\geq B+2$, comparing coefficients in
$F(a,b;z)=\sum_{m\in\Z}\varphi(a_m b)z^{-m-1}$ gives
\[
\varphi(a_{-n}b)=[z^{n-1}]F(a,b;z) =0.
\]
Taking $q=B+2$, the definition of $C_q(V)$ in
Section~\ref{subsec:conventions} proves the last assertion.
\end{proof}

\begin{lemma}
\label{lem:descentweightcutoff}
There is an integer $N_\chi\geq0$ such that
\[
\varphi(V_m)=0,\qquad \chi(R_m)=0
\]
for all $m>N_\chi$.
\end{lemma}
\begin{proof}
Choose $q\geq2$ as in Lemma~\ref{lem:descentlaurent}.
Since $V$ is $C_2$-cofinite, Lemma~\ref{lem:cq} gives
$\dim V/C_q(V)<\infty$.
We can therefore choose an integer
$N_\chi\geq0$ such that
\[
V_m\subseteq C_q(V)\qquad(m>N_\chi).
\]
Since $\varphi(C_q(V))=0$, this gives $\varphi(V_m)=0$.
Every class in $R_m$ has a representative
in $U\cap V_m$. Note that
$\varphi|_U=\chi$, and therefore $\chi(R_m)=0$ for $m>N_\chi$.
\end{proof}

\begin{lemma}
\label{lem:descentvanishing}
Every unital character $R\to\C$ vanishes on $R_+$.
\end{lemma}
\begin{proof}
Apply Lemma~\ref{lem:descentweightcutoff} to the arbitrary
character $\chi$ fixed above. Let $u\in U$ be homogeneous of
positive weight $w$ and write $[u]=u+C_2(U)\in R_w$.
Multiplication in $R$ respects weights, so $[u]^r\in R_{rw}$
for every $r\geq1$. Choose $r$ with $rw>N_\chi$. Then
\[
\chi([u])^r=\chi([u]^r) =0,
\]
and hence $\chi([u])=0$. Every element of $R_+$ is a finite sum
of homogeneous classes of this kind, so $\chi(R_+)=0$.
\end{proof}

\begin{proof}[Proof of Theorem~\ref{thm:descent}]
Proposition~\ref{prop:finitegeneration} and
Lemma~\ref{lem:descentvanishing} verify the hypotheses of
Theorem~\ref{lem:gradedalgebra}, which proves the
$C_2$-cofiniteness of $U$.
\end{proof}

\section{Permutation subgroups}
\label{sec:permutationsubgroups}\label{subsec:permutationsubgroups}

For a simple vertex operator algebra of CFT type, cofiniteness of the symmetric group
fixed point algebra implies cofiniteness for every permutation subgroup.
The argument uses Theorem~\ref{thm:external-galois} and
Proposition~\ref{prop:external-module-c2}.

\begin{theorem}\label{prop:permutationsubgroups}
Let $V$ be a simple vertex operator algebra of CFT type and let $n\ge1$.
If $U^{(n)}$ is $C_2$-cofinite, then $T^n(V)^G$
is $C_2$-cofinite for every subgroup $G\le S_n$.
\end{theorem}
\begin{proof}
Put $X=T^n(V)$, $A=X^{S_n}$ and $B=X^G$.
Then $A\subseteq B\subseteq X$, with a common conformal
vector. Simplicity of $X$ follows from
Proposition~\ref{prop:tensorsimple}. Apply Theorem~\ref{thm:external-galois} to the
finite image of $S_n$ in $\operatorname{Aut}(X)$:
\[
X=\bigoplus_{\chi}M_\chi\otimes X_\chi,
\qquad
B=\bigoplus_{\chi}(M_\chi)^G\otimes X_\chi.
\]
The second equality follows because $G$ acts on the first
factor and trivially on the second. Thus $B$ is a finite direct
sum of irreducible $A$-modules. Since $A$ is $C_2$-cofinite,
Proposition~\ref{prop:external-module-c2} gives $\dim X_\chi/C_2^A(X_\chi)<\infty$.
It follows that
\[
\dim B/C_2^A(B)
=\sum_\chi\dim(M_\chi)^G\,
\dim X_\chi/C_2^A(X_\chi)<\infty.
\]
Finally, $a,b\in B$ for every $a\in A$ and $b\in B$,
so $C_2^A(B)\subseteq C_2(B)$ and
\[
\dim B/C_2(B)\le\dim B/C_2^A(B)<\infty.
\]
This proves the assertion, including nonfaithful permutation
actions by the use of the finite image.
\end{proof}

The proofs below establish cofiniteness for the symmetric
permutation orbifolds $U^{(n)}$. Applying this theorem in
each tensor degree gives the corresponding conclusions for
$T^n(V)^G$ with $G\le S_n$.

\Needspace{16\baselineskip}
\section{Theorem~\ref*{thm:criterion}: The Virasoro criterion}\label{sec:conformalcriterion}

We develop the quotient descriptions and translation bounds needed for
Theorem~\ref{thm:criterion}. The arguments involving symmetric operations and quotients begin with a vertex operator algebra $V$ of CFT type;
simplicity and self contragredience are imposed only in the statements
that use them.

Section~\ref{sec:symmetricoperations} studies the symmetric operations,
and Section~\ref{sec:symmetricquotients} identifies the two quotient
spaces. Section~\ref{subsec:gradedtools} obtains the weight one
translation bounds through Heisenberg modules. These results are
combined in Section~\ref{sec:criterionproof} to prove the criterion.

\subsection{Symmetric operations and diagonal values}
\label{sec:symmetricoperations}

Let $V$ be a vertex operator algebra of CFT type. We use $D$,
$DV$, and $V_+$ as defined in Section~\ref{subsec:conventions},
and the product $a*b=a_{-1}b$. All products retain their indicated bracketing.
The goal of this section is Proposition~\ref{prop:diagonalspan}:
modulo $DV$, every value of
$B_{k-1}(a_1,\ldots,a_{k-1};a_k)$ with all arguments in $V_+$
is a finite linear combination of values with all arguments equal.
We first establish symmetry of the operations $B_k$ and the
stability and functoriality of their spans $T_k(V)$, and then
prove full symmetry modulo $DV$ and the explicit formula for
these linear combinations. Section~\ref{sec:symmetricquotients}
uses this result to compare the spans with powers and describe
the quotients needed for Theorem~\ref{thm:criterion}.
Throughout, $a*b=a_{-1}b$ satisfies the identities in
Proposition~\ref{prop:external-product-identities}.

\begin{definition}\label{def:symmetricoperations}
Let $V$ be a vertex operator algebra of CFT type. In $B_k(a_1,\ldots,a_k;b)$,
the first $k$ arguments belong to $V_+$ and the last argument belongs
to $V$. Set $B_0(;b)=b$ and define
\begin{equation}
\label{eq:B}
B_{k+1}(a_0,\ldots,a_k;b) ={}-B_k(a_1,\ldots,a_k;a_0*b) -\sum_{j=1}^k B_k(a_1,\ldots,a_0*a_j,\ldots,a_k;b).
\end{equation}
The sum is empty for $k=0$. In particular,
\[
B_1(a;b)=-a*b,\qquad
B_2(a_1,a_2;b)=a_2*(a_1*b)+(a_1*a_2)*b.
\]
For $k\geq1$, put
\[
T_k(V)=\Span_{\C}\{B_k(a_1,\ldots,a_k;b):
 a_1,\ldots,a_k\in V_+,\ b\in V\}\subset V_+.
\]
Write $L_a(b)=a*b$.
\end{definition}

\begin{lemma}\label{lem:Ldelta}
For $a\in V_+$, use the symmetric algebra of
Section~\ref{subsec:permutationbackground} and define an operator
on $\Sym(V_+)\otimes V$ by
\begin{equation}
\label{eq:deltadefinition}
\delta_a((c_1\cdots c_k)\otimes v) ={}\sum_{j=1}^k (c_1\cdots(a*c_j)\cdots c_k)\otimes v +(c_1\cdots c_k)\otimes(a*v).
\end{equation}
It preserves symmetric degree and satisfies
\begin{equation}\label{eq:Ldeltacommutators}
[L_a,L_b]=L_{a*b-b*a},\qquad
[\delta_a,\delta_b]=\delta_{a*b-b*a}.
\end{equation}
For $b\in V_+$ and $S\in\Sym(V_+)\otimes V$, let $bS$
denote multiplication by $b$ in the symmetric factor. Then
\begin{equation}\label{eq:deltaproduct}
\delta_a(bS)=(a*b)S+b\delta_aS.
\end{equation}
\end{lemma}
\begin{proof}
The first commutator identity follows by subtracting
\eqref{eq:product-associator} with $a,b$ interchanged from
the same identity in its original order.
In the second, actions on different tensor
factors commute; the actions on each individual factor have
commutator $L_{a*b-b*a}$. This proves the second identity.
The defining sum is symmetric in $c_1,\ldots,c_k$, so the
operator is well defined.

To prove \eqref{eq:deltaproduct}, fix $a,b\in V_+$ and take
\[
S=(c_1\cdots c_k)\otimes v,
\qquad c_1,\ldots,c_k\in V_+,\quad v\in V,\quad k\geq0.
\]
Here $c_1\cdots c_k=1\in\Sym^0(V_+)$ when $k=0$.
By definition of multiplication in the symmetric factor,
\[
bS=(bc_1\cdots c_k)\otimes v,
\qquad
(a*b)S=((a*b)c_1\cdots c_k)\otimes v.
\]
Applying \eqref{eq:deltadefinition} to $bS$ gives
\[
\begin{aligned}
\delta_a(bS) &=((a*b)c_1\cdots c_k)\otimes v +\sum_{j=1}^k (bc_1\cdots(a*c_j)\cdots c_k)\otimes v +(bc_1\cdots c_k)\otimes(a*v)\\
&=(a*b)S+b\delta_aS.
\end{aligned}
\]
The first term is the action on $b$; the remaining terms are
$b\delta_aS$ by \eqref{eq:deltadefinition} applied to $S$.
The same calculation holds for $k=0$, with an empty sum.
Since these tensors span $\Sym(V_+)\otimes V$, linearity
proves \eqref{eq:deltaproduct} for every $S$.
\end{proof}

\begin{lemma}\label{lem:Bsymmetry}
The operation $B_k$ is symmetric in its first $k$ arguments.
\end{lemma}
\begin{proof}
Assume this for $B_j$, $0\leq j\leq k$, where $k\geq1$;
the initial cases $B_0,B_1$ are immediate. Write
\[
r_j((a_1\cdots a_j)\otimes b)=B_j(a_1,\ldots,a_j;b).
\]
The recurrence \eqref{eq:B} becomes
\begin{equation}\label{eq:rsymmetricrecurrence}
r_j(aS)=-r_{j-1}(\delta_aS)
\qquad(1\leq j\leq k).
\end{equation}
With its first argument fixed, $B_{k+1}$ is already symmetric in
the other $k$ arguments. To exchange its first two arguments $a,b$,
take $S\in\Sym^{k-1}(V_+)\otimes V$ and set $c=a*b-b*a$.
The difference of the two values is
\[
-r_k(\delta_a(bS))+r_k(\delta_b(aS))=-r_k(cS)-r_k(b\delta_aS)+r_k(a\delta_bS) =r_{k-1}\bigl(\delta_cS-[\delta_a,\delta_b]S\bigr) =0.
\]
Here we used \eqref{eq:deltaproduct}, then
\eqref{eq:rsymmetricrecurrence} and
\eqref{eq:Ldeltacommutators}. These exchanges generate all
permutations, completing the induction.
\end{proof}

\begin{lemma}\label{lem:Bderivation}
Let $\partial:V\to V$ be a derivation of the product
$u*v=u_{-1}v$. This means that $\partial$ is $\C$-linear and
satisfies the Leibniz rule
\[
\partial(u*v)=(\partial u)*v+u*(\partial v)
\qquad(u,v\in V).
\]
Assume also that $\partial(V_+)\subseteq V_+$. Then, for every
integer $k\geq0$, $a_1,\ldots,a_k\in V_+$, and $b\in V$,
\begin{equation}
\label{eq:Bderivation}
\partial B_k(a_1,\ldots,a_k;b) ={}\sum_{i=1}^k B_k(a_1,\ldots,\partial a_i,\ldots,a_k;b) +B_k(a_1,\ldots,a_k;\partial b).
\end{equation}
Consequently $\partial T_k(V)\subseteq T_k(V)$ for every $k\geq1$.
\end{lemma}
\begin{proof}
We prove \eqref{eq:Bderivation} by induction on $k$. The case
$k=0$ is
\[
\partial B_0(;b)=\partial b=B_0(;\partial b).
\]
Assume the identity holds for $k$, and fix
$a_0,\ldots,a_k\in V_+$ and $b\in V$. The Leibniz rule gives
\[
\partial(a_0*b)=(\partial a_0)*b+a_0*(\partial b),\qquad \partial(a_0*a_j)=(\partial a_0)*a_j+a_0*(\partial a_j) \qquad(1\leq j\leq k).
\]
Apply $\partial$ to \eqref{eq:B} and use the induction
hypothesis on each $B_k$ value. Then we have
\[
\begin{alignedat}{1}
\partial B_{k+1}(a_0,\ldots,a_k;b)&=-\Bigl[B_k(a_1,\ldots,a_k;(\partial a_0)*b) +\sum_{j=1}^k B_k(a_1,\ldots,(\partial a_0)*a_j,\ldots,a_k;b)\Bigr]\\
&\mathrel{\phantom{=}}-\sum_{i=1}^k\Bigl[ B_k(a_1,\ldots,\partial a_i,\ldots,a_k;a_0*b)\\
&\mathrel{\phantom{=}}+B_k(a_1,\ldots,a_0*(\partial a_i),\ldots,a_k;b)\\
&\mathrel{\phantom{=}}+\sum_{j=1}^{i-1} B_k(a_1,\ldots,a_0*a_j,\ldots,\partial a_i,\ldots,a_k;b)\\
&\mathrel{\phantom{=}}+\sum_{j=i+1}^k B_k(a_1,\ldots,\partial a_i,\ldots,a_0*a_j,\ldots,a_k;b)\Bigr]\\
&\mathrel{\phantom{=}}-\Bigl[B_k(a_1,\ldots,a_k;a_0*(\partial b)) +\sum_{j=1}^k B_k(a_1,\ldots,a_0*a_j,\ldots,a_k;\partial b)\Bigr]\\
&=B_{k+1}(\partial a_0,a_1,\ldots,a_k;b) +\sum_{i=1}^k B_{k+1}(a_0,a_1,\ldots,\partial a_i,\ldots,a_k;b)\\
&\mathrel{\phantom{=}}+B_{k+1}(a_0,\ldots,a_k;\partial b).
\end{alignedat}
\]
This proves \eqref{eq:Bderivation} for $k+1$.

For $k\geq1$, the hypotheses $\partial a_i\in V_+$ and
$\partial b\in V$ show that every term on the right hand side
of \eqref{eq:Bderivation} belongs to $T_k(V)$. Therefore
\[
\partial T_k(V) =\Span_{\C}\{\partial B_k(a_1,\ldots,a_k;b): a_1,\ldots,a_k\in V_+,\ b\in V\} \subseteq T_k(V).\qedhere
\]
\end{proof}

\begin{lemma}\label{lem:Tstability}
For $k\geq1$, $T_k(V)$ is preserved by $D$ and by the
coefficient $v_0$ of every vertex operator $Y(v,z)$, and
$T_{k+1}(V)\subseteq T_k(V)$.
\end{lemma}
\begin{proof}
The operators $D$ and $v_0$ are derivations of $*$ by
\eqref{eq:productderivations}. They preserve $V_+$: $D$
raises weight by one, while $v_0V_m\subseteq V_{d+m-1}$
for $v\in V_d$, $d,m>0$; a weight zero $v$ has $v_0=0$.
Lemma~\ref{lem:Bderivation} gives the two stability assertions.
Every term in the recurrence for $B_{k+1}$ is a $B_k$ value
with positive first arguments, giving the last inclusion.
\end{proof}

For a graded vertex algebra homomorphism $f:V\to\widetilde V$,
the same recurrence gives
\begin{equation}
\label{eq:functoriality}
f(B_k^V(a_1,\ldots,a_k;b)) =B_k^{\widetilde V}(f(a_1),\ldots,f(a_k);f(b)),\qquad f(T_k(V))\subseteq T_k(\widetilde V).
\end{equation}
Also $f(Da)=Df(a)$.

All congruences below are in the vector space $V/DV$.

\begin{proposition}\label{prop:fullsymmetry}
For every $k\geq1$, the class
\[
B_{k-1}(a_1,\ldots,a_{k-1};a_k)+DV,
\qquad a_1,\ldots,a_k\in V_+,
\]
is symmetric in all $k$ arguments. Moreover,
\begin{equation}\label{eq:Bvacuum}
B_k(a_1,\ldots,a_k;\one)
\equiv-kB_{k-1}(a_1,\ldots,a_{k-1};a_k)\pmod{DV}.
\end{equation}
\end{proposition}
\begin{proof}
Let $\pi:V\to V/DV$ be the quotient map and abbreviate
\[
F_k(a_1,\ldots,a_k)
=\pi(B_{k-1}(a_1,\ldots,a_{k-1};a_k)).
\]
We prove both assertions together by induction. For $k=1$,
they follow from $B_1(a;\one)=-a$. For $k=2$,
$F_2(a,b)=-\pi(a*b)$ is symmetric by \eqref{eq:trace}.
The vacuum relation for $k=2$ follows from
\[
\pi(B_2(a,b;\one))=\pi(b*a+a*b)
=2\pi(a*b)=-2F_2(a,b).
\]

Let $k\geq3$ and assume both assertions for $k-1$.
Fix $a_1,\ldots,a_k\in V_+$ and put
\[
P=\pi(B_k(a_1,\ldots,a_k;\one)),\qquad
C_i=F_k(a_1,\ldots,a_{i-1},a_{i+1},\ldots,a_k,a_i).
\]
By Lemma~\ref{lem:Bsymmetry}, the order of the other arguments
in $C_i$ does not matter. For distinct $i,j\in\{1,\ldots,k\}$, write
\[
\{r_1<\cdots<r_{k-2}\}
=\{1,\ldots,k\}\setminus\{i,j\}.
\]
Thus the remaining arguments, in this order, are
$a_{r_1},\ldots,a_{r_{k-2}}$. By symmetry in the first arguments,
\[
C_i=\pi\bigl(B_{k-1}(a_{r_1},\ldots,a_{r_{k-2}},a_j;a_i)\bigr),\qquad C_j=\pi\bigl(B_{k-1}(a_i,a_{r_1},\ldots,a_{r_{k-2}};a_j)\bigr).
\]
Apply \eqref{eq:B} to $P$ with first argument $a_i$. Using
$a_i*\one=a_i$ and the induction hypothesis
\eqref{eq:Bvacuum} in degree $k-1$ gives
\[
\begin{alignedat}{1}
P&=\pi\bigl(B_k(a_i,a_{r_1},\ldots,a_{r_{k-2}},a_j;\one)\bigr)\\
&=-\pi\bigl(B_{k-1}(a_{r_1},\ldots,a_{r_{k-2}},a_j;a_i)\bigr) -\pi\bigl(B_{k-1}(a_{r_1},\ldots,a_{r_{k-2}},a_i*a_j;\one)\bigr)\\
&\mathrel{\phantom{=}}-\sum_{s=1}^{k-2}\pi\bigl( B_{k-1}(a_{r_1},\ldots,a_i*a_{r_s},\ldots,a_{r_{k-2}},a_j;\one) \bigr)\\
&=-C_i+(k-1)F_{k-1}(a_{r_1},\ldots,a_{r_{k-2}},a_i*a_j)\\
&\mathrel{\phantom{=}}+(k-1)\sum_{s=1}^{k-2} F_{k-1}(a_{r_1},\ldots,a_i*a_{r_s},\ldots,a_{r_{k-2}},a_j).
\end{alignedat}
\]
In each summand indexed by $s$, the entry $a_{r_s}$ is replaced
by $a_i*a_{r_s}$; the other entries retain the displayed order.
The induction hypothesis applies because every argument, including
each such product, belongs to $V_+$.

Next apply \eqref{eq:B} to the displayed expression for $C_j$,
again with first argument $a_i$ but with last argument $a_j$:
\[
\begin{aligned}
C_j &=\pi\bigl(B_{k-1}(a_i,a_{r_1},\ldots,a_{r_{k-2}};a_j)\bigr)\\
&=-\pi\bigl(B_{k-2}(a_{r_1},\ldots,a_{r_{k-2}};a_i*a_j)\bigr) -\sum_{s=1}^{k-2}\pi\bigl( B_{k-2}(a_{r_1},\ldots,a_i*a_{r_s},\ldots,a_{r_{k-2}};a_j) \bigr)\\
&=-F_{k-1}(a_{r_1},\ldots,a_{r_{k-2}},a_i*a_j) -\sum_{s=1}^{k-2} F_{k-1}(a_{r_1},\ldots,a_i*a_{r_s},\ldots,a_{r_{k-2}},a_j).
\end{aligned}
\]
The last equality is the definition of $F_{k-1}$. Comparing
the two expansions, and then interchanging $i$ and $j$, yields
\[
P=-C_i-(k-1)C_j,\qquad P=-C_j-(k-1)C_i.
\]
Subtracting these equalities gives
\[
0 =\bigl(-C_i-(k-1)C_j\bigr) -\bigl(-C_j-(k-1)C_i\bigr) =(k-2)(C_i-C_j).
\]
Since $k\geq3$, it follows that $C_i=C_j$ for all $i,j$.
For every $\sigma\in S_k$, symmetry in the first $k-1$ arguments
therefore gives
\[
F_k(a_{\sigma(1)},\ldots,a_{\sigma(k)}) =C_{\sigma(k)} =C_k=F_k(a_1,\ldots,a_k).
\]
This proves full symmetry. Finally,
\[
\begin{aligned}
\pi\bigl(B_k(a_1,\ldots,a_k;\one)\bigr)&=P\\
&=-C_i-(k-1)C_j\\
&=-kC_k\\
&=-kF_k(a_1,\ldots,a_k)\\
&=-k\,\pi\bigl(B_{k-1}(a_1,\ldots,a_{k-1};a_k)\bigr),
\end{aligned}
\]
which is \eqref{eq:Bvacuum} and completes the induction.
\end{proof}

\begin{proposition}\label{prop:diagonalspan}
Let $\pi:V\to V/DV$ be the quotient map. For $k\geq1$, define
\[
F_k(a_1,\ldots,a_k)
=\pi\bigl(B_{k-1}(a_1,\ldots,a_{k-1};a_k)\bigr),
\qquad a_1,\ldots,a_k\in V_+.
\]
For $I=\{1,\ldots,k\}$ and $J\subseteq I$, put
$x_J=\sum_{j\in J}a_j$, with $x_{\varnothing}=0$. Then
\[
F_k(a_1,\ldots,a_k)
=\frac1{k!}\sum_{J\subseteq I}(-1)^{k-|J|}F_k(x_J,\ldots,x_J).
\]
Consequently,
\[
\Span\{F_k(a_1,\ldots,a_k):a_1,\ldots,a_k\in V_+\} =\Span\{F_k(x,\ldots,x):x\in V_+\}.
\]
\end{proposition}
\begin{proof}
By Proposition~\ref{prop:fullsymmetry}, $F_k$ is symmetric
in all arguments; it is also multilinear by the definition of $B_{k-1}$.
Multilinearity gives
\[
\sum_{J\subseteq I}(-1)^{k-|J|}F_k(x_J,\ldots,x_J) =\sum_{i_1,\ldots,i_k\in I} \left(\sum_{\{i_1,\ldots,i_k\}\subseteq J\subseteq I} (-1)^{k-|J|}\right) F_k(a_{i_1},\ldots,a_{i_k}).
\]
For $K=\{i_1,\ldots,i_k\}$, the inner coefficient is
\[
\sum_{K\subseteq J\subseteq I}(-1)^{k-|J|} =\sum_{r=0}^{k-|K|}\binom{k-|K|}{r}(-1)^{k-|K|-r} =\begin{cases} 0,&K\ne I,\\ 1,&K=I. \end{cases}
\]
Thus only the tuples $(\sigma(1),\ldots,\sigma(k))$ with
$\sigma\in S_k$ remain. By symmetry,
\[
\sum_{J\subseteq I}(-1)^{k-|J|}F_k(x_J,\ldots,x_J)
=\sum_{\sigma\in S_k}F_k(a_{\sigma(1)},\ldots,a_{\sigma(k)})
=k!\,F_k(a_1,\ldots,a_k).
\]
Dividing by $k!\ne0$ proves the asserted formula. Every $x_J$
belongs to $V_+$, so this formula expresses each value of $F_k$
as a finite linear combination of its diagonal values.
The term indexed by the empty set is zero because $k\geq1$. Conversely,
every diagonal value is itself a value of $F_k$, proving
the equality of spans.
\end{proof}

\subsection{From tensor product quotients to quotients of \texorpdfstring{$V$}{V}}
\label{sec:symmetricquotients}

Let $V$ be a vertex operator algebra of CFT type. Retain $V_+$,
$D$, the product $a*b=a_{-1}b$, and the spaces $T_k(V)$ from
Section~\ref{sec:symmetricoperations}. For $x\in V_+$, use the
left associated powers
\[
x^{[0]}=\one,\qquad x^{[k+1]}=x^{[k]}*x,
\qquad
\cT_k(V)=\Span_{\C}\{x^{[k]}:x\in V_+\}\quad(k\geq1).
\]
Write $R_a(b)=b*a$. Operator powers mean composition,
so $R_x^k\one=x^{[k]}$.

The first three lemmas concern arbitrary base algebras $U$ and
their modules. We then evaluate diagonal values modulo $DV$ and
prove $DV+T_k(V)=DV+\cT_k(V)$ using
Proposition~\ref{prop:diagonalspan}. Together with Abe's multiplication
formula \eqref{eq:Phimultiply}, these facts identify the two quotient
spaces used in Theorem~\ref{thm:criterion}. The first concerns
relative $C_1$ finiteness; the second also removes derivatives and
the vacuum. Both inverse maps are constructed by reducing symmetric
degree, with Proposition~\ref{prop:fullsymmetry} used in the second case.

\begin{lemma}\label{lem:stronggeneration}
Let a vertex operator algebra $U$ of CFT type be strongly
generated by a homogeneous subspace
$S\subset U_+$ that is closed under derivatives.
For every nonnegatively graded $U$-module $M$,
\[
\Crel{U}{M}=\Span\{s_{-1}m:s\in S,\ m\in M\}.
\]
\end{lemma}
\begin{proof}
Put $A_S(M)=\Span\{s_{-1}m:s\in S,\ m\in M\}$.
The inclusion $A_S(M)\subseteq C_1^U(M)$ follows from
$S\subset U_+$. For $i\ge0$, translation covariance gives
\[
(D^is)_{-1}=i!s_{-1-i},\qquad
s_{-1-i}m=\frac1{i!}(D^is)_{-1}m\in A_S(M),
\]
because $D^is\in S$. We prove the converse in the stronger
form $u_{-m}M\subseteq A_S(M)$ for all homogeneous
$u\in U_+$ and all $m\ge1$, by induction on $\wt(u)$.
Strong generation and the derivative conversion just proved
write $u$ as a finite sum of vectors $t_{-1}a$ with
$t\in S$ homogeneous and $a\in U$ homogeneous.
For each such term,
$\wt(a)=\wt(u)-\wt(t)<\wt(u)$.

For a homogeneous vector $w\in M(r)$, the iterate identity gives
\[
(t_{-1}a)_{-m}w
=\sum_{i\ge0}t_{-1-i}(a_{-m+i}w)
 +\sum_{i\ge0}a_{-m-1-i}(t_iw).
\]
The first sum vanishes for $i>r+\wt(a)+m-1$, and the second
for $i>r+\wt(t)-1$, by the nonnegative module grading.
Every term in the first sum is in $A_S(M)$ by derivative
closure. If $\wt(a)>0$, every term in the second sum is in
$A_S(M)$ by induction, since its coefficient index is negative.
If $a\in\C\one$, that sum is zero because
$-m-1-i\le-2$ and $\one_k=0$ for $k\ne-1$.
This proves the stronger assertion, and its case $m=1$ gives
$C_1^U(M)=A_S(M)$.
\end{proof}

\begin{lemma}\label{lem:strongindecomposables}
Let a vertex operator algebra $U$ of CFT type be strongly
generated by a homogeneous subspace
$S\subset U_+$ that is closed under derivatives.
The equality
\begin{equation}\label{eq:J}
C_2(U)+\Span\{u_{-1}v:u,v\in U_+\}
=DU+\Span\{s_{-1}v:s\in S,\ v\in U_+\}
\end{equation}
holds, and both sides lie in $U_+$.
\end{lemma}
\begin{proof}
Define
\[
J_S=DU+\Span\{s_{-1}v:s\in S,\ v\in U_+\}.
\]
Since $Du=u_{-2}\one$ for $u\in U$, one has
$DU\subseteq C_2(U)$. Also $S\subseteq U_+$, so
\[
J_S\subseteq C_2(U)+\Span\{u_{-1}v:u,v\in U_+\}.
\]
For the reverse inclusion, we prove $u_{-m}v\in J_S$ for homogeneous $u,v\in U_+$
and $m\ge1$ by induction on $\wt(u)$.
As in the proof of Lemma~\ref{lem:stronggeneration}, strong generation
and derivative closure reduce the proof to $u=t_{-1}a$, where
$t\in S$ and $a\in U$ are homogeneous and
$\wt(a)=\wt(u)-\wt(t)<\wt(u)$.
The iterate identity gives
\[
(t_{-1}a)_{-m}v
=\sum_{i\ge0}t_{-1-i}(a_{-m+i}v)
 +\sum_{i\ge0}a_{-m-1-i}(t_iv).
\]
When $i=0$ in the first sum,
\[
\wt(a_{-m}v)=\wt(a)+\wt(v)+m-1>0,
\]
so $t_{-1}(a_{-m}v)$ is a defining positive product in $J_S$.
When $i\ge1$, put $y=a_{-m+i}v$. If $y$ has positive weight,
\[
t_{-1-i}y=\frac1{i!}(D^it)_{-1}y\in J_S.
\]
If $y=\lambda\one$, then
\[
t_{-1-i}y=\frac{\lambda}{i!}D^it\in DU\qquad(i\ge1).
\]
A negative weight $y$ is zero. These cases exhaust the first sum.
In the second sum put $y=t_iv$. A positive weight $y$ is
handled by induction on the first argument $a$; when
$y=\lambda\one$,
\[
a_{-m-1-i}y=\frac{\lambda}{(m+i)!}D^{m+i}a\in DU,
\]
since $m+i\ge1$. If $a\in\C\one$, this product is zero.
Thus all positive products, and all $u_{-2}v$ with $v\in U_+$,
belong to $J_S$. The remaining $C_2$ generators are
$u_{-2}\one=Du$, while $\one_{-2}v=0$.
Therefore
\[
C_2(U)+\Span\{u_{-1}v:u,v\in U_+\}\subseteq J_S.
\]
Together with the first inclusion, this proves \eqref{eq:J}.
The grading and $D\one=0$ give $J_S\subseteq U_+$.
\end{proof}

\begin{lemma}\label{lem:indecomposable-preimage}
Let $U$ be a vertex operator algebra of CFT type, and put
\[
R_U=U/C_2(U),\qquad
J=C_2(U)+\Span\{u_{-1}v:u,v\in U_+\}.
\]
For the natural quotient map
\[
q:U_+\twoheadrightarrow(R_U)_+,\qquad u\longmapsto u+C_2(U),
\]
one has
\[
q^{-1}\bigl((R_U)_+^2\bigr)=J.
\]
The map $q$ induces a graded vector space isomorphism
\[
U_+/J\longrightarrow(R_U)_+/(R_U)_+^2,\qquad u+J\longmapsto q(u)+(R_U)_+^2.
\]
\end{lemma}
\begin{proof}
For homogeneous $a,b\in U$, a nonzero vector $a_{-2}b$ has
weight $\wt(a)+\wt(b)+1>0$. Thus $C_2(U)$ is graded and lies
in $U_+$, so $q$ is a graded surjection with
\[
\ker q=C_2(U).
\]
For $u,v\in U_+$, Proposition~\ref{prop:external-c2-algebra} gives
\[
q(u_{-1}v)=q(u)q(v)\in(R_U)_+^2.
\]
Together with $q(C_2(U))=0$, this proves
\[
J\subseteq q^{-1}\bigl((R_U)_+^2\bigr).
\]

Conversely, let $x\in U_+$ satisfy $q(x)\in(R_U)_+^2$.
By the definition of the square, there are finitely many
$\lambda_j\in\C$ and $\alpha_j,\beta_j\in(R_U)_+$ such that
\[
q(x)=\sum_{j=1}^r\lambda_j\alpha_j\beta_j.
\]
Since $q$ is surjective, choose $u_j,v_j\in U_+$ with
$q(u_j)=\alpha_j$ and $q(v_j)=\beta_j$. Then
\[
q\left(x-\sum_{j=1}^r\lambda_j(u_j)_{-1}v_j\right)=q(x)-\sum_{j=1}^r\lambda_jq(u_j)q(v_j) =0.
\]
The kernel calculation therefore gives
\[
x-\sum_{j=1}^r\lambda_j(u_j)_{-1}v_j\in C_2(U),
\]
and hence $x\in J$. This proves the reverse inclusion and the
asserted preimage equality.

The composite map
\[
U_+\xrightarrow{q}(R_U)_+
\longrightarrow(R_U)_+/(R_U)_+^2
\]
is surjective and has kernel $q^{-1}((R_U)_+^2)=J$.
It therefore induces the stated bijection on $U_+/J$.
Since $q$ and the final quotient map preserve the grading,
so does the induced isomorphism.
\end{proof}

\begin{lemma}\label{lem:Bpowers}
For $k\geq0$, $x\in V_+$, and $b\in V$,
\begin{equation}\label{eq:Bcongruence}
B_k(x,\ldots,x;b)\equiv(-1)^kk!\,R_x^kb\pmod{DV}.
\end{equation}
\end{lemma}
\begin{proof}
For $b\in V_+$ and $k\geq1$, Proposition~\ref{prop:fullsymmetry}
allows $b$ to be placed first modulo $DV$. Applying \eqref{eq:B} and
using that symmetry again on its $k$ terms gives
\[
B_k(x,\ldots,x;b) \equiv B_k(b,x,\ldots,x;x) \equiv-kB_{k-1}(x,\ldots,x;b*x)\pmod{DV}.
\]
For $b=\one$, the relation between the first and last
expressions follows directly from \eqref{eq:Bvacuum}.
By linearity this relation holds for every $b\in V$.
Starting with $B_0(;b)=b$, induction gives
\eqref{eq:Bcongruence}.
\end{proof}

\begin{proposition}\label{prop:Tpower}
For every $k\geq1$,
\[
DV+T_k(V)=DV+\cT_k(V).
\]
\end{proposition}
\begin{proof}
By Proposition~\ref{prop:diagonalspan}, the classes of
$B_{k-1}(a_1,\ldots,a_{k-1};a_k)$ with all arguments in $V_+$
are spanned by their diagonal classes.
By \eqref{eq:Bcongruence}, those diagonal classes are
$(-1)^{k-1}(k-1)!x^{[k]}+DV$. Therefore
\begin{equation}
\label{eq:positiveBpowers}
DV+\Span\{B_{k-1}(a_1,\ldots,a_{k-1};a_k):a_i\in V_+\} =DV+\cT_k(V).
\end{equation}

For a generator $B_k(a_1,\ldots,a_k;b)$ of $T_k(V)$,
split $b$ into its vacuum and positive weight parts.
The vacuum part belongs to the right hand side of
\eqref{eq:positiveBpowers} by \eqref{eq:Bvacuum}.
For the positive weight part, \eqref{eq:B} expresses the
generator as a sum of $B_{k-1}$ values with all arguments positive.
Thus $T_k(V)\subseteq DV+\cT_k(V)$.
Conversely, \eqref{eq:Bcongruence} with $b=\one$ gives
$\cT_k(V)\subseteq DV+T_k(V)$.
\end{proof}

We retain $U^{(n)}$, $W^{(n)}$, $R^{(n)}$, $(R^{(n)})_+$, and
Abe's normalization of $\phi_k^{(n)}$ and $\eta^{(n)}$ from
Section~\ref{subsec:permutationbackground} and Proposition~\ref{prop:external-abe}. Proposition~\ref{prop:external-abe-generation} gives strong
generation by $\eta^{(n)}(V_+)$ after the vacuum is removed.
For $n\geq2$, the symmetric power decomposition is
\[
\bigoplus_{k=0}^{n-1}\Sym^k(V_+)\otimes V \longrightarrow W^{(n)},\qquad (a_1\cdots a_k)\otimes b \longmapsto\phi_k^{(n-1)}(a_1,\ldots,a_k)\otimes b.
\]
Since $\eta^{(n)}(a)=\eta^{(n-1)}(a)\otimes\one
+\one^{\otimes(n-1)}\otimes a$, Abe's formula
\eqref{eq:Phimultiply} gives
\begin{equation}
\label{eq:relativeAbe}
\begin{alignedat}{1}
\eta^{(n)}(a)_{-1} \bigl(\phi_k^{(n-1)}(a_1,\ldots,a_k)\otimes b\bigr)&=\phi_{k+1}^{(n-1)}(a,a_1,\ldots,a_k)\otimes b\\
&\mathrel{\phantom{=}}+\sum_{j=1}^k \phi_k^{(n-1)}(a_1,\ldots,a*a_j,\ldots,a_k)\otimes b\\
&\mathrel{\phantom{=}}+\phi_k^{(n-1)}(a_1,\ldots,a_k)\otimes(a*b).
\end{alignedat}
\end{equation}
Here $a,a_1,\ldots,a_k\in V_+$, $b\in V$, and
$0\leq k\leq n-1$, with $\phi_n^{(n-1)}=0$.

\begin{theorem}\label{prop:seedquotient}
For every $n\geq2$, there is a graded vector space isomorphism
\[
W^{(n)}/\Crel{U^{(n)}}{W^{(n)}}\cong V/T_n(V).
\]
On spanning classes, it sends
$[\phi_k^{(n-1)}(a_1,\ldots,a_k)\otimes b]$ to
$B_k(a_1,\ldots,a_k;b)+T_n(V)$, for $0\leq k\leq n-1$.
\end{theorem}
\begin{proof}
Write $[w]=w+C_1^{U^{(n)}}(W^{(n)})$ for $w\in W^{(n)}$,
and let
\[
\pi:V\longrightarrow V/T_n(V),\qquad
q:W^{(n)}\longrightarrow W^{(n)}/C_1^{U^{(n)}}(W^{(n)}),
\quad q(w)=[w],
\]
be the quotient maps. The proof constructs $\alpha$ and $\rho$
and makes the following square commute. The induced maps
$\overline\alpha$ and $\overline\rho$ on the antidiagonal
are shown to be mutually inverse:
\[
\begin{tikzcd}[row sep={40mm,between origins},column sep={40mm,between origins}]
V \arrow[r,"\alpha"] \arrow[d,"\pi"']
  & W^{(n)}/C_1^{U^{(n)}}(W^{(n)})
    \arrow[dl,dashed,shift left=0.7ex,"\overline\rho"] \\
V/T_n(V) \arrow[ur,dashed,shift left=0.7ex,"\overline\alpha"]
  & W^{(n)} \arrow[u,"q"'] \arrow[l,"\rho"]
\end{tikzcd}
\]
\smallskip\noindent\emph{(1) Prove that $\alpha$ is surjective.}
Since $\eta^{(n)}(V_+)$ strongly generates $U^{(n)}$ and is
closed under derivatives, Lemma~\ref{lem:stronggeneration} gives
\[
\begin{alignedat}{1}
C_1^{U^{(n)}}(W^{(n)})&=\Span\bigl\{\eta^{(n)}(a)_{-1} (\phi_k^{(n-1)}(a_1,\ldots,a_k)\otimes b):\\
&\mathrel{\phantom{=}}a,a_1,\ldots,a_k\in V_+,\ b\in V,\ 0\leq k\leq n-1\bigr\}.
\end{alignedat}
\]
Define a linear map
\[
\alpha:V\longrightarrow W^{(n)}/C_1^{U^{(n)}}(W^{(n)}),\qquad b\longmapsto[\one^{\otimes(n-1)}\otimes b].
\]
We first express every spanning class in terms of $\alpha$.
For $k=0$, this is the identity
$[\phi_0^{(n-1)}\otimes b]=\alpha(B_0(;b))$.
Assume that, for some $k$ with $0\leq k\leq n-2$,
\[
[\phi_k^{(n-1)}(a_1,\ldots,a_k)\otimes b]
=\alpha\bigl(B_k(a_1,\ldots,a_k;b)\bigr)
\]
holds for all $a_1,\ldots,a_k\in V_+$ and $b\in V$.
For $a\in V_+$, solving \eqref{eq:relativeAbe} for its term
of degree $k+1$ and applying this induction hypothesis to
each term of degree $k$ gives
\[
\begin{alignedat}{1}
[\phi_{k+1}^{(n-1)}(a,a_1,\ldots,a_k)\otimes b]&=-\sum_{j=1}^k [\phi_k^{(n-1)}(a_1,\ldots,a*a_j,\ldots,a_k)\otimes b]\\
&\mathrel{\phantom{=}}-[\phi_k^{(n-1)}(a_1,\ldots,a_k)\otimes(a*b)]\\
&=\alpha\!\left( -\sum_{j=1}^k B_k(a_1,\ldots,a*a_j,\ldots,a_k;b) -B_k(a_1,\ldots,a_k;a*b)\right)\\
&=\alpha\bigl(B_{k+1}(a,a_1,\ldots,a_k;b)\bigr),
\end{alignedat}
\]
where the last equality is \eqref{eq:B}. Thus, for
$0\leq k\leq n-1$,
\begin{equation}
\label{eq:Brelative}
[\phi_k^{(n-1)}(a_1,\ldots,a_k)\otimes b] =\alpha\bigl(B_k(a_1,\ldots,a_k;b)\bigr) =[\one^{\otimes(n-1)}\otimes B_k(a_1,\ldots,a_k;b)].
\end{equation}
These classes span the quotient, so $\alpha$ is surjective.

\smallskip\noindent\emph{(2) Show that $\alpha$ descends to the quotient.}
At $k=n-1$, the term $\phi_n^{(n-1)}$ is zero.
Equations \eqref{eq:relativeAbe} and \eqref{eq:Brelative} then give
\[
\begin{aligned}
0
&=[\eta^{(n)}(a)_{-1}
  (\phi_{n-1}^{(n-1)}(a_1,\ldots,a_{n-1})\otimes b)]\\
&=\alpha\!\left(
 \sum_{j=1}^{n-1}B_{n-1}(a_1,\ldots,a*a_j,\ldots,a_{n-1};b)
 +B_{n-1}(a_1,\ldots,a_{n-1};a*b)\right)\\
&=-\alpha\bigl(B_n(a,a_1,\ldots,a_{n-1};b)\bigr).
\end{aligned}
\]
Since these $B_n$ values span $T_n(V)$, we have
$T_n(V)\subseteq\ker\alpha$. Hence there is a surjective map
\[
\overline\alpha:V/T_n(V)
 \longrightarrow W^{(n)}/C_1^{U^{(n)}}(W^{(n)}),
\qquad \overline\alpha(\pi(b))=\alpha(b).
\]

Thus $\alpha=\overline\alpha\circ\pi$.

\smallskip\noindent\emph{(3) Construct $\rho$ and check that it is well defined.}
To construct its inverse, use the symmetric power decomposition
of $W^{(n)}$ to define
\[
\rho(\phi_k^{(n-1)}(a_1,\ldots,a_k)\otimes b)
=\pi\bigl(B_k(a_1,\ldots,a_k;b)\bigr),
\qquad 0\leq k\leq n-1.
\]
This specifies a linear map $\rho:W^{(n)}\to V/T_n(V)$:
multilinearity and Lemma~\ref{lem:Bsymmetry} give
\[
B_k(a_{\sigma(1)},\ldots,a_{\sigma(k)};b)
=B_k(a_1,\ldots,a_k;b)
\qquad(\sigma\in S_k,\ 1\leq k\leq n-1),
\]
so the rule respects the symmetric power relations.

\smallskip\noindent\emph{(4) Show that $\rho$ descends to the quotient.}
For $0\leq k\leq n-2$, applying $\rho$ to
\eqref{eq:relativeAbe} yields
\[
\begin{alignedat}{1}
\rho\bigl(\eta^{(n)}(a)_{-1} (\phi_k^{(n-1)}(a_1,\ldots,a_k)\otimes b)\bigr)&=\pi\Bigl(B_{k+1}(a,a_1,\ldots,a_k;b)\\
&\mathrel{\phantom{=}}+\sum_{j=1}^k B_k(a_1,\ldots,a*a_j,\ldots,a_k;b)\\
&\mathrel{\phantom{=}}+B_k(a_1,\ldots,a_k;a*b)\Bigr)\\
&=0
\end{alignedat}
\]
by \eqref{eq:B}. At $k=n-1$, the same sum has no
$B_n$ term coming from a degree-$n$ tensor, and therefore
\[
\rho\bigl(\eta^{(n)}(a)_{-1} (\phi_{n-1}^{(n-1)}(a_1,\ldots,a_{n-1})\otimes b)\bigr)=-\pi\bigl(B_n(a,a_1,\ldots,a_{n-1};b)\bigr) =0.
\]
It follows that $\rho(C_1^{U^{(n)}}(W^{(n)}))=0$.
Define the induced linear map by
\[
\overline\rho:W^{(n)}/C_1^{U^{(n)}}(W^{(n)}) \longrightarrow V/T_n(V),\qquad \overline\rho([w])=\rho(w)\qquad(w\in W^{(n)}),
\]
where $[w]=w+C_1^{U^{(n)}}(W^{(n)})$. This definition is
independent of the representative: if $[w]=[w']$, then
\[
w-w'\in C_1^{U^{(n)}}(W^{(n)}) \quad\Longrightarrow\quad \rho(w)-\rho(w')=\rho(w-w') =0.
\]
In particular, for $0\leq k\leq n-1$,
\[
\overline\rho\bigl([\phi_k^{(n-1)}(a_1,\ldots,a_k)\otimes b]\bigr) =\pi\bigl(B_k(a_1,\ldots,a_k;b)\bigr) =B_k(a_1,\ldots,a_k;b)+T_n(V).
\]
Thus $\rho=\overline\rho\circ q$.

\smallskip\noindent\emph{(5) Verify the inverse maps and conclude.}
For every $b\in V$, the defining formulas for
$\overline\alpha$, $\alpha$, $\overline\rho$, and $\rho$ give
\[
\begin{aligned}
(\overline\rho\circ\overline\alpha)(\pi(b))
&=\overline\rho\bigl(\alpha(b)\bigr)\\
&=\overline\rho\bigl([\one^{\otimes(n-1)}\otimes b]\bigr)\\
&=\rho(\one^{\otimes(n-1)}\otimes b)\\
&=\rho(\phi_0^{(n-1)}\otimes b)\\
&=\pi\bigl(B_0(;b)\bigr)\\
&=\pi(b),
\end{aligned}
\]
where $\phi_0^{(n-1)}=\one^{\otimes(n-1)}$ and $B_0(;b)=b$.
Conversely, for $0\leq k\leq n-1$,
$a_1,\ldots,a_k\in V_+$, and $b\in V$,
\[
\begin{aligned}
(\overline\alpha\circ\overline\rho)
  \bigl([\phi_k^{(n-1)}(a_1,\ldots,a_k)\otimes b]\bigr)
&=\overline\alpha\bigl(
  \rho(\phi_k^{(n-1)}(a_1,\ldots,a_k)\otimes b)\bigr)\\
&=\overline\alpha\bigl(
  \pi(B_k(a_1,\ldots,a_k;b))\bigr)\\
&=\alpha\bigl(B_k(a_1,\ldots,a_k;b)\bigr)\\
&=[\one^{\otimes(n-1)}\otimes B_k(a_1,\ldots,a_k;b)]\\
&=[\phi_k^{(n-1)}(a_1,\ldots,a_k)\otimes b].
\end{aligned}
\]
The first four equalities use the definitions of
$\overline\rho$, $\rho$, $\overline\alpha$, and $\alpha$,
respectively; the last is \eqref{eq:Brelative}.
Since $\pi$ is surjective and the classes in the second
calculation span $W^{(n)}/C_1^{U^{(n)}}(W^{(n)})$, linearity gives
\[
\overline\rho\circ\overline\alpha
=\operatorname{id}_{V/T_n(V)},\qquad
\overline\alpha\circ\overline\rho
=\operatorname{id}_{W^{(n)}/C_1^{U^{(n)}}(W^{(n)})}.
\]
Thus $\overline\rho$ is the isomorphism in the statement.
Finally, the recurrence gives
\[
B_k(V_{r_1},\ldots,V_{r_k};V_s)
\subseteq V_{r_1+\cdots+r_k+s},
\]
which is also the weight of the corresponding symmetric tensor.
Both maps therefore preserve the grading.
\end{proof}

For the second quotient, recall
$R^{(n)}=U^{(n)}/C_2(U^{(n)})$.

\begin{theorem}\label{prop:positivequotient}
For every $n\geq1$, $\eta^{(n)}$ induces an isomorphism
\begin{equation}\label{eq:positivequotient}
V_+/(DV+\cT_{n+1}(V))\ \cong\ (R^{(n)})_+/(R^{(n)})_+^2.
\end{equation}
\end{theorem}
\begin{proof}
By Lemma~\ref{lem:indecomposable-preimage}, the preimage of $(R^{(n)})_+^2$ under $U^{(n)}_+\to(R^{(n)})_+$ is
\[
\begin{aligned}
J_n
&=C_2(U^{(n)})+
  \Span\{u_{-1}v:u,v\in U^{(n)}_+\}\\
&=D_{(0)}U^{(n)}+
  \Span\{\eta^{(n)}(a)_{-1}v:a\in V_+,\ v\in U^{(n)}_+\},
\end{aligned}
\]
where the second equality is \eqref{eq:J} in
Lemma~\ref{lem:strongindecomposables}. Hence the canonical map
\[
U^{(n)}_+/J_n\longrightarrow(R^{(n)})_+/(R^{(n)})_+^2,
\qquad u+J_n\longmapsto (u+C_2(U^{(n)}))+(R^{(n)})_+^2,
\]
is an isomorphism. Write $[u]=u+J_n$, and let
\[
\pi:V_+\longrightarrow V_+/(DV+\cT_{n+1}(V)),\qquad
q:U^{(n)}_+\longrightarrow U^{(n)}_+/J_n,\quad q(u)=[u],
\]
be the quotient maps. Define
$\alpha:V_+\to U^{(n)}_+/J_n$ by $\alpha(a)=[\eta^{(n)}(a)]$.
As in Theorem~\ref{prop:seedquotient}, the proof constructs
$\rho$ and makes the following square commute. The induced maps
$\overline\alpha$ and $\overline\rho$ on the antidiagonal
are shown to be mutually inverse:
\[
\begin{tikzcd}[row sep={40mm,between origins},column sep={40mm,between origins}]
V_+ \arrow[r,"\alpha"] \arrow[d,"\pi"']
  & U^{(n)}_+/J_n
    \arrow[dl,dashed,shift left=0.7ex,"\overline\rho"] \\
V_+/(DV+\cT_{n+1}(V))
    \arrow[ur,dashed,shift left=0.7ex,"\overline\alpha"]
  & U^{(n)}_+ \arrow[u,"q"'] \arrow[l,"\rho"]
\end{tikzcd}
\]
\smallskip\noindent\emph{(1) Prove that $\alpha$ is surjective.}
For $k=1$, the equality $\phi_1^{(n)}=\eta^{(n)}$ and
$B_0(;a_1)=a_1$ give
$[\phi_1^{(n)}(a_1)]=\alpha(B_0(;a_1))$.
Assume that, for some $k$ with $1\leq k<n$,
\[
[\phi_k^{(n)}(a_1,\ldots,a_k)]
=\alpha\bigl(B_{k-1}(a_1,\ldots,a_{k-1};a_k)\bigr)
\]
holds for all $a_1,\ldots,a_k\in V_+$.
For $a\in V_+$, since
$\eta^{(n)}(a)_{-1}\phi_k^{(n)}(a_1,\ldots,a_k)\in J_n$,
Abe's formula \eqref{eq:Phimultiply} and this induction
hypothesis give
\[
\begin{alignedat}{1}
[\phi_{k+1}^{(n)}(a,a_1,\ldots,a_k)]&=-\sum_{j=1}^k [\phi_k^{(n)}(a_1,\ldots,a*a_j,\ldots,a_k)]\\
&=\alpha\Bigl( -\sum_{j=1}^{k-1} B_{k-1}(a_1,\ldots,a*a_j,\ldots,a_{k-1};a_k)\\
&\mathrel{\phantom{=}}-B_{k-1}(a_1,\ldots,a_{k-1};a*a_k)\Bigr)\\
&=\alpha\bigl(B_k(a,a_1,\ldots,a_{k-1};a_k)\bigr)
\end{alignedat}
\]
by \eqref{eq:B}. Thus
\begin{equation}\label{eq:Bphi}
\begin{aligned}
[\phi_k^{(n)}(a_1,\ldots,a_k)]
&=\alpha\bigl(B_{k-1}(a_1,\ldots,a_{k-1};a_k)\bigr)\\
&=[\eta^{(n)}(B_{k-1}(a_1,\ldots,a_{k-1};a_k))],
\qquad 1\leq k\leq n.
\end{aligned}
\end{equation}
Because $U^{(n)}_+$ is spanned by these positive degree
symmetric tensors, $\alpha$ is surjective.

\smallskip\noindent\emph{(2) Show that $\alpha$ descends to the quotient.}
Translation covariance gives, for every $a\in V$,
\[
\alpha(Da)=[\eta^{(n)}(Da)] =[D_{(0)}\eta^{(n)}(a)] =0,
\]
so $DV\subseteq\ker\alpha$. At symmetric degree $n$,
the boundary convention $\phi_{n+1}^{(n)}=0$ gives
\[
\begin{aligned}
0 &=[\eta^{(n)}(a)_{-1}\phi_n^{(n)}(a_1,\ldots,a_n)]\\
&=\sum_{j=1}^n [\phi_n^{(n)}(a_1,\ldots,a*a_j,\ldots,a_n)]\\
&=\alpha\Bigl( \sum_{j=1}^{n-1} B_{n-1}(a_1,\ldots,a*a_j,\ldots,a_{n-1};a_n) +B_{n-1}(a_1,\ldots,a_{n-1};a*a_n)\Bigr)\\
&=-\alpha\bigl(B_n(a,a_1,\ldots,a_{n-1};a_n)\bigr).
\end{aligned}
\]
All arguments in this expression belong to $V_+$.
Equation~\eqref{eq:positiveBpowers}, with $k=n+1$, now yields
\[
DV+\cT_{n+1}(V) =DV+\Span\{B_n(a_0,\ldots,a_{n-1};a_n):a_i\in V_+\} \subseteq\ker\alpha.
\]
Thus $\alpha$ induces a surjection
\[
\overline\alpha:V_+/(DV+\cT_{n+1}(V))
 \longrightarrow U^{(n)}_+/J_n,
\qquad \overline\alpha(\pi(a))=[\eta^{(n)}(a)].
\]

Thus $\alpha=\overline\alpha\circ\pi$.

\smallskip\noindent\emph{(3) Construct $\rho$ and check that it is well defined.}
Define the inverse on the symmetric power summands by
\[
\rho(\phi_k^{(n)}(a_1,\ldots,a_k))
=\pi\bigl(B_{k-1}(a_1,\ldots,a_{k-1};a_k)\bigr),
\qquad 1\leq k\leq n.
\]
For every $\sigma\in S_k$, Proposition~\ref{prop:fullsymmetry}
gives
\[
B_{k-1}(a_{\sigma(1)},\ldots,a_{\sigma(k-1)};a_{\sigma(k)}) -B_{k-1}(a_1,\ldots,a_{k-1};a_k) \in DV\subseteq\ker\pi.
\]
Together with multilinearity, this shows that $\rho$ is a
well defined linear map on $U^{(n)}_+$.

\smallskip\noindent\emph{(4) Show that $\rho$ descends to the quotient.}
For $1\leq k<n$, \eqref{eq:Phimultiply} and \eqref{eq:B} give
\[
\begin{alignedat}{1}
\rho\bigl(\eta^{(n)}(a)_{-1}\phi_k^{(n)}(a_1,\ldots,a_k)\bigr)&=\pi\Bigl(B_k(a,a_1,\ldots,a_{k-1};a_k)\\
&\mathrel{\phantom{=}}+\sum_{j=1}^{k-1} B_{k-1}(a_1,\ldots,a*a_j,\ldots,a_{k-1};a_k)\\
&\mathrel{\phantom{=}}+B_{k-1}(a_1,\ldots,a_{k-1};a*a_k)\Bigr)\\
&=0.
\end{alignedat}
\]
For $k=n$, the degree-$n+1$ term is absent, so
\[
\rho\bigl(\eta^{(n)}(a)_{-1}\phi_n^{(n)}(a_1,\ldots,a_n)\bigr)=-\pi\bigl(B_n(a,a_1,\ldots,a_{n-1};a_n)\bigr) =0,
\]
where the last equality uses \eqref{eq:positiveBpowers}.
Since the $\phi_k^{(n)}$ with $1\leq k\leq n$ span $U^{(n)}_+$,
these calculations prove
\[
\rho\bigl(\Span\{\eta^{(n)}(a)_{-1}v:
 a\in V_+,\ v\in U^{(n)}_+\}\bigr)=0.
\]

It remains to check the derivative part of $J_n$.
Translation acts on each tensor factor, so
\[
D_{(0)}\phi_k^{(n)}(a_1,\ldots,a_k)
=\sum_{j=1}^k\phi_k^{(n)}(a_1,\ldots,Da_j,\ldots,a_k).
\]
Applying $\rho$ and then Lemma~\ref{lem:Bderivation} gives
\[
\begin{alignedat}{1}
\rho\bigl(D_{(0)}\phi_k^{(n)}(a_1,\ldots,a_k)\bigr)&=\pi\Bigl(\sum_{j=1}^{k-1} B_{k-1}(a_1,\ldots,Da_j,\ldots,a_{k-1};a_k)\\
&\mathrel{\phantom{=}}+B_{k-1}(a_1,\ldots,a_{k-1};Da_k)\Bigr)\\
&=\pi\bigl(DB_{k-1}(a_1,\ldots,a_{k-1};a_k)\bigr)\\
&=0.
\end{alignedat}
\]
Also $D\one=0$. Therefore $\rho(D_{(0)}U^{(n)})=0$, and hence
$\rho(J_n)=0$. Empty sums in these calculations are zero;
in particular, the argument includes $n=1$ and $k=1$.
We obtain an induced map
\[
\overline\rho:U^{(n)}_+/J_n
 \longrightarrow V_+/(DV+\cT_{n+1}(V)),
\qquad \overline\rho([u])=\rho(u).
\]
Thus $\rho=\overline\rho\circ q$.

\smallskip\noindent\emph{(5) Verify the inverse maps and conclude.}
For every $a\in V_+$, the defining formulas for
$\overline\alpha$, $\alpha$, $\overline\rho$, and $\rho$ give
\[
(\overline\rho\circ\overline\alpha)(\pi(a)) =\overline\rho\bigl(\alpha(a)\bigr) =\overline\rho\bigl([\eta^{(n)}(a)]\bigr) =\rho(\eta^{(n)}(a)) =\rho(\phi_1^{(n)}(a)) =\pi\bigl(B_0(;a)\bigr) =\pi(a),
\]
where $\phi_1^{(n)}(a)=\eta^{(n)}(a)$ and $B_0(;a)=a$.
Conversely, for $1\leq k\leq n$ and
$a_1,\ldots,a_k\in V_+$,
\[
\begin{aligned}
(\overline\alpha\circ\overline\rho)
  \bigl([\phi_k^{(n)}(a_1,\ldots,a_k)]\bigr)
&=\overline\alpha\bigl(
  \rho(\phi_k^{(n)}(a_1,\ldots,a_k))\bigr)\\
&=\overline\alpha\bigl(
  \pi(B_{k-1}(a_1,\ldots,a_{k-1};a_k))\bigr)\\
&=\alpha\bigl(B_{k-1}(a_1,\ldots,a_{k-1};a_k)\bigr)\\
&=[\eta^{(n)}(B_{k-1}(a_1,\ldots,a_{k-1};a_k))]\\
&=[\phi_k^{(n)}(a_1,\ldots,a_k)].
\end{aligned}
\]
The first four equalities use the definitions of
$\overline\rho$, $\rho$, $\overline\alpha$, and $\alpha$,
respectively; the last is \eqref{eq:Bphi}.
Since $\pi$ is surjective and the classes in the second
calculation span $U^{(n)}_+/J_n$, linearity gives
\[
\overline\rho\circ\overline\alpha
=\operatorname{id}_{V_+/(DV+\cT_{n+1}(V))},\qquad
\overline\alpha\circ\overline\rho
=\operatorname{id}_{U^{(n)}_+/J_n}.
\]
Composing $\overline\alpha$ with the canonical identification
$U^{(n)}_+/J_n\cong(R^{(n)})_+/(R^{(n)})_+^2$ proves
\eqref{eq:positivequotient} with the stated map induced by
$\eta^{(n)}$.
\end{proof}

\subsection{Heisenberg modules and weight one translation bounds}
\label{subsec:gradedtools}

For each $n\ge2$, we seek an exponent $r_n$ such that
$D^{r_n}V_1\subseteq T_n(V)$ for every simple,
self contragredient vertex operator algebra $V$ of CFT type.
The Heisenberg decomposition and $C_1$-cofiniteness first give
a translation relation in $\mathcal H$. The invariant bilinear
form then provides a basis of $V_1$ to which this relation
transfers by graded vertex algebra homomorphisms.
We first establish the module lemmas for a general base algebra $U$,
then apply them to $(\mathcal H^{\otimes n})^{O(n,\C)}$, and finally
transfer the resulting relation to the target algebra $V$.

\begin{lemma}\label{lem:imageclosure}
Finite sums and quotients of $C_1$-cofinite $U$-modules are
$C_1$-cofinite. If a $C_1$-cofinite graded $U$-submodule $Z$ of an
extension $V$ satisfies
\[
P_1=Z,\qquad P_{j+1}=Z\diamond P_j,\qquad P_r=V
\]
for some finite $r$, then $V$ is $C_1$-cofinite as a $U$-module.
Here all $P_j$ are assumed to be $\N$-gradable weak modules.
\end{lemma}
\begin{proof}
For a surjective $U$-module map $p:M\twoheadrightarrow N$,
\[
p(C_1^U(M))
=\Span\{u_{-1}p(m):u\in U_+,\ m\in M\}=C_1^U(N).
\]
It induces a surjection $M/C_1^U(M)\twoheadrightarrow N/C_1^U(N)$.
For a finite direct sum one has
\[
C_1^U\!\left(\bigoplus_{i=1}^s M_i\right)
=\bigoplus_{i=1}^s C_1^U(M_i),
\]
\[
\dim\bigl((\bigoplus_i M_i)/C_1^U(\bigoplus_i M_i)\bigr)
=\sum_i\dim\bigl(M_i/C_1^U(M_i)\bigr).
\]
For $U$-submodules $M_1,\ldots,M_s\subseteq M$, their sum is a quotient
of their direct sum under $(m_1,\ldots,m_s)\mapsto\sum_i m_i$.
These identities prove both closure assertions.

For graded $U$-submodules $M,N\subset V$, the commutator identity
gives, for $u\in U$, $m\in M$, $n\in N$, and $k,\ell\in\Z$,
\[
u_k(m_\ell n)=m_\ell(u_kn)
+\sum_{j\ge0}\binom kj(u_jm)_{k+\ell-j}n\in M\diamond N.
\]
Thus $M\diamond N$ is a $U$-submodule. The restriction of $Y_V$
is an intertwining operator with this target, and it is
surjective because its coefficients span $M\diamond N$.
Assume inductively that $P_j$ is $C_1$-cofinite.
Theorem~\ref{thm:image} at $m=1$, applied to $Z,P_j$, gives
\[
\dim P_{j+1}/C_1^U(P_{j+1})<f_1(Z,P_j)<\infty.
\]
Starting with $P_1=Z$ and repeating this for $j=1,\ldots,r-1$
proves $C_1$-cofiniteness of $P_r=V$.
\end{proof}

For the next two lemmas, let $U$ be a conformal graded vertex
algebra, and write $D=(\omega_U)_0$ on its modules.

\begin{lemma}\label{lem:cyclic-transfer}
Let $M$ be a graded $U$-module with
$\dim M/C_1^U(M)<\infty$. For every $d\in\N$ and $w\in M_d$,
there is an integer $r\ge0$ such that
\begin{equation*}\label{eq:cyclic-transfer}
D^kw\in C_1^U(M)\qquad(k\ge r).
\end{equation*}
\end{lemma}
\begin{proof}
The subspace $C_1^U(M)$ is graded, because
$\wt(a_{-1}m)=\wt(a)+\wt(m)$ for homogeneous $a,m$.
Finite dimensionality of the quotient gives $t\ge0$ such that, for every
$i>t$,
\[
M_i\subseteq C_1^U(M).
\]
Since $\wt(D^kw)=d+k$, choose
$r=\max\{0,t-d+1\}$. Then, for $k\ge r$,
\[
D^kw\in M_{d+k}\subseteq C_1^U(M).\qedhere
\]
\end{proof}

Let $N$ be a graded $U$-module with
\[
N=\bigoplus_{\alpha\in\Lambda}N_\alpha,\qquad
N_\alpha\text{ an irreducible }U\text{-module},\qquad
\dim\bigl(N_\alpha/C_1^U(N_\alpha)\bigr)<\infty
\quad(\alpha\in\Lambda).
\]
For $v\in N$, write
\[
M_v:=U\cdot v =\Span_{\C}\bigl\{u^1_{j_1}\cdots u^s_{j_s}v: s\in\N,\ u^i\in U,\ j_i\in\Z\ (1\le i\le s)\bigr\},
\]
where $s=0$ gives $v$.

\begin{lemma}\label{lem:cyclic-submodule}
For every $d\in\N$ and $v\in N_d$, one has
$\dim\bigl(M_v/C_1^U(M_v)\bigr)<\infty$, and there is an integer $r_v\ge0$ such that,
for every integer $k\ge r_v$,
\[
D^kv\in C_1^U(M_v).
\]
\end{lemma}
\begin{proof}
Fix $d\in\N$ and $v\in N_d$.
Since $N$ is an algebraic direct sum, there is a finite
$F\subseteq\Lambda$ such that
\[
v\in N_F:=\bigoplus_{\alpha\in F}N_\alpha,
\qquad M_v=U\cdot v\subseteq N_F.
\]
A finite direct sum of simple modules is semisimple, so
$N_F=M_v\oplus N'_F$ for a $U$-submodule $N'_F$.
The corresponding $U$-linear projection
$p:N_F\twoheadrightarrow M_v$ satisfies
\[
\begin{alignedat}{1}
p(C_1^U(N_F))&=C_1^U(M_v),\\
\dim\bigl(M_v/C_1^U(M_v)\bigr)& \le\dim\bigl(N_F/C_1^U(N_F)\bigr) =\sum_{\alpha\in F}\dim\bigl(N_\alpha/C_1^U(N_\alpha)\bigr)<\infty,
\end{alignedat}
\]
by Lemma~\ref{lem:imageclosure}. Homogeneity of $v$ makes $M_v$
graded, and $D=(\omega_U)_0$ preserves $M_v$.
The remaining assertion follows from Lemma~\ref{lem:cyclic-transfer}.
\end{proof}

We apply the Heisenberg decomposition and cofiniteness results of
Section~\ref{subsec:heisenbergmodules}, together with
Lemma~\ref{lem:cyclic-submodule}, to the module for the stabilizer of a tensor position.

\begin{proposition}\label{cor:heis-relative}
For every $n\ge2$ there is an integer $r_n\ge0$, depending only
on $n$, such that
\[
D_{(0)}^{r_n}h^{(n)}\in
C_1^{U^{(n)}(\mathcal H)}(W^{(n)}(\mathcal H)).
\]
\end{proposition}
\begin{proof}
Put $P=(\mathcal H^{\otimes n})^{O(n,\C)}$,
$w=h^{(n)}=\one^{\otimes(n-1)}\otimes h$, and $M=P\cdot w$.
Theorem~\ref{thm:heisenbergmodules} gives
\[
\mathcal H^{\otimes n}=\bigoplus_\lambda L_\lambda\otimes M_\lambda.
\]
Here $P$ acts only on the $M_\lambda$ factor. Choosing a basis
of the finite dimensional space $L_\lambda$ identifies
$L_\lambda\otimes M_\lambda$ with $d_\lambda=\dim L_\lambda$
copies of the irreducible $P$-module $M_\lambda$.
Each copy is an irreducible $P$-submodule of
$\mathcal H^{\otimes n}$, so Proposition~\ref{prop:heisenberg-wick} gives
\[
M_\lambda/C_1^P(M_\lambda)
=M_\lambda/(M_\lambda)_{\mathrm{Wick}},\qquad
\dim M_\lambda/C_1^P(M_\lambda)<\infty.
\]
The equality holds term by term: a positive weight normally
ordered product is exactly $a_{-1}m$ in \eqref{eq:c1}.

Because the decomposition is an algebraic direct sum, there is
a finite set $F$ of indices such that
\[
w\in N_F:=\bigoplus_{\lambda\in F}L_\lambda\otimes M_\lambda.
\]
Every summand is $P$-stable, so
\[
M=P\cdot w\subseteq N_F
\cong\bigoplus_{\lambda\in F}M_\lambda^{\oplus d_\lambda}
\qquad\text{as }P\text{-modules}.
\]
The finite direct sum $N_F$ is semisimple, so its $P$-submodule
$M$ is a direct summand. Since $w$ has weight one, applying
Lemma~\ref{lem:cyclic-submodule} with $U=P$, $N=N_F$,
$d=1$, and $v=w$, and using the estimate in its proof, gives
\[
\dim\bigl(M/C_1^P(M)\bigr)
\le\sum_{\lambda\in F}d_\lambda
\dim\bigl(M_\lambda/C_1^P(M_\lambda)\bigr)<\infty.
\]
Permutation matrices are orthogonal, so $P\subset U^{(n)}(\mathcal H)$.
If $g\in S_{n-1}$ fixes the last position, then $gw=w$ and
$ga=a$ for all $a\in P$. On every product of vertex operator coefficients,
\[
g(a^1_{k_1}\cdots a^t_{k_t}w)
=(ga^1)_{k_1}\cdots(ga^t)_{k_t}gw
=a^1_{k_1}\cdots a^t_{k_t}w.
\]
Thus $M\subseteq W^{(n)}(\mathcal H)$.
The common conformal vector and translation are
\[
\Omega=\frac12\sum_{i=1}^n(h^{(i)})_{-1}h^{(i)}\in P,
\qquad D_{(0)}=\Omega_0=\sum_{i=1}^nD^{(i)}.
\]
In particular $D_{(0)}$ preserves $M$ and
\[
D_{(0)}^r w
=\one^{\otimes(n-1)}\otimes D^rh.
\]
Lemma~\ref{lem:cyclic-transfer}, applied to the graded
$P$-module $M$, gives an integer $r_0\ge0$ such that
\[
D_{(0)}^r w\in C_1^P(M)\qquad(r\ge r_0).
\]
Set $r_n=\max\{n-1,r_0\}$. Then
\[
D_{(0)}^{r_n}w\in C_1^P(M)
\subseteq C_1^{U^{(n)}(\mathcal H)}(W^{(n)}(\mathcal H)).
\]
Here the inclusion follows from $P_+\subseteq U^{(n)}(\mathcal H)_+$
and $M\subseteq W^{(n)}(\mathcal H)$ by the definition of $C_1$.
Every space and every bound in this construction is determined
by $n$, which proves the stated uniformity.
\end{proof}

\begin{lemma}\label{lem:weightoneform}
Let $V$ be a simple, self contragredient vertex operator algebra of CFT type.
Then $L(1)V_1=0$, and the form $\beta$ defined by
$a_1b=\beta(a,b)\one$ on $V_1$ is symmetric and nondegenerate.
There is a basis $h_i$ of $V_1$ with $(h_i)_1h_i=\one$,
$(h_i)_0h_i=0$, and $(h_i)_jh_i=0$ for $j\geq2$.
\end{lemma}
\begin{proof}
Self contragredience and Proposition~\ref{prop:external-invariant-form}
give a nondegenerate invariant symmetric form $(\ ,\ )$ on $V$.
Applying invariance to the conformal vector gives
\[
(L(n)u,v)=(u,L(-n)v).
\]
In particular, $L(0)$ is self adjoint, so distinct weight spaces
are orthogonal. The restriction of the form to each $V_d$ is
nondegenerate. Since $V_0=\C\one$, normalize it by
$(\one,\one)=1$.
For $a\in V_1$, write $L(1)a=\lambda\one$. Then
\[
\lambda=(L(1)a,\one)=(a,L(-1)\one) =0,
\]
and hence $L(1)V_1=0$.

For $a,b\in V_1$, invariance and the vacuum identity now give
\[
(Y(a,z)b,\one) =(b,Y(-z^{-2}a,z^{-1})\one) =-z^{-2}(b,e^{z^{-1}D}a)=-z^{-2}(b,a).
\]
The last equality uses orthogonality of different weights.
Taking the coefficient of $z^{-2}$ yields
\[
\beta(a,b)=(a_1b,\one)=-(b,a).
\]
Thus $\beta$ is symmetric and nondegenerate. A nondegenerate
symmetric form over $\C$ has an orthonormal basis: choose a
vector of nonzero square, scale it to square one, and repeat on
its nondegenerate orthogonal complement. This gives $h_i$ with
$\beta(h_i,h_i)=1$.
For $j\ge2$,
\[
\wt((h_i)_jh_i)=1-j<0,
\qquad (h_i)_jh_i=0.
\]
Finally, the coefficient of $z^{-1}$ in skew symmetry gives
\[
(h_i)_0h_i=-(h_i)_0h_i+D((h_i)_1h_i)
=-(h_i)_0h_i+D\one=-(h_i)_0h_i.
\]
Therefore $(h_i)_0h_i=0$, proving all assertions.
\end{proof}

\begin{lemma}\label{lem:weightonebase}
Let $V$ be simple, self contragredient and of CFT type.
For $r_n$ from Proposition~\ref{cor:heis-relative},
\begin{equation}\label{eq:weightonebase}
D^{r_n}V_1\subseteq T_n(V).
\end{equation}
\end{lemma}
\begin{proof}
Apply Theorem~\ref{prop:seedquotient} to $\mathcal H$.
Proposition~\ref{cor:heis-relative} gives
\[
0=[D_{(0)}^{r_n}h^{(n)}]\longmapsto D^{r_n}h+T_n(\mathcal H),
\quad\text{hence}\quad D^{r_n}h\in T_n(\mathcal H).
\]
Choose the basis $h_1,\ldots,h_t$ from
Lemma~\ref{lem:weightoneform}. For each $i$, the identities
$(h_i)_jh_i=\delta_{j,1}\one$ for $j\ge0$ allow us to apply
Proposition~\ref{prop:heisenberg-homomorphism}, giving a graded
vertex algebra homomorphism $f_i:\mathcal H\to V$ with
$f_i(h)=h_i$ and $f_iD=Df_i$. By \eqref{eq:functoriality},
\[
D^{r_n}h_i=f_i(D^{r_n}h)\in f_i(T_n(\mathcal H))
\subseteq T_n(V).
\]
For $v=\sum_i\alpha_i h_i$ this implies
$D^{r_n}v=\sum_i\alpha_iD^{r_n}h_i\in T_n(V)$.
If $V_1=0$ the assertion is empty.
\end{proof}

\subsection{Proof of Theorem~\ref*{thm:criterion}}
\label{sec:translationcriterion}\label{sec:closure}\label{sec:criterionproof}

This subsection establishes the translation bounds and relative
finiteness needed for Theorem~\ref{thm:criterion}.
Fix $n\ge2$ and retain $U^{(n)}$, $W^{(n)}$, $(R^{(n-1)})_+$ and
$T_n(V)$ from Sections~\ref{subsec:permutationbackground},
\ref{sec:symmetricoperations} and~\ref{sec:symmetricquotients}.
The two quotient identifications proved there give
\begin{align}
 W^{(n)}/C_1^{U^{(n)}}(W^{(n)})
 &\cong V/T_n(V),\label{eq:seed-summary}\\
 V/(\C\one+T_n(V)+DV)
 &\cong(R^{(n-1)})_+/(R^{(n-1)})_+^2.\label{eq:cotangent-summary}\notag
\end{align}
For the second identity, first use
$DV+T_n(V)=DV+\cT_n(V)$ and then apply
Theorem~\ref{prop:positivequotient} in degree $n-1$.
The first map sends
$[\one^{\otimes(n-1)}\otimes a]$ to $a+T_n(V)$.
All these maps preserve weights.

For homogeneous $v\in V_d$, the standard skew symmetry
identity \eqref{eq:skew} also gives
\begin{equation}\label{eq:v0omega}
v_0\omega=(d-1)Dv
+\sum_{j\geq2}\frac{(-1)^{j+1}}{j!}D^jL(j-1)v.
\end{equation}
For clarity, extract $z^{-1}$ from
$Y(v,z)\omega=e^{zD}Y(\omega,-z)v$. The general term is
$(-1)^{j+1}D^jL(j-1)v/j!$. The $j=0$ and $j=1$ terms are
$-Dv$ and $D L(0)v=dDv$, respectively; their sum is
$(d-1)Dv$. For $j>d+1$, $L(j-1)v$ has negative weight and
vanishes. This proves \eqref{eq:v0omega}, including its
truncation, directly from skew symmetry.

\begin{proposition}\label{prop:translation}
Let $V$ be a simple, self contragredient vertex operator algebra of CFT type.
If $D^r\omega\in T_n(V)$ for some $r\ge0$, then
\begin{equation*}\label{eq:translationbound}
D^NV\subseteq T_n(V)
\end{equation*}
for
\[
N=\begin{cases}
r+1,&V_1=0,\\
\max\{r_n,r+1\},&V_1\ne0.
\end{cases}
\]
\end{proposition}
\begin{proof}
By Lemma~\ref{lem:Tstability}, $DT_n(V)\subseteq T_n(V)$.
Since $N-1\ge r$, our hypothesis implies
\begin{equation}\label{eq:conformaldN}
D^{N-1}\omega=D^{N-1-r}(D^r\omega)\in T_n(V).
\end{equation}
The vacuum satisfies $D^N\one=0$. If $V_1\ne0$, then
$N\ge r_n$ and \eqref{eq:weightonebase} give
\[
D^NV_1=D^{N-r_n}(D^{r_n}V_1)\subseteq T_n(V).
\]
If $V_1=0$, this step is unnecessary.

Assume now that $d\ge2$, that $v\in V_d$, and that
$D^NV_e\subseteq T_n(V)$ for every $e<d$.
The identity \eqref{eq:v0omega} gives
\[
v_0\omega=(d-1)Dv+
\sum_{j\ge2}\frac{(-1)^{j+1}}{j!}D^jL(j-1)v.
\]
Translation covariance gives $[D,v_0]=0$. Applying
$D^{N-1}$ and rearranging therefore yields
\begin{equation*}
(d-1)D^Nv=v_0D^{N-1}\omega-
\sum_{j\ge2}\frac{(-1)^{j+1}}{j!}
D^{N+j-1}L(j-1)v.
\end{equation*}
The first term belongs to $T_n(V)$ by
\eqref{eq:conformaldN} and stability under the operators $v_0$.
For a nonzero summand in the sum,
\[
L(j-1)v\in V_{d-j+1},\qquad 0\le d-j+1<d,
\]
and
\[
D^{N+j-1}L(j-1)v
=D^{j-1}\bigl(D^NL(j-1)v\bigr)\in T_n(V)
\]
by the induction hypothesis and $D$-stability. Terms of
negative weight are zero. The sum is thus finite and entirely
contained in $T_n(V)$. Since $d-1\ne0$, division proves
$D^Nv\in T_n(V)$. Induction proves the assertion in all
weights, and every vector of $V$ is a finite sum of homogeneous
vectors. The integer $N$ is independent of $d$ and $v$.
\end{proof}

A graded $\C[D]$-module is a vector space
$Q=\bigoplus_{d\ge0}Q_d$ with $DQ_d\subseteq Q_{d+1}$.

\begin{lemma}\label{lem:gradedD}
Let $Q$ be a nonnegatively graded $\C[D]$-module.
If $\dim Q/DQ<\infty$ and $D^NQ=0$ for some integer $N\ge1$,
then $Q$ is finite dimensional.
\end{lemma}
\begin{proof}
Choose homogeneous $q_1,\ldots,q_s$ whose classes form a basis
of $Q/DQ$. Write $e_i=\deg q_i$ and set
$Q'=\sum_i\C[D]q_i$.
For $q\in Q_d$, the homogeneous basis of $Q/DQ$ gives
\[
q=\sum_{e_i=d}\lambda_iq_i+Dq',\qquad q'\in Q_{d-1},
\]
where $Q_{-1}=0$. Induction on $d$ proves $Q_d\subseteq Q'$,
and hence
\begin{equation*}\label{eq:gradedDgeneration}
Q=\sum_{i=1}^s\C[D]q_i.
\end{equation*}
Since $D^NQ=0$, each $\C[D]q_i$ is a quotient of the
finite dimensional space $\C[D]/(D^N)$. Their finite sum is
therefore finite dimensional. If $s=0$, the degree induction
gives $Q=0$.
\end{proof}

Set
\begin{equation*}
Q_n=V/(\C\one+T_n(V)).
\end{equation*}
Because $D\one=0$ and $DT_n(V)\subseteq T_n(V)$, $D$
descends to a degree one endomorphism of $Q_n$.
Also,
\begin{equation}\label{eq:QnD}
Q_n/DQ_n=V/(\C\one+T_n(V)+DV)
\cong(R^{(n-1)})_+/(R^{(n-1)})_+^2.
\end{equation}

We use this quotient to prove the induction step in tensor
degree.

\begin{theorem}\label{thm:seedfinite}
Let $V$ be a vertex operator algebra of CFT type, and fix $n\ge2$.
Suppose $U^{(n-1)}$ is $C_2$-cofinite and
$D^NV\subseteq T_n(V)$ for some $N\ge1$. Then $W^{(n)}$ is
$C_1$-cofinite as a $U^{(n)}$-module.
If $V$ is also $C_2$-cofinite, then $U^{(n)}$ is
$C_2$-cofinite.
\end{theorem}
\begin{proof}
The algebra $R^{(n-1)}$ is finite dimensional, so
\eqref{eq:QnD} gives $\dim Q_n/DQ_n<\infty$.
The hypothesis $D^NV\subseteq T_n(V)$ implies $D^NQ_n=0$.
Lemma~\ref{lem:gradedD} therefore makes $Q_n$ finite dimensional.
Since $T_n(V)\subseteq V_+$, the vacuum class together with
finitely many lifts from $Q_n$ spans $V/T_n(V)$.
Equation~\eqref{eq:seed-summary} now proves the first assertion.
It remains to propagate this relative finiteness to $X=T^n(V)$.

For $1\le i\le n$, let $H_i$ be the stabilizer of the $i$th
position and put $W_i=X^{H_i}$. Thus $W_n=W^{(n)}$.
Choose $\sigma_i\in S_n$ with $\sigma_i(n)=i$.
For $u\in U^{(n)}$, $w\in W_n$ and $k\in\Z$,
\[
\sigma_i(u_kw)=(\sigma_i u)_k\sigma_i w=u_k\sigma_i w.
\]
Consequently $\sigma_i:W_n\to W_i$ is an isomorphism of
$U^{(n)}$-modules and maps their relative $C_1$ spaces onto
one another. In particular every $W_i$ is $C_1$-cofinite.
By Lemma~\ref{lem:imageclosure}, the sum
\[
Z=\sum_{i=1}^nW_i
\]
is $C_1$-cofinite as a $U^{(n)}$-module.

Define $P_1=Z$ and $P_{j+1}=Z\diamond P_j$.
The coefficient image argument in the proof of
Lemma~\ref{lem:imageclosure} shows that each $P_j$ is a
$U^{(n)}$-submodule of $X$. It inherits the nonnegative grading
of $X$, so all the $P_j$ are $\N$-gradable weak modules.

Every vector supported in a single tensor position $a^{(i)}$ belongs to $W_i\subset Z$.
For $a_1,\ldots,a_n\in V$, the tensor product vertex operator
gives, with right associated evaluation,
\begin{equation}\label{eq:puretensorP}
a_1\otimes\cdots\otimes a_n
=(a_1^{(1)})_{-1}\cdots(a_{n-1}^{(n-1)})_{-1}a_n^{(n)}.
\end{equation}
Indeed, at each stage the tensor being acted on has a vacuum in the
new position, and $a_{-1}\one=a$ in that position.
Starting with $a_n^{(n)}\in P_1$, the successive products in
\eqref{eq:puretensorP} lie in $P_2,\ldots,P_n$.
Pure tensors span $X$, whereas every $P_j\subset X$, so
$P_n=X$.
Applying Lemma~\ref{lem:imageclosure} to $Z$ and this chain gives
\[
\dim X/C_1^{U^{(n)}}(X)<\infty.
\]
By Lemma~\ref{lem:tensorc2}, $X$ is $C_2$-cofinite when $V$
is. Since $U^{(n)}=X^{S_n}$, Theorem~\ref{thm:descent}
applies to $(X,S_n)$ and proves the last assertion.
\end{proof}

\begin{proof}[Proof of Theorem~\ref{thm:criterion}]
Let $V$ have central charge $c$ and satisfy the hypotheses of
the theorem. The Virasoro construction in
Section~\ref{subsec:virasoro-heisenberg} gives a graded homomorphism
\[
f:\mathcal V_c\longrightarrow V,\qquad \omega_c\longmapsto\omega, \qquad fD=Df.
\]
By functoriality \eqref{eq:functoriality}, the assumed relation gives
\[
D^r\omega=f(D^r\omega_c)
\in f(T_n(\mathcal V_c))\subseteq T_n(V).
\]
The inclusions from Lemma~\ref{lem:Tstability} give
\[
T_n(V)\subseteq T_{n-1}(V)\subseteq\cdots\subseteq T_2(V).
\]
Hence the same relation holds in every $T_j(V)$ with
$2\le j\le n$. Proposition~\ref{prop:translation} supplies
an integer $N_j\ge1$ such that $D^{N_j}V\subseteq T_j(V)$.
Start with $U^{(1)}=V$, which is $C_2$-cofinite by hypothesis.
If $U^{(j-1)}$ is $C_2$-cofinite, apply
Theorem~\ref{thm:seedfinite} with $(n,N)=(j,N_j)$ to
obtain the cofiniteness of $U^{(j)}$. Finite induction proves
this for all $1\le j\le n$. Applying
Theorem~\ref{prop:permutationsubgroups} in each degree
gives $C_2$-cofiniteness of $T^j(V)^G$ for every
$G\le S_j$, as asserted.
\end{proof}

\begin{corollary}\label{cor:criterionall}
Fix $c\in\C$. Suppose that, for every $n\ge2$, there exists
an integer $r=r(c,n)\ge0$ such that
\[
D^r\omega_c\in T_n(\mathcal V_c).
\]
Then $T^n(V)^G$ is $C_2$-cofinite for every $n\ge1$,
every subgroup $G\le S_n$, and every simple,
self contragredient, $C_2$-cofinite vertex operator algebra $V$ of CFT type
and central charge $c$.
\end{corollary}
\begin{proof}
For $n=1$ the assertion is the assumed cofiniteness of $V$.
For $n\ge2$, apply Theorem~\ref{thm:criterion} using the
relation for that $n$. The exponent need not be uniform as
$n$ varies.
\end{proof}

\Needspace{20\baselineskip}
\section{Corollary~\ref*{ex:fivefold}: The \texorpdfstring{$C_2$}{C2}-cofiniteness of \texorpdfstring{$T^5(V)^{S_5}$}{T5(V) invariant under S5}}
\label{sec:t5}

We now specialize to tensor degree five and prove that the
Virasoro criterion holds at every central charge.
The new input is the relation
$D^{16}\omega_c\in T_5(\mathcal V_c)$, obtained from a finite
integer calculation in Appendix~\ref{app:t5-verification}.
The symmetric operations, their functoriality,
and the quotient identifications are those already established in
Section~\ref{sec:conformalcriterion}. In particular, we will apply
Theorem~\ref{thm:criterion} directly instead of repeating the
general descent argument. We give the algebraic product and
specialization calculations explicitly below and invoke the
spanning certificate from the appendix.

Section~\ref{subsec:t5-products} develops the integral product recursion,
Section~\ref{subsec:t5-relation} specializes the polynomial relation,
and Section~\ref{subsec:t5-cofiniteness} applies the criterion.

\subsection{Virasoro products and polynomial specialization}
\label{subsec:specialization}
\label{subsec:t5-products}

The PBW construction in Section~\ref{subsec:virasoro-heisenberg} can be
performed over $R=\mathbb Q[c]$. Its weight spaces and
specializations are
\begin{equation*}
(\mathcal V_R)_d=
\bigoplus_{\substack{\lambda_1\ge\cdots\ge\lambda_k\ge2\\
                     |\lambda|=d}}Rv_\lambda,
\qquad
\mathcal V_R\otimes_{R,\,c\mapsto\alpha}\C=\mathcal V_\alpha,
\end{equation*}
where $v_\lambda=L(-\lambda_1)\cdots L(-\lambda_k)\one$.
The empty partition represents the vacuum. Put $t=c/2$ and
make the base change $c\mapsto2t$. The Virasoro bracket becomes
\begin{equation}\label{eq:t5-integralbracket}
[L(m),L(n)]=(m-n)L(m+n)
+\frac{m^3-m}{6}t\,\delta_{m+n,0}.
\end{equation}
The coefficient $(m^3-m)/6$ is an integer, since
$m^3-m=(m-1)m(m+1)$ is a product of three consecutive integers.
Thus reordering Virasoro operators introduces only coefficients
in $\mathbb Z[t]$.

More explicitly, if $v_p=L(-q)v_\mu$ with $q\ge2$, then
\begin{equation}
\label{eq:algorithm-L}
L(m)L(-q)v_\mu ={}L(-q)L(m)v_\mu+(m+q)L(m-q)v_\mu +\delta_{m,q}\frac{m^3-m}{6}t\,v_\mu.
\end{equation}
This follows by writing $L(m)L(-q)=[L(m),L(-q)]+L(-q)L(m)$
and substituting $n=-q$ in \eqref{eq:t5-integralbracket}.
The vacuum values are
\[
L(m)\one=
\begin{cases}
v_{(-m)},&m\le-2,\\
0,&m\ge-1.
\end{cases}
\]
Translation is equally explicit. Since $D\one=0$ and
$[D,L(-p)]=(p-1)L(-p-1)$, induction on $r$ gives
\begin{equation}\label{eq:t5-translation-factorial}
D^rL(-p)\one
=\frac{(p+r-2)!}{(p-2)!}L(-p-r)\one
\qquad(p\ge2,\ r\ge0).
\end{equation}
Indeed, applying $D$ to the expression for $r$ multiplies its
coefficient by $p+r-1$, which produces the expression for $r+1$.
For $p=2$, this yields
\[
D^r\omega=r!L(-r-2)\one.
\]

We next express the coefficients of a product of arbitrary PBW
vectors using these Virasoro actions. Write $v_a=L(-p)u$,
where $p\ge2$, and put $s=D^{p-2}\omega/(p-2)!$.
Differentiating the vertex operator expansion gives
\[
Y(D^r\omega,z) =\partial_z^r\sum_{k\in\Z}\omega_kz^{-k-1},\qquad (D^r\omega)_n =(-1)^r n(n-1)\cdots(n-r+1)\omega_{n-r}.
\]
Using $\omega_k=L(k-1)$, we obtain
\begin{equation}\label{eq:tailderivativecoeff}
s_n=(-1)^{p-2}\binom n{p-2}L(n-p+1),
\qquad s_{-1}u=L(-p)u.
\end{equation}
For the second equality, use
$\binom{-1}{p-2}=(-1)^{p-2}$.
The product identity \eqref{eq:iterate} therefore reads
\begin{equation}\label{eq:tailiterate}
(s_{-1}u)_n b
=\sum_{j\ge0}s_{-1-j}(u_{n+j}b)
 +\sum_{j\ge0}u_{n-1-j}(s_jb).
\end{equation}
Substitution of \eqref{eq:tailderivativecoeff}, together with
$(-1)^{p-2}\binom{-1-j}{p-2}=\binom{p+j-2}{p-2}$, gives
\begin{equation}
\label{eq:t5-explicit-product}
\begin{alignedat}{1}
(L(-p)u)_n b&=\sum_{j\ge0}\binom{p+j-2}{p-2} L(-p-j)(u_{n+j}b)\\
&\mathrel{\phantom{=}}+\sum_{j\ge0}(-1)^{p-2}\binom j{p-2} u_{n-1-j}\bigl(L(j-p+1)b\bigr).
\end{alignedat}
\end{equation}
Both sums are finite. If $u$ and $b$ have weights $a$ and $b_0$,
then their inner vectors in the two sums have respective weights
\[
a+b_0-n-j-1,\qquad b_0+p-j-1.
\]
Nonnegative weights therefore allow the upper bounds
$j\le a+b_0-n-1$ and $j\le b_0+p-1$, respectively.
A negative upper bound gives an empty sum. The initial case is
$\one_n b=\delta_{n,-1}b$; every occurrence of a coefficient
of $u$ has a first argument of smaller PBW length than $L(-p)u$.
Thus \eqref{eq:t5-explicit-product} determines each product by
finite induction. In particular, the normally ordered product is
obtained by setting $n=-1$.

The same finite identities compute every normally ordered product
occurring in the defining recursion \eqref{eq:B}.

For bookkeeping of polynomial degrees, put
\begin{equation*}
F_d=\Span_{\mathbb Q}\{t^ev_\lambda:
e\ge0,\ 2e+\ell(\lambda)\le d\},
\qquad F_{-1}=0.
\end{equation*}
The parameter $t$ has filtration degree two and conformal weight
zero; a Virasoro operator has filtration degree one. Thus
$\wt(t^ev_\lambda)=|\lambda|$.

\begin{lemma}\label{lem:t5-virasorofiltration}
All products of integral PBW vectors have coefficients in
$\mathbb Z[t]$, and
\[
(F_d)_mF_e\subseteq F_{d+e}\qquad(m\in\mathbb Z).
\]
In particular, if
$d=\sum_{i=1}^5\ell(\lambda^{(i)})+\ell(\mu)$, then
\begin{equation}\label{eq:t5-Bfiltration}
B_5(v_{\lambda^{(1)}},\ldots,v_{\lambda^{(5)}};v_\mu)
\in F_d.
\end{equation}
\end{lemma}
\begin{proof}
In \eqref{eq:t5-integralbracket}, replacing a pair of operators
by $(m-n)L(m+n)$ lowers filtration degree from two to one;
replacing it by a scalar multiple of $t$ preserves degree two.
Therefore reordering never increases the degree of an operator
word. Acting on the vacuum removes every remaining operator
$L(m)$ with $m\ge-1$, so $L(m)F_e\subseteq F_{e+1}$.

We prove by induction on $k=\ell(\lambda)$ that every
coefficient of $Y(v_\lambda,z)$ maps $F_e$ into $F_{e+k}$
and preserves the integral span of the vectors $t^av_\mu$.
For $k=0$, the vertex operator is the identity. For $k>0$, write
$v_\lambda=L(-p)u$ with $\ell(u)=k-1$ and use
\eqref{eq:t5-explicit-product}. In its first sum, the
coefficient of $u$ increases filtration degree by at most $k-1$,
and the following Virasoro action increases it by at most one.
In its second sum these operations occur in the reverse order,
with the same total bound $k$. All scalar factors are integer
binomial coefficients, and the Virasoro bracket has coefficients
in $\mathbb Z[t]$. This proves the induction statement.

Multiplication of the two input vectors by $t^a$ and $t^b$
multiplies their product by $t^{a+b}$. Including these degrees
gives $(F_d)_mF_e\subseteq F_{d+e}$ for general inputs.
Finally, each term of the recursion \eqref{eq:B} is an iterated
$-1$ product in which every input occurs once. Applying the
product bound at each step gives \eqref{eq:t5-Bfiltration}.
\end{proof}

These formulas also justify specialization: each product has
finitely many PBW coefficients in $\mathbb Q[t]$, and evaluation
$t\mapsto c/2$ commutes term by term with
\eqref{eq:t5-integralbracket}, \eqref{eq:t5-explicit-product},
translation, and the operations $B_k$.

\subsection{The polynomial relation and its specialization}
\label{subsec:t5-relation}

The finite calculation in Appendix~\ref{app:t5-verification}
establishes the spanning statement of
Lemma~\ref{lem:t5-virasorocertificate}.
We now explain how this implies a relation for every central charge.

\begin{proposition}\label{prop:t5-virasoro}
For every $c\in\mathbb C$,
\[
D^{16}\omega_c\in T_5(\mathcal V_c).
\]
The relation can be written as a finite $\mathbb Q[c]$-linear
combination of $B_5$ values with PBW arguments of total weight
eighteen and total PBW length at most nine. The polynomial
coefficients can be chosen to have degree at most ten.
Consequently $D^N\omega_c\in T_5(\mathcal V_c)$ for every $N\ge16$.
\end{proposition}
\begin{proof}
In the notation of Appendix~\ref{app:t5-verification},
\eqref{eq:t5-weight18combination} supplies rational constants
$\alpha_{q,j}$ with
\[
L(-18)\one=\sum_{q,j}\alpha_{q,j}t^jb_q(t).
\]
Group the terms with the same tuple $q$ and define
$p_q(t)=\sum_j\alpha_{q,j}t^j\in\mathbb Q[t]$.
The range of exponents in \eqref{eq:t5-generators} gives
$\deg p_q\le10$. Thus
\begin{equation}\label{eq:t5-polynomialrelation}
L(-18)\one=
\sum_{q\in\mathcal Q_{18}}p_q(t)
B_5(v_{\lambda^{(1)}},\ldots,v_{\lambda^{(5)}};v_\mu).
\end{equation}
Here $q=(\lambda^{(1)},\ldots,\lambda^{(5)};\mu)$ has
nonempty first five partitions, total weight eighteen, and
total length at most nine. Every first argument is therefore
of positive weight.

For a fixed $c\in\C$, let $v_\lambda^{(c)}$ be the corresponding
PBW vector in $\mathcal V_c$. Specialization is the map
\[
\operatorname{ev}_c(t^ev_\lambda)=(c/2)^ev_\lambda^{(c)}.
\]
The formulas of Section~\ref{subsec:t5-products} show that
this map preserves the vertex products and the operations $B_5$.
Applying it to \eqref{eq:t5-polynomialrelation} gives
\begin{equation}\label{eq:t5-specializedrelation}
L(-18)\one=
\sum_{q\in\mathcal Q_{18}}p_q(c/2)
B_5(v_{\lambda^{(1)}}^{(c)},\ldots,
    v_{\lambda^{(5)}}^{(c)};v_\mu^{(c)})
\in T_5(\mathcal V_c).
\end{equation}
If $p_q(t)=\sum_j\alpha_{q,j}t^j$, then
$p_q(c/2)=\sum_j2^{-j}\alpha_{q,j}c^j$.
These are polynomials in $c$ with rational coefficients;
their denominators are constant integers. Hence specialization
is valid for every $c$, with no excluded values.

Taking $p=2$ and $r=16$ in
\eqref{eq:t5-translation-factorial} gives
\[
D^{16}\omega_c=16!L(-18)\one.
\]
Multiply \eqref{eq:t5-specializedrelation} by $16!$ to obtain
the asserted $\mathbb Q[c]$-linear relation. Finally, for
$N\ge16$, Lemma~\ref{lem:Tstability} gives
\[
D^N\omega_c=D^{N-16}(D^{16}\omega_c)
\in D^{N-16}T_5(\mathcal V_c)\subseteq T_5(\mathcal V_c).\qedhere
\]
\end{proof}

\subsection{Proof of Corollary~\ref*{ex:fivefold}}
\label{subsec:t5-cofiniteness}

\begin{theorem}\label{thm:t5-main}
Let $V$ be a simple, self contragredient, $C_2$-cofinite vertex operator algebra
of CFT type. Then
\begin{equation*}
\dim_{\mathbb C}\bigl(U^{(5)}/C_2(U^{(5)})\bigr)
=\dim_{\mathbb C}\bigl(T^5(V)^{S_5}/C_2(T^5(V)^{S_5})\bigr)
<\infty.
\end{equation*}
More generally, $T^j(V)^G$ is $C_2$-cofinite for every
$1\le j\le5$ and every subgroup $G\le S_j$.
\end{theorem}
\begin{proof}
The construction in Section~\ref{subsec:virasoro-heisenberg}
gives a graded vertex algebra homomorphism
\[
f:\mathcal V_{c_V}\longrightarrow V,\qquad \omega\longmapsto\omega_V, \qquad f(\one)=\one,\quad fD=Df.
\]
Proposition~\ref{prop:t5-virasoro} and functoriality
\eqref{eq:functoriality} therefore imply
\begin{equation*}\label{eq:t5-transfer}
D^{16}\omega_V=f(D^{16}\omega)
\in f(T_5(\mathcal V_{c_V}))\subseteq T_5(V).
\end{equation*}
Proposition~\ref{prop:t5-virasoro} verifies the hypothesis of Theorem~\ref{thm:criterion} at $c=c_V$, with
$n=5$ and $r=16$. That theorem gives $C_2$-cofiniteness of
$T^j(V)^G$ for every $1\le j\le5$ and every $G\le S_j$.
\end{proof}

\begin{corollary}\label{cor:t5-rationality}
If $V$ is strongly rational, then $T^5(V)^G$ is strongly
rational for every subgroup $G\le S_5$.
\end{corollary}
\begin{proof}
The tensor product $T^5(V)$ is strongly rational.
Theorem~\ref{thm:t5-main} gives $C_2$-cofiniteness of
$T^5(V)^G$. Apply Theorem~\ref{thm:external-orbifold-rationality}
to the finite image of $G$ in $\operatorname{Aut}(T^5(V))$
to obtain the assertion.
\end{proof}

\appendix

\Needspace{7\baselineskip}
\section{Exact verification of the fivefold relation}
\label{app:t5-verification}

\subsection{Finite coefficient space}
\label{subsec:t5-certificate}

The Virasoro product formulas and filtration are established in
Section~\ref{subsec:t5-products}. This appendix constructs
the finite coefficient space, gives its exact spanning certificate,
and records the implementation and verification program.

Let $\mathcal Q_{18}$ consist of tuples
$q=(\lambda^{(1)},\ldots,\lambda^{(5)};\mu)$ with
nonempty first five partitions, in nondecreasing lexicographic
order, and
\[
\sum_{i=1}^5|\lambda^{(i)}|+|\mu|=18.
\]
The last partition may be empty. Write
\[
d(q)=\sum_{i=1}^5\ell(\lambda^{(i)})+\ell(\mu),
\qquad
b_q(t)=B_5(v_{\lambda^{(1)}},\ldots,v_{\lambda^{(5)}};v_\mu).
\]
Necessarily $5\le d(q)\le9$. Consider the rational vector
space and the finite indexed collection
\begin{align}
E&=\Span_{\mathbb Q}\{t^ev_\lambda:
|\lambda|=18,\ 2e+\ell(\lambda)\le26\},
\label{eq:t5-coefficientspace}\\
\mathcal G&=\bigl(t^jb_q(t):
q\in\mathcal Q_{18},\ 0\le j\le
\lfloor(26-d(q))/2\rfloor\bigr).
\label{eq:t5-generators}
\end{align}
The powers of $t$ in \eqref{eq:t5-coefficientspace} are
separate basis vectors over $\mathbb Q$. No rational function
in $t$ is inverted.

\begin{lemma}\label{lem:t5-enumeration}
Every vector in $\mathcal G$ belongs to $E$, and
\[
\dim_{\mathbb Q}E=9\cdot88+87+80+51+9=1019,
\qquad |\mathcal Q_{18}|=249,\qquad |\mathcal G|=2569.
\]
\end{lemma}
\begin{proof}
The containment follows from \eqref{eq:t5-Bfiltration}.
Let $p_{\le r}(w)$ count partitions of $w$ into at most
$r$ parts, each at least two. These numbers are given by
\[
p_{\le r}(w)
=[z^w]\sum_{\ell=0}^r
\frac{z^{2\ell}}{\prod_{i=1}^{\ell}(1-z^i)}.
\]
For weight eighteen, every partition has length at most nine,
and
\[
p_{\le9}(18)=88,\quad p_{\le8}(18)=87,\quad
p_{\le6}(18)=80,\quad p_{\le4}(18)=51,\quad
p_{\le2}(18)=9.
\]
The exponents $e=0,\ldots,8$ in
\eqref{eq:t5-coefficientspace} each permit all $88$ partitions;
$e=9,10,11,12$ give the remaining four counts. This proves
the dimension formula.

The input enumeration, followed by the allowed powers of $t$,
is
\[
\begin{array}{c|r|c|r}
d(q)&\text{PBW tuples}&\text{allowed }j&\text{rows}\\ \hline
5&18&0,\ldots,10&198\\
6&73&0,\ldots,10&803\\
7&95&0,\ldots,9&950\\
8&51&0,\ldots,9&510\\
9&12&0,\ldots,8&108\\ \hline
&249&&2569
\end{array}
\]
For example, the number of tuples of length $d$ is the
coefficient of $y^5z^{18}x^d$ in
\[
\prod_{\lambda\ne\varnothing}
\frac{1}{1-yz^{|\lambda|}x^{\ell(\lambda)}}
\prod_{p\ge2}\frac{1}{1-xz^p},
\]
where $\lambda$ ranges over partitions with parts at least
two. The first product selects the unordered first five
arguments, and the second selects $\mu$.
All products may be truncated at weight eighteen.
This gives the displayed counts; multiplying by the number
of allowed $j$ gives the row counts.
\end{proof}

\subsection{Exact coefficient representation}
\label{subsec:t5-verification}

Put $t=c/2$. A vector is represented by a finite coefficient array
\[
(e,\lambda)\longmapsto a_{e,\lambda}\in\mathbb Z,
\qquad
v=\sum_{e,\lambda}a_{e,\lambda}t^e
 L(-\lambda_1)\cdots L(-\lambda_r)\one,
\]
where $\lambda_1\ge\cdots\ge\lambda_r\ge2$.
The empty partition represents the vacuum.
Addition is coefficientwise, and multiplication by $t^j$
replaces $(e,\lambda)$ by $(e+j,\lambda)$.
Only nonzero coefficients need to be retained.
In particular, $t$ is an independent polynomial variable, not
a value substituted for the central charge.

The action on a PBW vector is computed from
\eqref{eq:algorithm-L}, and general products from
\eqref{eq:t5-explicit-product} with its stated weight bounds.
The vacuum coefficients are $\one_n b=\delta_{n,-1}b$.
Each product reduces the PBW length of its first argument;
commuting Virasoro operators uses a shorter PBW tail, with reordering of operators with negative indices reducing inversions. These reductions
terminate at the vacuum values of
Section~\ref{subsec:t5-products}.

The operations $B_k$ are evaluated by \eqref{eq:B}, starting
with $B_0(;b)=b$. Each call reduces the number of first arguments.
Lemma~\ref{lem:Bsymmetry} permits these arguments to be sorted
before a value is cached.

\subsection{Enumerating the finite matrix}
\label{subsec:t5-matrix}

Partitions of weight $w$ with parts at least two are
enumerated recursively. Given an upper bound $p_{\max}$ on
the next part, choose $p=\min(w,p_{\max}),\ldots,2$ in that order, replace
the remaining weight by $w-p$, and bound the next part by $p$.
A prescribed length bound is reduced by one at each step.
The empty partition occurs exactly when the remaining weight
is zero. Thus every admissible partition is generated once,
in nonincreasing order.

For the first five arguments of $B_5$, choose a nondecreasing
list of nonempty partitions. Repetition is allowed, as required
by symmetry. If $h$ arguments remain, their total remaining
weight must be at least $2h$ and their total remaining PBW
length at least $h$. These two bounds restrict the enumeration.
After the first five arguments are selected, the last argument
is any partition of the remaining weight and may be empty.
The resulting tuples are exactly $\mathcal Q_{18}$ in
Section~\ref{subsec:t5-certificate}.

For each tuple $q$, compute $b_q(t)$ and its input length
$d(q)$. The filtration estimate \eqref{eq:t5-Bfiltration}
implies
\[
b_q(t)=\sum_{e,\lambda}a_{e,\lambda}t^ev_\lambda,
\qquad |\lambda|=18,\quad 2e+\ell(\lambda)\le d(q).
\]
For $0\le j\le\lfloor(26-d(q))/2\rfloor$, shifting the
exponent by $j$ gives a row representing $t^jb_q(t)$ in $E$.
The inequality $2(e+j)+\ell(\lambda)\le26$ ensures that
no coefficient falls outside the chosen column basis.
Columns are ordered by the tuple
\[
(2e+\ell(\lambda),\,\ell(\lambda),\,\lambda,\,e)
\]
in decreasing lexicographic order. The enumeration in
Lemma~\ref{lem:t5-enumeration} gives $1019$ columns and
$2569$ rows. Rows are grouped by increasing input length
$d(q)=5,\ldots,9$. Within each group, the first five
partitions are enumerated lexicographically, the last partition
uses the recursion above for partitions, and the allowed
exponents $j$ are taken in increasing order for each tuple.
These orders, together with the recorded zero based row
indices, determine the certified square minor uniquely.

\subsection{Modular elimination and the spanning certificate}
\label{subsec:t5-elimination}

Only after every polynomial coefficient has been constructed
over $\mathbb Z$ do we reduce the integer coordinate
matrix $M$ modulo $p=101$. Process the rows in their fixed order. A nonzero working row $v$ has leading index
$j=\min\{k:v_k\ne0\}$. If a pivot row $u^{(j)}$ with
$u^{(j)}_j=1$ is already stored, replace
\[
v\longleftarrow v-v_j u^{(j)}\pmod{101}.
\]
This kills its leading entry. If no such pivot is stored,
normalize the row by
\[
u^{(j)}\longleftarrow v_j^{-1}v\pmod{101}
\]
and record the index of the original input row.
Distinct pivot positions prove independence of the recorded
rows. Here all $1019$ column positions become pivots, so the
selected index set $I$ has size $1019$.

Independently recompute the determinant of the original
square submatrix $M_I$ by row elimination, processing columns
in order. If the successive nonzero pivots are
$z_1,\ldots,z_{1019}$ and there are $s$ row exchanges, then
\[
\det M_I\equiv(-1)^s\prod_{j=1}^{1019}z_j\pmod{101}.
\]
Row replacement by a multiple of another row does not change
the determinant; a row exchange changes its sign.
The resulting residue is $31$, which is nonzero.
The following statement records the resulting conclusion in characteristic zero.

\begin{lemma}\label{lem:t5-virasorocertificate}
The vectors in $\mathcal G$ span $E$ over $\mathbb Q$.
\end{lemma}
\begin{proof}
Let $M$ be their $2569\times1019$ integer coordinate matrix
in the basis \eqref{eq:t5-coefficientspace}.
Use the ordering in Section~\ref{subsec:t5-matrix}, and let
$I$ be the row indices selected by the elimination just described
over $\mathbb F_{101}$. The exact calculation, implemented in
Section~\ref{subsec:t5-program}, gives
\begin{equation*}
|I|=1019,\qquad
\operatorname{rank}_{\mathbb F_{101}}M=1019,\qquad
\det M_I\equiv31\pmod{101}.
\end{equation*}
The program constructs the coefficients over $\mathbb Z[t]$
and verifies the displayed rank and determinant by elimination over a finite field. Since reduction modulo a prime cannot increase
the rank of an integer matrix,
\[
1019=\operatorname{rank}_{\mathbb F_{101}}M
\le\operatorname{rank}_{\mathbb Q}M\le\dim_{\mathbb Q}E=1019.
\]
Thus the rows span $E$, and in particular
\begin{equation}\label{eq:t5-weight18combination}
L(-18)\one\in E=\Span_{\mathbb Q}\mathcal G.
\end{equation}
The coefficients of such an expression in the vectors
$t^jb_q(t)$ are rational constants, so their denominators are
independent of $t$. The rank certificate proves the required
membership without computing these coefficients or a matrix
inverse.
\end{proof}

Reduction modulo $101$ certifies an integer determinant; it
never substitutes a value for $c$. Specialization of the
resulting polynomial relation is justified in
Proposition~\ref{prop:t5-virasoro}.

\subsection{Implementation notes}
\label{subsec:t5-implementation}

The program implements the integer calculation of
Sections~\ref{subsec:t5-certificate}--\ref{subsec:t5-elimination}.
A vector is stored as a dictionary with keys \texttt{(e,p)},
where \texttt{p} is a PBW partition tuple, and integer values
giving the coefficients of $t^ev_p$. Addition combines equal
keys, multiplication by a power of $t$ shifts \texttt{e},
and entries with zero coefficient are deleted.

The routine \texttt{L(m,p)} implements \eqref{eq:algorithm-L}.
The routine \texttt{vertexcoefficient(a,n,b)} uses
\eqref{eq:tailderivativecoeff} and \eqref{eq:tailiterate}, with
the weight bounds derived in Section~\ref{subsec:t5-products}.
Normal products correspond to \texttt{n=-1}.
The routine \texttt{B(args,b)} implements \eqref{eq:B}.
Symmetry permits its first arguments to be sorted for caching.

The function \texttt{partitions(w,maxpart,maxlen)} enumerates
the partitions with the specified weight, part and length bounds.
The routine \texttt{generators} applies the bounds on remaining
weight and length and selects the nondecreasing first five
partitions. Before assembling a matrix row, the code checks
that every dictionary key belongs to the prescribed column basis.
The column order is the one fixed in
Section~\ref{subsec:t5-matrix}.

The routine \texttt{choose\_independent\_rows} selects the pivot
rows modulo $101$ and records their original indices.
The routine \texttt{determinant\_mod} independently checks the
determinant of the selected square submatrix.
With \texttt{--audit}, the program also computes all $1931$
normal products by an independent recursion on the right PBW
argument and compares their full integer dictionaries with those
obtained from \eqref{eq:tailiterate}.
These comparisons take place before modular reduction.

\subsection{Verification program}
\label{subsec:t5-program}

The following Python program uses only the standard library and
implements the exact verification described in
Sections~\ref{subsec:t5-certificate}--\ref{subsec:t5-elimination}. Save it as
\texttt{verify\_t5\_d16.py} and run
\begin{verbatim}
python3 verify_t5_d16.py --audit
\end{verbatim}
The output file \texttt{t5\_d16\_certificate.json} records the
selected row indices, their PBW labels, and the column labels.
No external coefficient data are needed. The program was run with
Python 3.12.14.

\begingroup
\scriptsize
\begin{verbatim}
"""Exact polynomial certificate for the weight-eighteen relation.
Run: python3 verify_t5_d16.py --audit
Only the Python standard library is used; t=c/2 stays formal.
"""
from functools import lru_cache
from math import comb

STAR_PAIRS=set()


def addto(out, vec, factor=1, shift=0):
    if not factor:
        return
    for (e, p), v in vec.items():
        key = e + shift, p
        out[key] = out.get(key, 0) + factor * v
        if not out[key]:
            del out[key]


def binom(n, k):
    if k < 0:
        return 0
    if n < 0:
        return (-1)**k * comb(k-n-1, k)
    return comb(n, k) if k <= n else 0


@lru_cache(None)
def L(m, p):
    """L_m acting on PBW vector p."""
    if not p:
        return {(0, (-m,)): 1} if m <= -2 else {}
    if m <= -2 and -m >= p[0]:
        return {(0, (-m,) + p): 1}
    q, rest = p[0], p[1:]
    out = {}
    for (e, s), v in L(m, rest).items():
        addto(out, L(-q, s), v, e)
    addto(out, L(m-q, rest), m+q)
    if m == q:
        addto(out, {(0, rest): 1}, (m**3-m)//6, 1)
    return out


def Lvec(m, v):
    out = {}
    for (e, p), z in v.items():
        addto(out, L(m, p), z, e)
    return out


@lru_cache(None)
def vertexcoefficient(a, n, b):
    """a_n b, from the normal-ordering recursion in the manuscript."""
    if not a:
        return {(0, b): 1} if n == -1 else {}
    if sum(a) + sum(b) - n - 1 < 0:
        return {}
    p, u = a[0], a[1:]
    if not u:
        coefficient = (-1)**(p-2) * binom(n, p-2)
        return {key: coefficient*v for key, v in L(n-p+1, b).items()} if coefficient else {}
    # a = (D^(p-2) omega/(p-2)!)_{-1} u.
    out = {}
    for j in range(max(0, sum(u)+sum(b)-n)):
        tmp = vertexcoefficient(u, n+j, b)
        if tmp:
            addto(out, Lvec(-p-j, tmp), binom(p+j-2, p-2))
    for j in range(max(0, sum(b)+p)):
        coefficient = (-1)**(p-2) * binom(j, p-2)
        if coefficient:
            for (e, s), v in L(j-p+1, b).items():
                addto(out, vertexcoefficient(u, n-1-j, s), coefficient*v, e)
    return out


@lru_cache(None)
def B(args, b):
    """B_k on basis arguments. First arguments are a sorted tuple."""
    if not args:
        return {(0, b): 1}
    a, rest = args[0], args[1:]
    out = {}
    STAR_PAIRS.add((a,b))
    for (e, p), v in vertexcoefficient(a, -1, b).items():
        addto(out, B(rest, p), -v, e)
    for i, x in enumerate(rest):
        STAR_PAIRS.add((a,x))
        for (e, p), v in vertexcoefficient(a, -1, x).items():
            newargs = tuple(sorted(rest[:i] + (p,) + rest[i+1:]))
            addto(out, B(newargs, b), -v, e)
    return out


@lru_cache(None)
def partitions(w, maxpart=None, maxlen=None):
    if w == 0:
        return ((),)
    if w < 2 or maxlen == 0:
        return ()
    maxpart = min(w, maxpart if maxpart is not None else w)
    ans = []
    for p in range(maxpart, 1, -1):
        for tail in partitions(w-p, p, None if maxlen is None else maxlen-1):
            ans.append((p,)+tail)
    return tuple(ans)


def generators(w, k, degree):
    """All c-independent PBW B_k inputs with weight w, length<=degree."""
    options = [(sum(p), len(p), p) for r in range(2, w-2*(k-1)+1)
               for p in partitions(r, maxlen=degree-k+1)]
    options.sort(key=lambda z: z[2])
    def firsts(args, start, wt, lng, left):
        if left == 0:
            for b in partitions(wt, maxlen=lng):
                yield args, b
            return
        for i in range(start, len(options)):
            wp, lp, p = options[i]
            if wp <= wt-2*(left-1) and lp <= lng-left+1:
                yield from firsts(args+(p,), i, wt-wp, lng-lp, left-1)
    yield from firsts((), 0, w, degree, k)


def column_order(key):
    e,p=key
    return 2*e+len(p),len(p),p,e


def choose_independent_rows(matrix, prime):
    """Sparse forward elimination; retain indices of independent input rows."""
    pivots={}
    chosen=[]
    for i,row in enumerate(matrix):
        v={j:a%prime for j,a in row.items() if a%prime}
        while v:
            lead=min(v)
            if lead not in pivots:
                inv=pow(v[lead],-1,prime)
                pivots[lead]={j:a*inv%prime for j,a in v.items()}
                chosen.append(i)
                break
            z=v[lead]
            for j,a in pivots[lead].items():
                q=(v.get(j,0)-z*a)%prime
                if q:v[j]=q
                else:v.pop(j,None)
    return chosen


def determinant_mod(square, prime):
    """Independent, column-driven elimination of a square integer minor."""
    a=[{j:z%prime for j,z in row.items() if z%prime}
       for row in square]
    det=1
    for j in range(len(a)):
        pivot=next((i for i in range(j,len(a)) if a[i].get(j,0)),None)
        if pivot is None:return 0
        if pivot!=j:
            a[j],a[pivot]=a[pivot],a[j]
            det=-det
        z=a[j][j]
        det=det*z%prime
        inv=pow(z,-1,prime)
        top=[(k,v) for k,v in a[j].items() if k>j]
        for i in range(j+1,len(a)):
            scale=a[i].pop(j,0)*inv%prime
            if scale:
                for k,v in top:
                    q=(a[i].get(k,0)-scale*v)%prime
                    if q:a[i][k]=q
                    else:a[i].pop(k,None)
    return det%prime


def build_certificate():
    import json
    from pathlib import Path
    weight,degree,k,prime=18,26,5,101
    columns=sorted([(e,p) for e in range(degree//2+1)
                    for p in partitions(weight,maxlen=degree-2*e)],
                   key=column_order,reverse=True)
    assert len(columns)==1019
    index={key:i for i,key in enumerate(columns)}
    matrix=[];labels=[];counts={}
    for d in range(k,weight//2+1):
        counts[d]=0
        for a,b in generators(weight,k,d):
            if sum(map(len,a))+len(b)!=d:continue
            counts[d]+=1
            value=B(a,b)
            assert all(sum(p)==weight and 2*e+len(p)<=d
                       for e,p in value)
            for shift in range((degree-d)//2+1):
                row={index[e+shift,p]:z for (e,p),z in value.items() if z}
                matrix.append(row)
                labels.append((shift,a,b))
    assert counts=={5:18,6:73,7:95,8:51,9:12}
    assert len(matrix)==2569
    chosen=choose_independent_rows(matrix,prime)
    assert len(chosen)==1019
    residue=determinant_mod([matrix[i] for i in chosen],prime)
    assert residue==31
    certificate={
        'claim':'D^16 omega belongs to T_5(V_c) for every complex c',
        'parameter':'t=c/2',
        'weight':weight,'filtration_bound':degree,
        'matrix_size':[len(matrix),len(columns)],
        'counts_by_input_PBW_length':counts,
        'prime':prime,'rank_mod_prime':len(chosen),
        'minor_determinant_mod_prime':residue,
        'selected_zero_based_row_indices':chosen,
        'selected_row_labels':[labels[i] for i in chosen],
        'columns':columns,
    }
    path=Path(__file__).with_name('t5_d16_certificate.json')
    path.write_text(json.dumps(certificate,indent=2)+'\n')
    print('Input PBW lengths:',counts,flush=True)
    print('Integer matrix: 2569 x 1019; rank modulo 101: 1019',
          flush=True)
    print('Selected minor determinant modulo 101:',residue,flush=True)
    print('PASS: D^16 omega=16! L(-18)1 belongs to T_5(V_c)'
          ' for every complex c.',flush=True)


@lru_cache(None)
def right_coefficient(a,n,b):
    """Independent recursion on the right PBW argument."""
    if not a:return {(0,b):1} if n==-1 else {}
    if sum(a)+sum(b)-n-1<0:return {}
    if not b:
        if n>=0:return {}
        out={(0,a):1}
        for i in range(1,-n):
            out=Lvec(-1,out)
            assert all(x%i==0 for x in out.values())
            out={key:x//i for key,x in out.items()}
        return out
    p,tail=b[0],b[1:]
    d=sum(a)
    out=Lvec(-p,right_coefficient(a,n,tail))
    addto(out,right_coefficient(a,n-p,tail),n+(p-1)*(d-1))
    for j in range(2,d+2):
        for (e,s),v in L(j-1,a).items():
            addto(out,right_coefficient(s,n-p+1-j,tail),-binom(1-p,j)*v,e)
    return out


def audit_products():
    for a,b in sorted(STAR_PAIRS):
        assert vertexcoefficient(a,-1,b)==right_coefficient(a,-1,b),(a,b)
    assert len(STAR_PAIRS)==1931
    print('AUDIT PASS: all 1931 polynomial normal products agree'
          ' with the independent right-argument recursion.')


if __name__=='__main__':
    import argparse
    parser=argparse.ArgumentParser(description=__doc__)
    parser.add_argument('--audit',action='store_true',
                        help='cross-check all normal products')
    args=parser.parse_args()
    build_certificate()
    if args.audit:audit_products()
\end{verbatim}
\endgroup

\section*{AI use disclosure}
We used ChatGPT during the development and drafting of this paper
to explore questions, possible formulations, and proof strategies;
to assist with computations and literature searches; and to review arguments and improve the exposition.
This assistance included suggestions for intermediate claims, proof steps, and corrections.
All mathematical statements, proofs, computational results, and references were independently checked by the
authors, who take sole responsibility for the contents of the paper.

\end{document}